\documentclass[12pt]{article}
\usepackage{latexsym,amsfonts}

\newtheorem{Defn}{Definition }[section]
\newtheorem{Propn}[Defn]{Proposition }
\newtheorem{Lemma}[Defn]{Lemma }
\newtheorem{Thm}[Defn]{Theorem }
\newtheorem{Exl}[Defn]{Example }
\newtheorem{Cor}[Defn]{Corollary }

\batchmode

\def\choose#1#2{\left(\!\ba{c}#1\\ #2\ea\!\right)}

\def\as#1{\renewcommand{\arraystretch}{#1}}

\def\dyl{\displaystyle }
\def\Cdot{ \mbox{\begin{picture}(8,6)\put(4,3){\circle*{1.5}} \end{picture}}  }

\def\defeq{\stackrel{\rm def}{=\!\!=}}

\def\ub#1{^{(#1)}}

\def\lrb#1{_{#1]}}

\def\bc{\begin{center}}                    
\def\ec{\end{center}}

\def\r{_{r}}

\def\iy{_{\infty}}
\def\iiy{^{\infty}}

\def\d{^{\dag}}
\def\otsym{\otimes_{\!\!\!\mbox{ \tiny\rm sym}}}
\def\CA{{\cal A}}                
\def\CC{{\cal C}}                
\def\CD{{\cal D}}
\def\CE{{\cal E}}                
\def\CF{{\cal F}}

\def\CH{{\cal H}}
\def\CK{{\cal K}}
\def\CN{{\cal N}}

\def\CS{{\cal S}}
\def\CW{{\cal W}}                

\def\CM{{\cal M}}               
\def\CP{{\cal P}}               

\def\CL{{\cal L}}
\def\CB{{\cal B}}

\def\CR{{\cal R}}

\def\la{\left\langle}
\def\ra{\right\rangle}

\def\FD{{\frak D}}

\def\FF{{\frak F}}

\def\FH{{\frak H}}
\def\Fh{{\frak h}}

\def\FL{{\frak L}}
\def\FM{{\frak M}}

\def\FU{{\frak U}}

\def\BBC{{\Bbb C}}

\def\BBE{{\Bbb E}}

\def\BBN{{\Bbb N}}
\def\BBP{{\Bbb P}}

\def\BBQ{{\Bbb Q}}
\def\BBR{{\Bbb R}}

\def\CD{{\cal D}}
\def\Ga{\alpha}
\def\Gb{\beta}
\def\GD{\Delta}
\def\Gd{\delta}
\def\Ge{\epsilon}
\def\Gve{\varepsilon}

\def\Gg{\gamma}                   
                
\def\LGa{\mbox{\Large$\alpha$}}   
\def\LGk{\mbox{\Large$\kappa$}}

\def\Lnu{\mbox{\Large$\nu$}}

\def\Gs{\sigma}

\def\Gv{\varphi}                    

\def\GL{\Lambda}

\def\GW{\Omega}                    
 
\def\Gth{\theta}
\def\step{^{step}}

\def\ub#1{^{(#1)}}

\def\beqn{\begin{eqnarray}}
\def\eeqn{\end{eqnarray}}
\def\bb{\bigbreak}
\def\ni{\noindent}

\def\pr{^{\prime}}
\def\prpr{^{\prime\prime}}

\def\st{^{{\textstyle *}}}

\def\ol{\overline}
\def\ul{\underline}

\def\row{\rightarrow}

\def\roy{\rightarrow\infty} 
\def\nm#1{\left\|#1\right\|}      

\def\ip#1#2{\left\langle#1,#2\right\rangle}

\def\pf{{\it Proof. }} 

\def\hraa{\lhook\joinrel\longrightarrow\!\!\!\!\!\rightarrow}

\def\inv{^{-1}}

\def\mx#1#2#3#4{\left[\begin{array}{cc} #1&#2\\#3&#4 \end{array}\right]}

\def\ba{\begin{array}}
\def\ea{\end{array}}
\def\epf{\rule{1.5mm}{2mm}}

\def\btb{\begin{thebibliography}}
\def\etb{\end{thebibliography}}

\begin{document}

\title{The Ito Formula for Essentially Self-adjoint Quantum 
Semimartingales} \author{G. F. Vincent--Smith} 
\date{ \ }
\maketitle 

\section {Introduction.} 
It is a source of embarrassment that the quantum Ito formula 
\beqn\label{Ito1}
f(\hat{M}_t)&=&f(\hat{M}_0)+\int_0^t(Df(\hat{M}_s)(d\hat{M}_s)+D^2_I 
f(\hat{M}_s)(d\hat{M}_s,d\hat{M}_s))
\eeqn
in \cite{Vin2}, for a regular (bounded) quantum semimartingale $M$, does not directly imply the Ito formula for classical 
Brownian motion $W$:
\beqn\label{Ito2}
f({W }_t)=f({W }_0)
+\int_0^t f\pr({W }_s)\,d{W }_s+\frac12 \int_0^t 
f\prpr({W }_s)\,d\la {W } \ra_s
\eeqn
This is unsatisfactory from our present point of 
view, which is to derive as much as possible 
of the classical theory of Brownian motion from the theory of quantum 
stochastic processes.
The quantum Ito formula belongs to non-commutative probability.
It is a purely operator theoretic result and makes no 
reference to paths. It implies classical results by regarding 
random variables as multiplication operators and using the 
Wiener--Ito isomorphism between Weiner space and Fock space.
\bb\ni
In this article we find a remedy and show that 
(\ref{Ito2}), and its extension to gauge processes, 
includes the case where $\hat{M}$ is an essentially 
self-adjoint quantum semimartingale.
\bb\ni
The classical Ito formula (\ref{Ito2}) then follows 
from the quantum Ito formula (\ref{Ito1}) as in 
\cite[Section 7.]{Vin2} by taking 
$\hat{M}=B\equiv A+A\d$, the Brownian quantum semimartingale,
which is a process of unbounded essentially self-adjoint operators.
\bb\ni
In (\ref{Ito1}) the process $f(\hat{M})=\{\hat{M}_t:t\in [0,1]\}$ 
is defined by the Fourier functional calculus, 
$Df(X)(\cdot)$ is the differential of $f$ and 
$D^2_I f(X)(\cdot,\cdot)$ is the (unsymmetrised) 
`Ito' second differential of $f$ at the operator $X$. 
\bb\ni
The Ito formula (\ref{Ito2}) may be obtained indirectly from 
(\ref{Ito1}), quite simply, as follows. 
It was shown in \cite{Vin2} that (\ref{Ito1}) implies 
(\ref{Ito2}) when ${W }$ is 
replaced by $h_n(W)$ where $h_n$ is bounded twice continuously 
differentiable. The Ito formula for Brownian motion then follows by 
approximating the function $h(x)=x$ by an appropriate sequence of 
$h_n$, for example the sequence used in Section \ref{PCM}. 
Lifting such an argument to non-commuting essentially self-adjoint 
quantum semimartingales presents difficulties. 
\def\so{_{\rm so}}
\bb\ni
Let $\hat{M}=\{\hat{M}_t:t\in [0,1]\}$ be a symmetric quantum semimartingale:
\beqn\label{E 71.202}
\hat{M}_t=\hat{M}_0+\int_0^t (\hat{E}\,d\GL+\hat{F}\,dA
+\hat{F}\st\,dA\d+\hat{H}\,ds).
\eeqn
The integrands are adapted processes of operators in $\CB(\FH)$, 
where $\FH$ is Bose--Fock space over $L^2[0,1]$,
and $E_t$ and $H_t$ are self-adjoint for all $t$.
\bb\ni
$\hat{M}$ is {\it essentially self-adjoint} with {\it core } 
$\CD$ if $\hat{M}_t$ is essentially self-adjoint on $\CD$ 
for almost all $t$. In this case the {\it closure } of $\hat{M}$ is
$\ol{M}=\{\ol{M}_t:t\in [0,1]\}$ where $\ol{M}_t$ is the closure 
of $\hat{M}_t$. The process of bounded operators 
$$
e^{i\ol{M}}=\{e^{i\ol{M}_t}:t\in[0,1]\}
$$
may then be defined by the functional calculus.
\bb\ni
An essentially self-adjoint quantum semimartingale $\hat{M}$ 
{\it satisfies the quantum Duhamel formula} if $e^{i\ol{M}}$ is 
a regular quantum semimartingale and 
\beqn\label{EE 72.01} 
e^{i\ol{M}_t}= I+\int_0^t (\hat{E}_{{\rm exp}(iM)}\,d\GL
+\hat{F}_{{\rm exp}(iM)}\,dA
+\hat{G}_{{\rm exp}(iM)}\,dA\d
+\hat{H}_{{\rm exp}(iM)}\,ds),
\eeqn
where 
\begin{eqnarray}
\nonumber\hat{E}_{{\rm exp}(iM)}(t)&
=& e^{i(\ol{M}_t+\hat{E}_t)}-e^{i\ol{M}_t}
\\
\nonumber\hat{F}_{{\rm exp}(iM)}(t)&=& i\int_0^1 e^{i(1-u)\ol{M}_t} 
\hat{F}_t e^{iu(\ol{M}_t+\hat{E}_t)}\,du
\\
\label{EE 79.999}\hat{G}_{{\rm exp}(iM)}(t)&=& i\int_0^1 e^{i(1-
u)(\ol{M}_t+\hat{E}_t)} \hat{F}\st_t e^{iu \ol{M}_t}\,du
\\
\nonumber\hat{H}_{{\rm exp}(iM)}(t)&=& i\int_0^1 e^{i(1-u)\ol{M}_t} 
\hat{H}_t e^{iu \ol{M}_t}\,du
\\
\nonumber &&\hspace{0.5in} + i^2\int_0^1\int_0^1
ue^{i(1-u)\ol{M}_t}F_te^{iu(1-v)(\ol{M}_t+\hat{E}_t)} 
\hat{F}\st_t e^{iuv \ol{M}_t}\,du\,dv
\end{eqnarray} 
Let $f\in C^{2+}(\BBR)=\{f\in L^1(\BBR):p\mapsto p^2 \hat{f}(p) 
\mbox{ is in } L^1(\BBR)\}$ where $\hat{f}$ is the Fourier 
transform of $f$. If $p\hat{M}$ satisfies the quantum Duhamel formula 
for each real $p$ then $f(\ol{M})=\{f(\ol{M}_t):t\in [0,1]\}$ is a regular 
quantum semimartingale and {\it satisfies the quantum Ito 
formula}
\beqn\label{EE 72.02}
f(\ol{M}_t)=f(0)+\int_0^t(\hat{E}_{f(M)}\,d\GL
+\hat{F}_{f(M)}\,dA+\hat{G}_{f(M)}\,dA\d+\hat{H}_{f(M)}\,ds)
\eeqn
where, for $\hat{X}$ in $\{\hat{E},\hat{F},\hat{G}, \hat{H}\}$,
$$
\hat{X}_{f(M)}=\int_{-\infty}\iiy \hat{f}(p) \hat{X}_{{\rm exp}(ipM)}\,dp
$$
This is the case whenever 
$M$ is regular and self-adjoint \cite{Vin2}.
\bb\ni
The original formula of Hudson and Parthasarathy will be referred to as 
the {\em quantum Ito product formula}.
\bb\ni
It was shown in \cite{Vin2} that (\ref{EE 72.02}) implies 
the Ito formula for a 
large class of classical semimartingales, which may even have 
jumps. Unfortunately this class does not include Brownian motion.
\bb\ni
The main objectives in this article are to
\bb\ni
(a) identify a non-trivial class, $\CC$, containing $B$, of 
essentially self-adjoint quantum semimartingales;
\bb\ni
(b) show that if $\hat{M}$ belongs to $\CC$ 
then $p\hat{M}$ satisfies the quantum Duhamel formula for all real 
$p$ and $f(\hat{M})$ satisfies the quantum Ito formula for all 
$f\in C^{2+}(\BBR)$;
\bb\ni
(c) enlarge the class $\CC$ identified in (a) and (b) by showing 
that if $\hat{M}$ is an essentially self-adjoint quantum 
semimartingale satisfying the quantum Duhamel formula 
then so is the perturbation $\hat{M}+\hat{J}$ 
whenever $\hat{J}$ is a regular self-adjoint quantum semimartingale.
\bb\ni
We also touch briefly on the quantum  Stratonovich formula.
\bb\ni
We approach (a) and (b) {\em via} ``chaos matrices".
\bb\ni
Bose--Fock space may be defined as the direct sum of Hilbert spaces, 
sometimes called chaos spaces: 
$$ 
\FH = \FH^{0}\oplus  \FH^{1}\oplus \cdots\oplus \FH^{j}\oplus \cdots,. 
$$ 
$\FH^j$ is a closed subspace of the $j$-fold tensor product
$L^2[0,1]^{\otimes j}$. 
The symmetric tensors $\{f\otimes\cdots\otimes f:f\in L^2[0,1]\}$ 
are  in $\FH^j$. The {\it exponential vectors}, 
$\{e(f):f\in L^2[0,1]\}$ where
$$
e(f)=1\oplus f\oplus \frac{f^{\otimes 2}}{(2!)^{1/2}}+\cdots
+\frac{f^{\otimes j}}{(j!)^{1/2}}+\cdots,
$$
are total in $\FH$. 
\bb\ni 
A matrix $T$ with entries $T^i_j\in \CB(\FH^j,\FH^i)$ is called a chaos 
matrix. Every bounded operator in $\FH$ has a unique chaos matrix 
representation. The same is not true of unbounded 
operators.
\bb\ni
The chaos matrix approach is suited to Lindsay's 
construction of quantum stochastic integrals {\em via} 
the gradient process $\nabla$ and the 
Hitsuda--Skorohod process $\CS$ \cite{Lin1}. These are processes 
of unbounded operators with simple superdiagonal and subdiagonal 
chaos matrix representations. Lindsay's integrals are 
defined for non-adapted processes and are an extension of Hudson 
and Parthasarathy's integrals.
\bb\ni
The basic processes $\GL,A,A\d,t$ of quantum stochastic calculus
have uncomplicated chaos matrix representations. Representations
leading naturally to a definition of the 
chaos matrix quantum stochastic integral of a chaos matrix-valued 
process with respect to $d\GL$, $dA$, $dA\d$ and $ds$.
\bb\ni
The correspondence between operator-valued and matrix processes 
is preserved under the respective quantum stochastic integrations. 
If the operators $\hat{E}_t$, $\hat{F}_t$, $\hat{F}\st_t$ and 
$\hat{H}_t$ in the quantum semimartingale
(\ref{E 71.202}) are 
represented by chaos matrices, ${E}_t$, ${F}_t$, ${F}\st_t$ and ${H}_t$ and
\begin{eqnarray*}
{M}_t={M}_0+\int_0^t ({E}\,d\GL+{F}\,dA +{F}\st\,dA\d+{H}\,ds),
\end{eqnarray*}
then $M_t$ is the chaos matrix representation of $\hat{M}_t$ for 
$t\in [0,1]$.
\bb\ni
The entries of the chaos matrix processes $E,F,F\st$ and $H$ are 
themselves operator valued processes that belong to 
Lebesgue spaces:
$$
E^i_j\in L\iiy([0,1],\CB(\FH^j,\FH^i)\so),\quad
F^i_j, (F\st)^i_j\in L^2([0,1],\CB(\FH^j,\FH^i)\so), 
$$
$$
H^i_j\in L^1([0,1],\CB(\FH^j,\FH^i)\so).
$$
These spaces have Bochner--Lebesgue norms but measurability is 
with respect to the strong operator topology on 
$\CB(\FH^j,\FH^i)$.
\bb\ni
The norms are used to define a control 
matrix, $\LGk= \LGk(M)$:
$$
\LGk^i_j=i^{1/2}\nm{E^{i-1}_{j-1}}\iy j^{1/2}+
\nm{F^i_{j-1}}_2 j^{1/2}+i^{1/2}\nm{(F\st)^{i-1}_j}_2+\nm{H^i_j}_1,
\quad i,j=0,1,\ldots.
$$
The control matrix plays a r\^ole similar to that of the dominating 
function in Lebesgue's dominated convergence theorem.  
\bb\ni
Let $\ell^2$ be the Hilbert space of complex sequences 
$x=(x_0,x_1,\ldots,x_n,\ldots)^t$, written as infinite column 
vectors, with $\nm{x}_2=\left(\sum_{n=0}\iiy |x_n|^2\right)^{1/2}<\infty$. 
Denote by $\ell_{00}$ the dense subspace of sequences with finitely many 
non-zero entries.
\bb\ni
The matrix $\LGk$ represents an operator in $\ell^2$ (possibly with 
domain $\{0\}$.)  A vector $x\in \ell^2$ is analytic for ${\LGk}$ 
if $x\in \CD({\LGk}^n)$ for $n=1,2,\ldots$ and 
$$ 
\sum_{n=0}\iiy z^n\frac{\|{\LGk}^n x\|}{n!}
$$
has non-zero radius of convergence $r(x)$. The vector space of analytic 
vectors $x$ with $r(x)\geq \Ge$ is denoted $\CA_{\Ge}(\hat{\LGk})$.
\bb\ni
From the results of Sections \ref{QDF} and \ref{QIF} 
we obtain the following theorem 
which allows us to achieve the objectives (a) and (b) above. Let 
$\CE\subset \FH$ be the linear span of the exponential vectors.
\begin{Thm} \nonumber Let $\hat{M}$ be a symmetric quantum semimartingale 
whose chaos matrix representation $M$ has control matrix $\LGk$
and let $\Ge>0$. 
\bb\ni
If $\ell_{00}\subset \CA_{\Ge}({\LGk})$ then 
\bb\ni
{\rm (i)} $\hat{M}$ is essentially self-adjoint with core 
$\CE$.
\bb\ni
{\rm (ii)} $pM$ satisfies the quantum Duhamel formula for each real 
$p$.
\bb\ni
{\rm (iii)} $f(\ol{M})$ satisfies the quantum Ito formula for each 
$f\in C^{2+}(\BBR)$.
\end{Thm}
The conditions of this theorem are satisfied in particular
when ${E}$, ${F}$, ${F}\st$ and ${H}$ are adapted processes of 
$(2k+1)$-diagonal chaos matrices. It is shown in Section 
\ref{PDK} that there is a large class of such adapted processes 
defined by kernels.
\bb\ni
The proof of essential self-adjointness uses Nelson's analytic 
vector theorem. A chaos matrix version of 
the quantum Ito product formula is then proved 
for powers of $M$. A chaos matrix version of the quantum Duhamel 
formula and the quantum
Ito formula are then obtained in much the same way as for regular 
quantum semimartingales. The operator versions of the formulae 
then follow from the correspondence between quantum 
semimartingales and their chaos matrix representations.
\bb\ni
The proof of essential self-adjointness is in some respects an 
extension of the proof of essential self-adjointness, 
on a common core, of the field operators $\Phi(f)$ in quantum 
statistical mechanics \cite[Proposition 5.2.3]{BR2} (see also 
\cite[\S 19.3]{GJ}).
\bb\ni
We adopt two approaches to (c). The terms on each side of (\ref{EE 72.01})are replaced by their Duhamel expansions.
Using elementary calculus sides are then rearranged and shown to be 
equal. 
\bb\ni
Using a more sophisticated method we prove (\ref{EE 72.01}) 
in the special case that $M$ is obtained from a Brownian 
martingale {\em via} the Weiner--Ito transformation.
Suppose that $W$ is classical Brownian motion and ${\sf F}$ is a 
real-valued bounded adapted process in Wiener space and let 
${\sf M}$ be the classical martingale
\beqn\label{E 72.552}
{\sf M}_t=\int_0^t {\sf F}_s\,dW_s.
\eeqn
The Wiener-Ito isomorphism identifies $W$ with $A+A\d$, the Brownian 
quantum semimartingale, ${\sf F}$ with an adapted operator-valued process 
$\hat{F}$ and ${\sf M}$ may be identified with a quantum semimartingale 
$\hat{M}$. It follows from the classical Ito 
formula that the quantum semimartingale $\hat{M}$ satisfies the 
quantum Ito formula for each $f\in C^{2+}(\BBR)$.
\bb\ni
We then show that if  $\hat{J}$ is a regular 
self-adjoint quantum semimartingale $\hat{M}+\hat{J}$ satisfies 
the quantum Ito formula. The method is to 
approximate $\hat{M}$ by a sequence $h_n(\hat{M})$ of regular 
quantum semimartingales. The regular quantum semimartingale 
$\hat{N}\ub{n}=h_n(\hat{M})+\hat{J}$ then satisfies the 
quantum Ito formula. A limiting argument shows that $f(\hat{M})$
satisfies the quantum Ito formula.
\bb\ni
This article is rather long even though the ideas contained in 
it, outlined above, are simple. In our attempt to 
keep within a functional analytic category
we have had to start from scratch in some places.
\bb\ni
The infinite matrix approach to operator theory has not 
proved very  popular: it gives nothing new for bounded operators 
and the correspondence between unbounded operators and chaos 
matrices is unclear. There seem to be few results in the 
literature and we have had to prove those we need.
Intrinsic conditions must be 
found for a chaos matrix valued process to be adapted.
The cmx quantum stochastic integrals 
defined in this paper are processes of chaos matrices and a 
balance sheet  
must be kept to show that they represent the corresponding 
Hudson--Parthasarathy processes. 
Another complication is that although the integrands for 
quantum semimartingales belong 
to Lebesgue spaces they do not necessarily belong to the standard 
Bochner--Lebesgue spaces.
Finally the conditions for unbounded operators to be essentially 
self-adjoint are delicate and must be checked in detail. 
\bb\ni
To keep the article to its present length
we let the time interval be $[0,1]$, start 
our processes at $0$, consider functions $f(x)$ rather than 
$f(x,t)$, dispense with the initial space and do not consider the 
cases of finite and countably infinite multiplicity.
\bb\ni
The following is a special case of Corollary \ref{T 81.362}
which overcomes the difficulty with Brownian motion.
\bb\ni
\begin{Thm} \nonumber 
Let $\hat{F}$ in the Bochner--Lebesgue space $
L^2([0,1],\CB(\FH))$ be an adapted process 
and let $M$ be the quantum semimartingale
$$
\hat{M}_t=\int_0^t (\hat{F}_s\,dA_s+\hat{F}\st_s dA\d_s).
$$
If $F_t$ is the chaos matrix representation $\hat{F}_t$ then 
$t\mapsto F_t{\,}^i_j$ belongs to Bochner--Lebesgue space
$L^2([0,1],\CB(\FH^j,\FH^i))$.
\bb\ni
The control matrix $\LGk$ with entries 
$$
\LGk^i_j =i^{\frac12}\|F^{j}_{i-1} \|_2
+\|F^{i}_{j-1} \|_2 j^{\frac12} 
$$
represents a symmetric linear transformation $\hat{\LGk}$ in 
$\ell^2$. 
\bb\ni
If $\ell_{00}$ consists of analytic vectors for $\hat{\LGk}$ with radius of 
convergence greater than fixed $\Ge$ then
\bb\ni
{\rm (i)} $\hat{M}$ is a process of essentially self-adjoint operators;
\bb\ni
{\rm (ii)} If $f\in C^{2+}(\BBR)$ then
$$
f(\hat{M}_t)=f(0)+\int_0^t \left(Df(\hat{M}_s)(d\hat{M}_s)+
D^2_I f(\hat{M}_s)(d\hat{M}_s,d\hat{M}_s)\right)
$$
\end{Thm}
When $F=I$ this theorem implies the Ito formula for classical 
Brownian motion.
\bb\ni
The article is organised as follows. 
Section \ref{CMX} is a general discussion of chaos matrices and 
the operators they represent. 
Section \ref{BP} recalls the definitions and 
some properties of the basic processes occurring 
in the theories of Hudson--Parsatharathy and Lindsay and shows 
that they may be represented as processes of chaos matrices. 
Adapted chaos matrix valued processes are characterised in 
Section \ref{S 3}.
\bb\ni
The cmx quantum stochastic integral $M_t(E,F,G,H)$ of a quadruple 
of chaos matrix valued processes $(E,F,G,H)$ 
is defined in Section \ref{QSI}. The integrands $(E,F,G,H)$ are 
not required to be adapted. Bounds are found for the entries in $M_t$. 
When $(E,F,G,H)$ are adapted processes of bounded chaos matrices 
it is shown that $M_t$ represents the Hudson-Parthasarathy 
quantum stochastic integral of the corresponding operator 
processes. This shows that the two definitions of quantum 
stochastic integral are consistent. 
\bb\ni
Section \ref{SMX} is devoted to properties of scalar matrices. In Section 
\ref{HPP} the integrands are allowed to be unbounded processes 
and sufficient conditions are found 
for a cmx quantum stochastic integrals 
to represent the corresponding 
Hudson--Parthasarathy quantum stochastic integral.
\bb\ni
The remainder of the article is devoted to the quantum 
Ito formula. A product Ito formula is proved for adapted chaos matrix valued 
processes in Section \ref{IPF}, a quantum Duhamel formula is proved in 
Section \ref{QDF} and the quantum Ito formula in Section 
\ref{QIF}.
\bb\ni
We give, without proof, in Section \ref{QSF} a quantum 
Stratonovich formula for regular quantum semimartingales and pose 
the problem of generalising it to the irregular case.
\bb\ni
In Section \ref{PCM} it is shown that 
a classical Brownian quantum semimartingale perturbed by a 
regular quantum semimartingale still satisfies the 
quantum Ito formula. The methods of this section 
seem to work only when the perturbed processes are classical. 
\bb\ni
The Duhamel expansion is reviewed in Section \ref{DHE} and 
used in Section \ref{PQS} to show that the set of 
essentially self-adjoint quantum semimartingales which satisfy 
the quantum Duhamel formula is closed under 
perturbations by regular self-adjoint quantum semimartingales.
\bb\ni
In Section \ref{PROB} the article concludes with five open problems.

\section{\label{CMX} Chaos Matrices.} 
A {\it chaos decomposition} of a Hilbert space $\CH$ is a 
decomposition of $\CH$ into the direct sum
$$ 
\CH = \CH^{0}\oplus  \CH^{1}\oplus 
\cdots\oplus \CH^{j}\oplus \cdots
$$ 
of Hilbert spaces $\CH^0$, $\CH^1$, $\CH^2,\ldots$. The Hilbert space $\CH^j$ 
is called the {\it $j$}th {\it chaos} and $\pi^j$ denotes the 
orthogonal projection onto $\CH^j$. Each $\psi\in \CH$ is written as 
a column vector 
$\psi=(\psi^{0},\psi^{1},\cdots,\psi^{j},\cdots)^t$ called the {\it 
chaos representation} of $\psi$. 
\bb\ni 
Let $\CH_{00}=\{\psi\in \CH: \psi^j=0 \mbox{ whenever } 
j>N=N(\psi)\in \BBN\}$, the dense subspace of vectors with finite 
chaos representations. 
\bb\ni 
A {\it chaos matrix} (for the above decomposition of $\CH$) is a 
matrix $T$ with {\it entries} $ T^i_j\in 
\CB(\CH^j,\CH^i)$, $i,j=0,1,2,\ldots$. We often use `cmx' 
as shorthand for `chaos matrix'. For example, 
`representation'  is short for `chaos matrix representation'.
\bb\ni 
The {\it adjoint} of 
$T$ is $T\st\defeq[(T^j_i)\st]$. If $T=T\st$ then $T$ is a {\it 
symmetric} chaos matrix. 
\bb\ni
A chaos matrix $T$ is $(2k+1)$-{\it diagonal} if $T^i_j=0$ 
whenever $|i-j|>k$.
\bb\ni 
If $S$ and $T$ are chaos matrices 
and, for all $i,j=0,1,\ldots$, the series $\sum_{\nu=0}\iiy 
S^i_{\nu}T^{\nu}_j$ is convergent in the strong operator topology 
to $U^i_j$ in $\CB(\CH^j,\CH^i)$ then the chaos matrix 
$U=[U^i_j]$ is the {\it product} $ST$ of $S$ and $T$. 
This product will exist if, 
for example, either of $S$ and $T$ is $(2k+1)$-diagonal for some 
$k\in \BBN$.
\bb\ni
The existence of $ST$ does not necessarily imply the existence of $T\st 
S\st$. 
\begin{Lemma}\label{L 73.00}Let $S$ and $T$ be chaos matrices. If $ST$ and 
$T\st S\st$ both exist then $(ST)\st=T\st S\st$.
\end{Lemma}
\pf If $\phi\in \CH^i$ and $\psi\in \CH^j$ then
$$
\ip{ST\psi}{\phi}=\sum_{\nu=0}\iiy\ip{S^i_{\nu}T^{\nu}_j\psi}{\phi}
=\sum_{\nu=0}\iiy\ip{\psi}{(T\st)^j_{\nu}(S\st)^{\nu}_i\phi}
=\ip{\psi}{T\st S\st \phi}.\,\,\epf
$$
Let $\CL(\CH)$ denote the set of linear transformations in $\CH$. 
To each chaos matrix $T$ corresponds $\tilde{T}$ in $\CL(\CH)$ 
as follows.  Let 
$$ 
\FD(T)=\left\{\psi\in \CH: \ba{l} 
\sum_{j=0}\iiy T^i_j \psi^j\,\mbox{ is convergent to 
$\eta^i(\psi)\in \CH^i$ for all $i$} \\ \hspace{1.2in}\mbox{and} 
\\ \mbox{ $\eta(\psi) 
=(\eta^0(\psi),\eta^1(\psi),\eta^2(\psi),\ldots)^t$ belongs to 
$\CH$} \ea \right\}. 
$$ 
Define $\tilde{T}$ by letting $\CD(\tilde{T})=\FD(T)$ and putting
$$ 
\tilde{T}\psi=\eta(\psi)\qquad\mbox{ whenever $\psi\in 
\FD(T)$} 
$$ 
Care must be taken at this point. If 
$T$ is a chaos matrix and the product chaos matrix $T^2$ exists it is 
not necessarily true that 
$$ 
\FD(T^2)=\CD(\tilde{T}^2)=\{\psi\in 
\CD(\tilde{T}):\tilde{T}\psi\in \CD(\tilde{T})\}. 
$$ 
It may be that $\psi\in \FD(T^2)$ yet $\psi\not\in \FD(T)$. 
\bb\ni 
The {\it $C\iiy$-vectors} of a chaos matrix $T$ are
$$
C\iiy(T)=\left\{ \ba{cl}
\bigcap_{k=0}\iiy \FD(T^k)&\quad\mbox{ if $T^k$ exists 
for $k=1,2,\ldots$,} \\ 
\emptyset &\quad\mbox{ otherwise.} 
\ea \right. 
$$
A vector $\Gv\in \FH$ is {\it analytic} for $T$ 
if $\Gv\in C\iiy(T)$  and the complex power series 
\beqn\label{E 73.732} 
\sum_{k=0}\iiy \frac{\|T^k \Gv\|}{k!}z^k
\eeqn
has non-zero radius of convergence $r(\Gv)=r(\Gv,T)$.
Then $\CA(T)$ is the vector space of analytic vectors of $T$ and 
$\CA_r(T)=\{\Gv\in \CA(T):r(\Gv)\geq r\}$ for $0<r\leq \infty$. 
\bb\ni 
Let $T$ be a chaos matrix and let $\hat{T}$ belong to $\CL(\CH)$. 
Then $T$ is a {\it chaos matrix representation} of $\hat{T}$ if 
$$
{\rm(i)}\quad \CD(\hat{T})\subset \FD(T);\qquad {\rm(ii)}\quad
T\Gv=\hat{T}\Gv \mbox{ whenever }\Gv\in \CD(\hat{T}).
$$
Thus $T$ is a cmx representation of $\hat{T}$ if and only if 
$\hat{T}$ is an extension of $\tilde{T}$.
We shall not consider operators with more than one chaos-matrix 
representation. 
\bb\ni
The number operator $\hat{N}$ in $\CH$ has 
unique diagonal chaos-matrix  representation 
$$ 
N={\rm diag}[0,I_1,2I_2\ldots,jI_j\ldots], 
$$ 
where $I_j$ is the identity transformation
in $\CH^j$. Clearly $\FD(N)=\{\psi\in 
\CH:(0,\psi^1,2\psi^2,\ldots,j\psi^j,\ldots)^t\in 
\CH\}=\CD(\tilde{N})$ and $\tilde{N}=\hat{N}$. 
\bb\ni 
A chaos matrix $T$ is { \it bounded} if $\FD(T)=\CH$ 
and $\tilde{T}$ is bounded.  Let $\CB_{\rm cmx}(\CH)$ be the 
Banach space
$$ 
\CB_{\rm cmx}(\CH)= \{T:\mbox{ 
$T$ is a bounded chaos matrix for $\CH$}\} 
$$ 
with $\nm{T}\defeq\|\tilde{T}\|$ whenever $T\in \CB_{\rm 
cmx}(\CH)$. The completeness of $\CB_{\rm cmx}(\CH)$ follows from 
part (iii) of the following theorem.
\begin{Thm}\label{T 73.001} {\rm (i)} Each $\hat{T}\in \CB(\CH)$ 
has a unique cmx representation $T$ whose entries are given by 
\beqn\label{E 73.02} T^i_j \psi^j=(\hat{T}\psi^j)^i \qquad 
\psi^j\in \CH^j 
\eeqn 
Moreover $\nm{T^i_j}\leq \nm{T}$ for all $i,j=0,1,2,\ldots$.
\bb\ni
{\rm (ii)} A chaos matrix $T$ is bounded 
if and only if $\FD(T)=\CH$. 
\bb\ni
{\rm (iii)} If $T$ is a chaos matrix with $\CH_{00}\subset \FD(T)$ and 
$\tilde{T}$ is bounded then $T$ is bounded.
\bb\ni 
{\rm (iv)}  With the product and adjoint operations defined 
above $\CB_{\rm cmx}(\CH)$ is a von Neumann algebra. 
The mapping $T\mapsto \tilde{T}$ is an isometric algebraic $\ast$-isomorphism 
from $\CB_{\rm cmx}(\CH)$ onto $\CB(\CH)$. 
\end{Thm} 
\pf (i) It follows from the definition that any cmx 
representation of $\hat{T}$ satisfies (\ref{E 73.02}).
If $\psi^j\in \CH^j$ then $\|(\hat{T}\psi^j)^i\|\leq \|\hat{T}\psi^j\|\leq 
\|\hat{T}\|\cdot\|\psi^j\|$ so that (\ref{E 73.02}) defines $T^i_j\in 
\CB(\CH^j,\CH^i)$. Therefore the chaos matrix $T$ is well 
defined by (\ref{E 73.02}). 
\bb\ni 
$\hat{T}$ is bounded and if $\psi\in \CH$ then $\sum_{j=0}\iiy T_i^j 
\psi^j= \sum_{j=0}\iiy(\hat{T} \psi^j)^i=(\hat{T}\psi)^i=\eta^i$ 
where $\eta=\hat{T} \psi$. Therefore $\FD(T)=\CH$ and $T$ 
uniquely represents $\hat{T}$.
\bb\ni 
(ii) If $T$ is bounded then $\FD(T)=\CH$ by definition.
\bb\ni
Suppose, conversely, that $\FD(T)=\CH$. Define the linear 
transformation $\tilde{T}^i:\CH\row \CH^i$ by
$$ 
\tilde{T}^i\psi=\sum_{j=0}\iiy 
T^i_j\psi^j=\sum_{j=0}\iiy (T^i_j \circ\pi^j)\psi=\eta^i 
$$ 
whenever $\psi\in \CH$. Since $T^i_j \circ\pi^j$ is bounded it 
follows from 
the uniform boundedness principle that $\tilde{T}^i\in 
\CB(\CH,\CH^i)$. It follows from the definition of $\FD(T)$ that 
$\sum_{i=0}\iiy{}\tilde{T}{}^i\psi$ is convergent to 
$\tilde{T}\psi$ for each $\psi\in \CH$ and by the uniform 
boundedness principle $\tilde{T}\in \CB(\CH)$. 
\bb\ni
(iii) The bounded linear operator $\tilde{T}$ has a unique 
extension $\hat{T}\in \CB(\CH)$. Since $\CH_{00}\subset \FD(T)$ 
the operator $T^i_j$ is defined by (\ref{E 73.02}). Therefore $T$ 
represents $\hat{T}$ and $\FD(T)=\CH$.
\bb\ni 
(iv) If $S,T\in \CB_{\rm cmx}(\CH)$ and 
$\psi^j\in \CH^j$ then 
$(T^0_j\psi^j,T^1_j\psi^j,\cdots,T^n_j\psi^j,\cdots)^t 
=\tilde{T}\psi^j\in \CH$. Since $\FD(S)=\CH$ 
\beqn\label{E 73.002} 
(\tilde{S}(\tilde{T}\psi^j))^i&=&\sum_{\nu=0}\iiy {S^i_{\nu} 
T^{\nu}_j\psi^j}.
\eeqn 
Therefore the series converges for each $\psi^j\in \CH^j$. By the 
uniform boundedness theorem  $\sum_{\nu=0}\iiy {S^i_{\nu} 
T^{\nu}_j}$ is convergent in the strong operator topology to 
$U^i_j$ in $\CB(\CH^j,\CH^i)$. Therefore the product 
chaos matrix $ST$ exists. 
\bb\ni 
It follows from (\ref{E 73.002}) 
that, for $\Gv^i\in \CH^i$ and $\psi^j\in \CH^j$, 
$$ 
\ip{\tilde{S}\tilde{T}\psi^j}{\Gv^i}=\ip{U^i_j \psi^j}{\Gv^i} 
=\ip{(ST)^i_j \psi^j}{\Gv^i}. 
$$ 
Therefore $ST$ is the chaos matrix of 
$\tilde{S}\tilde{T}$ and $\widetilde{ST}=\tilde{S}\tilde{T}$. 
\bb\ni 
If $\psi^j\in \CH^j$ and $\Gv^i\in \CH^i$ then 
$$ 
\ip{\tilde{T}\psi^j}{\Gv^i}=\ip{T^i_j\psi^j}{\Gv^i}=\ip{\psi^j}{(T^i_j)\st\Gv^i}= 
\ip{\psi^j}{(T\st)^j_i\Gv^i}=\ip{\psi^j}{\widetilde{T\st}\Gv^i}. 
$$ 
It 
follows from (i) that $T\st$ is the cmx representation of 
$(\tilde{T})\st$ and $\widetilde{(T\st)}=(\tilde{T})\st$. The 
transformation $T\mapsto\tilde{T}$ is clearly linear and is 
therefore an isometric $\ast$-isomorphism from 
$\CB_{\rm cmx}(\CH)$ onto $\CB(\CH)$. \,\,
\epf 
\bb\ni 
\def\cmx{_{\rm cmx}}
Let $\CB(\CH)\so$ denote $\CB(\CH)$ with the strong operator 
topology and let $\CB\cmx(\CH)\so$ denote $\CB\cmx(\CH)$ with the 
quotient topology induced by the $\ast$-isomorphism $T\mapsto 
\tilde{T}$.
\bb\ni
An {\it operator process} is a family 
$\hat{T}=\{\hat{T}_t:t\in [0,1]\}$ of linear operators.   
If $\hat{T}_t$ is bounded for almost all $t\in 
[0,1]$  then $\hat{T}$ is {\it regular}. Two  operator processes 
$\hat{T}$ and $\hat{T}\pr$ are identified if 
$\hat{T}_t=\hat{T}\pr_t$ for almost all $t\in [0,1]$.
\bb\ni
A {\it chaos matrix process} ({\it cmx process}) is a family 
$T=\{T_t:t\in [0,1]\}$ of chaos matrices.  
If $T_t$ is a bounded chaos matrix for almost all $t\in 
[0,1]$ then $T$ is {\it regular}. Two cmx processes $T$ and 
$T\pr$ are identified if $T_t=T\pr_t$ for almost all $t\in [0,1]$.
\bb\ni 
A cmx process $T$ is a {\it cmx representation} of an operator process 
$\hat{T}$ if $T_t$ is a cmx representation of $\hat{T}_t$ for 
almost all $t\in[0,1]$.
\bb\ni
The integrands of the 
quantum semimartingales considered by Attal and Meyer are 
regular processes 
$\hat{F}=\{\hat{F}_t:t\in [0,1]\}$ where $\hat{F}_t$ 
is a bounded operator in the
Bose-Fock space $\FH$ for each $t\in [0,1]$. Since $\CB(\FH)$ is 
non-separable the measurability 
conditions on these integrands, defined with respect to the 
strong operator topology on $\CB(\FH)$, are too weak to ensure that $\hat{F}$ 
belongs to one of the Bochner--Lebesgue spaces 
$L^p([0,1],\CB(\FH))$. Nor are they strong enough to ensure 
that $\hat{F}$ is Pettis integrable. The integrand $\hat{F}$ belongs 
instead to one of the following Lebesgue spaces.
\bb\ni
If $1\leq p\leq\infty$ let $L^p([0,1],\CB(\CH)\so)$ the normed vector space 
of processes
$\hat{F}:[0,1]\row \CB(\CH)$ such that
\bb
(i) $t\mapsto \hat{F}_t\xi$ is measurable for each $\xi$ in some dense 
subset of $\CH$;
\bb
(ii) $t\mapsto \nm{\hat{F}_t}$ belongs to $L^p([0,1])$,
\bb\ni
with norm, $\nm{\cdot}_p$, given by
$$
{\renewcommand{\arraystretch}{2}
\nm{\hat{F}}_p=\left\{ \begin{array}{cl}
\dyl \left(\int_0^1 \nm{\hat{F}_t}^p\,dt\right)^{\frac1p}, &\qquad 
1\leq p <\infty, \\
\dyl \sup_{t\in [0,1]}\nm{\hat{F}_t},&\qquad p=\infty.
\end{array}\right. }
$$
Suppose $\xi$ and $\xi\ub{n}$ in $\CH$ are such that $\nm{\xi-
\xi\ub{n}}<1/2^n$. Then $t\mapsto \hat{F}_t\xi$ and $t\mapsto 
\hat{F}_t\xi\ub{n}$ belong to the Bochner--Lebesgue space 
$L^p([0,1],\CH)$ and $\nm{\hat{F}_{\Cdot}\xi-\hat{F}_{\Cdot}\xi\ub{n}}
<\nm{\hat{F}}_p/2^n$. 
Therefore $\hat{\hat{F}}_t\xi\ub{n}$ converges to $\hat{\hat{F}}_t\xi$ 
for almost all $t$ 
in $[0,1]$ \cite[Theorem 3.11]{Rud}. This means that (i) may be replaced by the 
condition
\bb
(i)$\pr$ $t\mapsto \hat{F}_t\xi$ is measurable for each $\xi$ in $\CH$.
\bb\ni
Condition (i) is used rather than (i)$\pr$ because the 
measurability condition on an integrand, $\hat{F}$, 
of a quantum semimartingale is that
$t\mapsto \hat{F}_t e(f)$ be measurable whenever the exponential 
vector $e(f)$ belongs to a given total set of exponential vectors 
in $\FH$.
\bb\ni
A process $\hat{F}$ satisfying (i)$\pr$ is {\it strongly measurable}. The 
composition of strongly measurable processes is strongly measurable.
\begin{Cor} \label{CC 73.01}If $\hat{F}, \hat{G}:[0,1]\row \CB(\CH))$ 
are strongly measurable then so is $\hat{F}\circ \hat{G}:
[0,1]\row \CB(\CH)$, where $(\hat{F}\circ \hat{G})_t=\hat{F}_t\circ \hat{G}_t$.  
\end{Cor}
\pf If $\xi\in \CH$ there exists a sequence $\Gv\ub{n}$ of 
$\CH$-valued step functions such that $\hat{G}_t\xi=\lim_{n\roy} 
\Gv_t\ub{n}$ almost everywhere. But $t\mapsto \hat{F}_t
\otimes\hat{G}v_t\ub{n}$ is 
measurable and converges almost everywhere to 
$\hat{F}_t\otimes \hat{G}_t\xi$ so 
that $t\mapsto \hat{F}_t \otimes\hat{G}_t\xi $ is measurable. \,\,\epf
\bb\ni
\begin{Lemma} If $1\leq p\leq\infty$ then 
$L^p([0,1],\CB(\CH)\so)$ is a Banach space.
\end{Lemma}
\pf We follow the line of argument in \cite{Rud}. 
Suppose $1\leq p<\infty$ and 
let $\hat{F}\ub{n}$ be a Cauchy sequence in $1\leq p\leq\infty$. By 
extracting a subsequence it may be assumed that $\nm{\hat{F}\ub{m}-
\hat{F}\ub{n}}<1/2^n$ whenever $m>n$. This implies that $\hat{F}\ub{n}_t$ 
is norm convergent in $\CB(\CH)$ to $\hat{F}_t$ for almost all $t$ in 
$[0,1]$. If $\xi \in \CH$ then  $\hat{F}\ub{n}_t\xi$ is almost everywhere 
convergent to $\hat{F}_t\xi$ and $t\mapsto \hat{F}_t\xi$ is measurable.
\bb\ni
By Fatou's lemma
\beqn\label{E 073.01}
\int_0^1 \nm{\hat{F}_t-\hat{F}\ub{n}_t}^p\,dt
\leq \liminf_{m\roy}\int_0^1\nm{\hat{F}\ub{m}_t-\hat{F}\ub{n}_t}\,dt
\leq 
\frac{1}{2^np},
\eeqn
and $t\mapsto \nm{\hat{F}_t-\hat{F}\ub{n}_t}$ belongs to $L^p[0,1]$. Thus $\hat{F}-
\hat{F}\ub{n}$ and therefore $\hat{F}$ belong to $L^p([0,1],\CB(\CH)\so)$. It 
follows from (\ref{E 073.01}) that $\hat{F}\ub{n}$ is norm convergent 
to $\hat{F}$ so that $L^p([0,1],\CB(\CH)\so)$ is a Banach space.
\bb\ni
A straightforward argument shows that  
$L\iiy([0,1],\CB(\CH)\so)$ is also a Banach space. \,\,\epf
\bb\ni
We now define the integral of $\hat{F}\in L^p([0,1],\CB(\CH)\so)$ for 
$1\leq p\leq\infty$. It follows from (ii) that $t\mapsto 
\nm{\hat{F}_t\xi}$ 
belongs to $L^p[0,1]$ and from (i)$\pr$ that $t\mapsto \hat{F}_t \xi$ 
belongs to the Bochner--Lebesgue space $ L^p([0,1],\CH)$ for each 
$\xi\in \CH$ \cite[II.2. Theorem 2]{DiU}. Define the indefinite 
integral, $I_{\Cdot}(\hat{F})$, of $\hat{F}$ by putting 
$$
 I_t(\hat{F})\xi\equiv\left(\int_0^t \hat{F}_s\,ds\right)\xi\defeq 
\int_0^t \hat{F}_s\xi\,ds,\qquad \xi\in \CH,\quad t\in [0,1].
$$
If $\eta \in \CH$ then
$$
|\ip{I_t(\hat{F}) \xi}{\eta}|=\left|\int_0^t\ip{\hat{F}_s \xi}{\eta}\,ds\right|
\leq \int_0^t\nm{\hat{F}_s}\,ds \nm{\xi}\,\nm{\eta}
$$
$$
\leq t^{\frac1q}\nm{\hat{F}}_p \nm{\xi}\,\nm{\eta}.
$$
This shows that the indefinite integral $I_t(\hat{F})$ 
belongs to $\CB(\CH)$ with $\nm{I_t(\hat{F})}\leq t^{1/q}\nm{\hat{F}}_p$. 
\bb\ni
The vector space of cmx processes is denoted $\CL_{\rm 
cmx}(\CH)$.
Define the Frechet spaces
$$
\ba{rcl}
\CL^p_{\rm cmx}(\CH)\so&=&\{T\in \CL_{\rm cmx}(\CH):T^i_j \in 
L^p([0,1],\CB(\CH^j,\CH^i)\so)\},
\\
\CL^p_{\rm cmx}(\CH)&=&\{T\in \CL_{\rm cmx}(\CH):T^i_j \in 
L^p([0,1],\CB(\CH^j,\CH^i))\},
\ea
\quad 1\leq p\leq \infty,
$$
with the topology defined by the seminorms 
$X\mapsto \nm{X^i_j}_p$ $i,j=0,1,\ldots$ where 
$\nm{\cdot}_p$ is the norm on $L^p([0,1],\CB(\CH^j,\CH^i)\so)$ defined above.
\bb\ni
Thus a sequence $S\ub{n}$ in $\CL^p_{\rm cmx}(\CH)$ is {\it convergent} to $S$ in 
$\CL^p_{\rm cmx}(\CH)$ if $(S\ub{n})^i_j$ 
is convergent in $L^p([0,1],\CB(\CH^j,\CH^i))$ to $S^i_j$ for all $i,j$. 
A series $S\ub{n}$ is {\it summable} ({\it absolutely summable})  
if the series $(S\ub{n})^i_j$ is summable (absolutely summable) for each $i,j$. 
\bb\ni
Suppose $p,q,r\in [1,\infty]$ with $1/p+1/q=1/r$. 
If $X\in \CL^p_{\rm cmx}(\CH)$ and $Y\in \CL^q_{\rm cmx}(\CH)$ 
then the formula
$$
(X^i_{\nu} Y^{\nu}_j)_t=X_t{\,}^i_{\nu}Y_t{\,}^{\nu}_j
$$
defines, for each $\nu\in \BBN$, a process $X^i_{\nu} Y^{\nu}_j\in 
L^r([0,1],\CB(\CH^j,\CH^i))$. 
If $\sum_{\nu=0}\iiy X^i_{\nu}Y^{\nu}_j$ converges to $Z^i_j$ in 
$L^r([0,1],\CB_{\rm so}(\CH^j,\CH^i))$ for each $i,j\in \BBN$ 
the cmx process $Z$ is the {\it product} of $X$ and $Y$. We put 
$Z=XY$.
\bb\ni
If $1\leq p\leq\infty$ define $L^p([0,1],\CB\cmx(\CH)\so)$ to be 
the normed vector space of cmx processes $F$ such that 
\bb
(i) $t\mapsto F_t{\,}^i_j \xi^j$ is measurable for each $\xi^j$ 
in a dense subset of $\CH^j$;
\bb
(ii) $t\mapsto \nm{F_t}$ belongs to $L^p[0,1]$.
\bb\ni
As in the case of processes of bounded operator condition (i) 
is equivalent to 
\bb
(i)$\pr$ $t\mapsto F_t{\,}^i_j \xi^j$ is measurable for each $\xi^j$ 
in $\CH^j$.
\bb\ni
We note the following inclusions between these spaces.
\begin{Propn} { \ }
\bb\ni
{\rm(i)} $L^p([0,1],\CB\cmx(\CH)\so)$ is a closed subspace of 
$\CL^p\cmx(\CH)\so$.
\bb\ni
{\rm(ii)} $L^p([0,1],\CB\cmx(\CH))$ is a closed subspace of 
$\CL^p\cmx(\CH)$.
\bb\ni
{\rm(iii)} $L^p([0,1],\CB\cmx(\CH))$ is a closed subspace of 
$L^p([0,1],\CB\cmx(\CH)\so)$.
\bb\ni
{\rm(iv)} $L^p([0,1],\CB(\CH))$ is a closed subspace of 
$L^p([0,1],\CB(\CH)\so)$
\end{Propn}
We also have the following correspondences between 
regular operator processes and regular cmx processes.
\bb\ni
If $\hat{X}\in L^p([0,1],\CB(\CH)\so)$ define the cmx process
$X$ by letting  $X_t$ be the cmx representation of $\hat{X}_t$ 
whenever $t\in [0,1]$.
\bb\ni
If $X\in L^p([0,1],\CB\cmx(\CH)\so)$ define the operator process
$\tilde{X}$ by putting $\tilde{X}_t=(\widetilde{X_t})$ whenever $t\in 
[0,1]$. 
\begin{Propn} \label{PP 73.101} 
If $1\leq p\leq \infty$ then the mapping $\hat{X}\mapsto 
X$ is an isometric isomorphism from
\bb
{\rm(i)} $L^p([0,1],\CB(\CH)\so)$ onto $L^p([0,1],\CB\cmx(\CH)\so)$;
\bb
{\rm(ii)} $L^p([0,1],\CB(\CH))$ onto $L^p([0,1],\CB\cmx(\CH))$.
\bb\ni
Moreover $\tilde{X}=\hat{X}$.
\end{Propn}
The process $X$ is the {\it cmx representation} of $\hat{X}$.
\bb\ni
The relationship between unbounded chaos matrices and the 
operators they represent is not 
straightforward. The chaos matrices are a vector space for 
entrywise addition and scalar multiplication while 
$\CL(\CH)$ is not a vector space. There are non-trivial chaos matrices which 
represent only the trivial operator with domain $\{0\}$. Although 
a chaos matrices always represents a linear transformation 
not every linear transformation has a chaos matrix 
representation. Nor is the correspondence between a 
chaos matrix $T$ and the operator $\tilde{T}$ always  one to one. 
\bb\ni 
However under some conditions symmetric chaos matrices 
represent unbounded essentially self-adjoint linear transformations.
\bb\ni
If $X$ is in $\CL^p_{\rm cmx}(\CH)\so$ then $\int_0^t X_s\,ds$ is 
the cmx process with entries
$$
\int_0^t X_s{\,}^i_j\,ds. 
$$
\begin{Thm}\label{T 73.002}
Let $X$ be a process in $L^p([0,1],\CB_{\rm cmx}(\CH)_{\rm 
so})$ where $1\leq p\leq\infty$. If $1/p+1/q=1$ and 
$i,j=0,1,2,\ldots$ then, for $0\leq t\leq 1$,
\bb\ni
{\rm(i)} $X^i_j\in L^p([0,1],\CB(\CH^j,\CH^i)_{\rm so})$ with 
$\nm{X^i_j}_p\leq \nm{X}_p$.
\bb\ni
{\rm(ii)} The linear transformation $Y_t{\,}^i_j$ defined by the 
formula
$$
Y_t{\,}^i_j \psi=\int_0^t X_u{\,}^i_j\psi\,du,\qquad \psi\in \CH^j,
$$
belongs to $\CB(\CH^j,\CH^i)$ and $\|Y_t{\,}^i_j\|\leq 
t^{1/q}\nm{X^i_j}_p \leq t^{1/q}\nm{X}_p$.
\bb\ni
{\rm(iii)} The series $\sum_{j=0}\iiy Y_t{\,}^i_j \psi^j$ is convergent for 
each $\psi\in \CH$. 
\bb\ni
{\rm(iv)} The chaos matrix $Y_t{\,}$ with entries $Y_t{\,}^i_j$ is 
bounded with 
$\nm{Y_t{\,}}\leq t^{1/q}\nm{X}_p$.
\bb\ni
{\rm(v)} If $\psi\in \CH$ then 
$$
\tilde{Y}_t\psi=\int_0^t \tilde{X}_u\psi\,du. 
$$   
\end{Thm}
\pf (i) If $\psi\in \CH^j$ then $X_u{\,}^i_j\psi=\pi^i \tilde{X}_u\psi$. 
Since $\pi^i$ is continuous and $u\mapsto \tilde{X}_u\psi$ is measurable 
it follows that $u\mapsto X_u{\,}^i_j\psi$ is measurable. The 
norm bound follows from the fact that $\nm{X_u{\,}^i_j}\leq 
\nm{X_u}$ for all $u\in [0,1]$.
$$
\nm{\int_0^t X_u{\,}^i_j\psi\,du}\leq
\int_0^t \nm{X_u{\,}^i_j\psi}\,du \leq 
\int_0^t\nm{X_u{\,}^i_j}\,du\nm{\psi}\leq t^{\frac1q}\nm{X}_p\nm{\psi}
\leqno{\rm (ii)}
$$
(iii) If $\Gv\in \CH_{00}$ then, interchanging finite sums and 
integrals,
\begin{eqnarray*}
\nm{\sum_{j=0}\iiy Y_t{\,}^i_j \Gv^j}
&=& \nm{\int_0^t\sum_{j=0}\iiy  X_u{\,}^i_j \Gv^j\,du} \\
&=& \nm{\int_0^t (\tilde{X}_u\Gv)^i\,du} \\
&=& \nm{\left(\int_0^t \tilde{X}_u\Gv\,du\right)^i} \\
&\leq & \int_0^t\nm{\tilde{X}_u}\,du\nm{\Gv} \\
&\leq& t^{\frac1q} \nm{X}_p \nm{\Gv}.
\end{eqnarray*}
Thus the series $\sum_{j=0}\iiy Y_t{\,}^i_j \psi^j$
satisfies the Cauchy condition for convergence.
\bb\ni
(iv) If $\psi\in \CH$ and 
$\psi\ub{n}=(\psi^0,\psi^1,\ldots,\psi^n,0,\ldots)$ then 
$$
\nm{\sum_{j=0}^n X_u{\,}^i_j \psi^j}
=\nm{(\tilde{X}_u \psi\ub{n})^i}\leq 
\nm{(\tilde{X}_u}\cdot\nm{ \psi\ub{n}}\leq 
\nm{(\tilde{X}_u}\cdot\nm{ \psi}.
$$
By the dominated convergence theorem
$$
\int_0^t\left(\sum_{j=0}\iiy X_u{\,}^i_j \psi^j\right)\,du=
\sum_{j=0}\iiy\int_0^t  X_u{\,}^i_j \psi^j\,du
= \sum_{j=0}\iiy Y_t{\,}^i_j \psi
$$
converges to $\eta^i$ in $\CH^i$. If $\phi\in \CH$ then
\begin{eqnarray}
\nonumber
\sum_{i=0}\iiy |\ip{\eta^i}{\phi^i}|
&=& 
\sum_{i=0}\iiy \left|\int\ip{X_u{\,}^i_j \psi^j}{\phi^i}\,du\right|
\\
\label{ II 73.01}&=&
\sum_{i=0}\iiy \left|\int\ip{(\tilde{X}_u\psi)^i}{\phi^i}\,du\right|
\\
&\leq &
\nonumber\nm{\psi}\cdot\nm{\phi}\int_0^t \nm{X_u}\,du, 
\end{eqnarray}
where the inequality follows from the dominated convergence 
theorem and the Cauchy--Schwarz inequality
$$
\sum_{i=0}\iiy |\ip{(\tilde{X}_u\psi)^i}{\phi}|\leq 
\nm{\tilde{X}_u\psi}\cdot \nm{\phi}\leq 
\nm{X_u}\cdot\nm{\psi}\cdot\nm{\phi}.
$$
This shows that $\FD(Y_t)=\CH$ and $Y_t$ is bounded with 
$\nm{Y_t}\leq t^{1/q}\nm{X}_p$.
\bb\ni
(v) It follows from (\ref{ II 73.01}) and the dominated 
convergence theorem that
$$
\ip{\tilde{Y}_t\psi}{\phi}=
\int_0^t \ip{\tilde{X}_u \psi}{\phi}\,du 
= 
\ip{\int_0^t \tilde{X}_u \psi\,du}{\phi} \,\,\epf
$$
If $\CK$ is a complex Hilbert space with chaos decomposition 
$\CK=\bigoplus_{j=0}\iiy \CK^j$ the forgoing analysis above may be 
carried out with $\CB(\CH)$, $\CB\cmx(\CH)$,... replaced by 
$\CB(\CH,\CK)$, $\CB\cmx(\CH,\CK)$,.... 
\bb\ni
We shall only consider two chaos decompositions in 
this article: 
\bb\ni 
1. {\it Complex Hilbert sequence space} $\ell^2$. This is the  
direct sum of one-dimensional Hilbert spaces. 
In this case chaos matrices are called {\it scalar matrices}.
\bb\ni
The scalar matrices have a partial ordering: 
$\LGk\prec \Lnu$ if $\LGk^i_j\leq \Lnu^i_j$ for all $i,j$.
\bb\ni 
2. {\it Boson Fock space}. Let $\FH_{t]}{\,}^{j}=L_{\rm sym}^2([0,t]^{j})$  
be the Hilbert space of totally symmetric square-integrable 
functions on $[0,t]^{j}$ and let 
\beqn\label{D 73.01}\FH_{t]}\defeq \FH^{0}\oplus  
\FH_{t]}{\,}^{1}\oplus \cdots\oplus \FH_{t]}{\,}^{j}\oplus \cdots. 
\eeqn 
Then $\FH_t$ is the Boson Fock space over $L^2[0,t]$ 
often denoted $\FF_+(L^2[0,t])$. This is one of several 
equivalent representations of $\FF_+(L^2[0,t])$.
The chaos decomposition (\ref{D 73.01}) is the {\em Fock 
decomposition} and will be used throughout this article.
If $t=1$ put $\Fh=\FH_{1]}$ $\FH^j=\FH_{1]}{}^j$
\bb\ni 
If $g\in L^2[0,1]$ let $e^j(g)=(j!)^{-
1/2}g\otimes\cdots\otimes g\in \FH^j$. It is an easy consequence 
of the Stone-Weierstrass theorem that $\CE^j$, the linear span 
of $\{e^j(g):g\in L^2[0,1]\}$ is dense in $\FH^j$. The {\it 
exponential vector} $e(g)\in \FH$ is defined by 
$$ 
e(g)\defeq 
(1,e^1(g), e^2(g),\cdots, e^j(g),\cdots)^t=(1, g, 
\frac{g^{\otimes2}}{(2!)^{1/2}}, \cdots, \frac{g^{\otimes 
j}}{(j!)^{1/2}},\cdots)^t. 
$$ 
The {\it exponential domain}, $\CE_S$, 
$S\subset L^2[0,1]$ is the linear span of
$\{ e(g):g\in L\} $. When $S=L^2[0,1]$ then we put 
$\CE_S=\CE$. 
\bb\ni
The Bose--Fock spaces $\FH_{(t}=\bigoplus_{j=0}\iiy\FH_{(t}{}^j$, 
with $\FH_{(t}{\,}^{j}=L_{\rm sym}^2((0,1]^{j})$ 
are similarly defined. 
\bb\ni
Both $\FH_{t]}$ and $\FH_{(t}$ are subspaces of $\FH$ and $\FH$ 
can be identified with $\FH_{t]}\otimes\FH_{(t}$.
There is a unitary map $\FU=\FU_t:\FH_{t]}\otimes 
\FH_{(t}\row \FH$ such that, whenever 
$g\in L^2[0,t]$ and $h\in L^2(t,1]$,
$$
\FU(e(g)\otimes e(h))=e(g+h).
$$
For each linear transformation $\hat{T}_{t]}$ in $\FH_{t]}$ whose 
domain contains the exponential domain generated by $L^2[0,t]$ 
the formula $\hat{T}_t=\FU\circ(\hat{T}_{t]}\otimes 
\hat{I}_{(t})\circ\FU\inv$
defines a linear operator in $\FH$, the {\it ampliation } of 
$T_{t]}$.
\bb\ni
In order to describe ampliation in terms of chaos matrices we 
briefly review the construction of general Bose-Fock space over 
a complex Hilbert space $H$.
\bb\ni 
Let $H^0=\BBC$ and let 
$H^j=H\otimes\cdots\otimes H$ the $j$-fold tensor product of $H$ 
with itself. The {\it Fock space } of $H$ is 
$$
\FF(H)=\bigoplus_{j=0}\iiy H^j.
$$
Let $S_j$ be the group of permutations of $\{1,2,\ldots,j\}$. 
There is a unique orthogonal projection $P$ in $\FF(H)$
such that
$$
P(f_1\otimes\cdots\otimes f_j)=(j!)\inv \sum_{\Gs\in S_j} 
f_{\Gs_1}\otimes\cdots\otimes f_{\Gs_j}
$$
whenever $f_1,\ldots,f_j\in H$. The {\it Bose--Fock space} 
$\FF_+(H)$ is defined by
$$
\FF_+(H)=P\FF(H).
$$
If $\CH=\FF_+(H)$ and $\CH^j=PH^j$ then the Fock decomposition
$$
\CH=\CH^0\oplus\CH^1\oplus\cdots\oplus\CH^j\oplus\cdots
$$
is a chaos decomposition of $\CH$. The {\it exponential vector} 
of $h\in h$ is 
$$
e(h)=1\oplus h\oplus \frac{h\ub{2}}{(2!)^{1/2}}\oplus\cdots\oplus 
\frac{h\ub{j}}{(j!)^{1/2}}\oplus,
$$
where $h\ub{j}=h\otimes\cdots\otimes h$, the $j$-fold tensor 
product of $h$ with itself. The set $\{e(h):h\in H\}$ of 
exponential vectors is total and linearly independent. 
If $h_1,\ldots,h_j\in H$ 
then $P(h_1\otimes\cdots\otimes h_j)$ is in the linear span of 
$\{h\ub{j}:h\in H\}$. It follows that $\{h\ub{j}:h\in H\}$ is a total 
subset of $\CH^j$.
\bb\ni
Suppose $H=H_1\oplus H_2$ is the direct sum of orthogonal 
subspaces $H_1$ and $H_2$ and let $\CH_1$ and $\CH_2$ be the 
corresponding Bose--Fock spaces.  There is a unique unitary 
transformation $\FU$ from $\CH_1\otimes \CH_2$ onto $\CH$  such 
that 
\beqn\label{EE 73.00}
\FU(e(h_1)\otimes e(h_2))=e(h_1+h_2)\quad\mbox{ whenever }h_1\in 
H_1 \mbox{ and }h_2\in H_2.
\eeqn
For $i=1,2$ let $\CH_i=\bigoplus_{j=0}\iiy \CH_i{}^j$ be the 
Fock decomposition.
Then $\CH_1{}^{j-r}\otimes 
\CH_2{}^r$ may be regarded as a subspace of $\FF(H)^j$ by
putting 
$$
h_1\ub{j-r}\otimes h_2\ub{r}=
\overbrace{h_1\otimes\cdots\otimes h_1}^{j-r\mbox{ times}}\otimes 
\overbrace{h_2\otimes\cdots\otimes h_2}^{r\mbox{ times}}
$$
whenever $h_i\in H_1$. For $\Ga\in\CH_1{}^{j-r}$ and $\Gb\in 
\CH_2{}^r$ define $\Ga\otsym\Gb=P(\Ga\otimes 
\Gb)$ and put $\CH_1{}^{j-r}\otsym \CH_2{}^r
=P(\CH_1{}^{j-r}\otimes \CH_2{}^r)$. The vectors $\{h_1{}\ub{j-
r}\otsym h_2{}\ub{r}:h_1\in H_1,\,\,h_2\in H_2\}$ are total in 
$\CH_1{}^{j-r}\otsym \CH_2{}^r$. The subspaces 
$\{\CH_1{}^{j-r}\otsym \CH_2{}^r:0\leq r\leq j\}$ 
of $\CH$ are orthogonal.
\begin{Propn}\label{PP 73.01}
$$
\CH^j=\bigoplus_{r=0}^j \left(\CH_1{}^{j-r}\otsym 
\CH_2{}^r\right) 
\leqno{\rm 
(i)}
$$
and is the closed linear span of the vectors
$$
\psi^j=\Ga^j\otsym \Gb^0+\cdots+\Ga^{j-r}\otsym 
\Gb^r\cdots+\Ga^0\otsym \Gb^j,\qquad \Ga^{j-r}\in \CH_1{}^{j-
r},\,\,\Gb^r\in \CH_2{}^r.
$$
$$
\FU(\CH_1{}^{j-r}\otimes \CH_2{}^r)=\CH_1{}^{j-r}\otsym 
\CH_2{}^r. \leqno{\rm (ii)}
$$
$$
\FU(\Ga^{j-r}\otimes \Gb^r)=\choose{j}{r}^{\frac12}\Ga^{j-r}\otsym \Gb^r.
\leqno{\rm (iii)}
$$
\end{Propn}
\pf
If $h_1\in H_1$, $h_2\in H_2$ and $a,b\in \BBR$ then, by 
(\ref{EE 73.00}),
\begin{eqnarray*}
\FU(e(ah_1)\otimes e(bh_2))&=& 
\FU\left(\sum_{j=0}\iiy\sum_{r=0}^j a^{j-r} b^r \frac{h_1\ub{j-
r}\otimes h_2\ub{r}}{((j-r)!r!)^{1/2}}\right)
\\
&=& \sum_{j=0}\iiy \frac{(ah_1+bh_2)\ub{j}}{(j!)^{1/2}}.
\end{eqnarray*}
Since  $P((ah_1+bh_2)\ub{j})=(ah_1+bh_2)\ub{j}$,
\begin{eqnarray}
\label{EE 73.01}\FU \left(\sum_{r=0}^j a^{j-r} b^r \frac{h_1\ub{j-
r}\otimes h_2\ub{r}}{((j-r)!r!)^{1/2}}\right)
&=& \frac{(ah_1+bh_2)\ub{j}}{(j!)^{1/2}} \\
\nonumber &=& \sum_{r=0}^j a^{j-r}b^r \choose{r}{j} h_1\ub{j-r}\otimes 
h_2\ub{r} \\
\label{EE 73.02}&=& \sum_{r=0}^j a^{j-r}b^r \choose{r}{j} h_1\ub{j-r}\otsym
h_2\ub{r} 
\end{eqnarray}
The vectors on the right hand side of Equation (\ref{EE 73.01}) 
are total in $\CH^j$ so that (i) follows from (\ref{EE 73.02}). 
Parts (ii) and (iii) follow from equating the coefficients of 
$a^{j-r} b^r$ in  (\ref{EE 73.01}) and (\ref{EE 73.02}).\,\,\epf
\bb\ni 
A subspace $S$ of $L^2[0,\infty]$ is {\it admissible} if
\bb\ni
(i)  $\chi_{[0,t]}g$ belongs to $S$ for all 
$t\in [0,1]$ whenever $g\in S$; 
\bb\ni
(ii) $S+iS$ is dense in $L^2[0,1]$, 
\bb\ni
Hudson and Parthasarathy \cite{HP1} use for their integrands processes of 
operators whose domains contain $\CE_S$ for some fixed admissible 
subspace $S$ of $L^2[0,1]$. Fortunately, with these domain 
restrictions the correspondence between chaos matrix $T$ 
and linear transformation $\tilde{T}$ is one to one. 
\begin{Lemma}\label{L 73.01} Let $\hat{T}$ be a linear 
transformation in $\CH$ with $\CE_S\subset\CD(\hat{T})$ 
for some admissible subspace $S$ of $L^2[0,1]$. If $T$ is a chaos 
matrix representation of $\hat{T}$ then
for each $g\in S$, the sum 
\beqn\label{E 73.001} \sum_{j=0}\iiy T^i_j 
e^j(g) 
\eeqn 
is absolutely convergent in $\FH^i$ for $i=0,1,2,\ldots$. 
\bb\ni 
Consequently $T$ is the only chaos matrix 
representation of $\hat{T}$. 
\end{Lemma} 
\pf Since $S$ is a 
complex vector space $e(zg)$ is in $\CE_S$ whenever $g\in S$ and 
$z\in \BBC$. Now $e^j(zg)=z^je(g)$ so that the power series 
$\sum_{j=0}\iiy z^j T^i_j e^j(g)$ is convergent, and therefore 
absolutely convergent, in $\FH^i$ for all $z\in \BBC$. 
\bb\ni 
If $S$ is a cmx representation  of $\hat{T}$ then 
$\sum_{j=0}\iiy z^j T^i_j e^j(g)=\sum_{j=0}\iiy z^j S^i_j e^j(g) 
$ for all $z\in \BBC$. Therefore $S^i_j e^j(g)= T^i_j e^j(g)$ for 
all $g\in S$. Since $S$ is admissible $S+iS$ is dense in 
$L^2[0,1]$ and $\CE_S^j=\{e^j(g):g\in S\}$ is total in $\FH^j$. 
The bounded operators $S^i_j$ and $T^i_j$ agree on $\CE_S^j$ and 
therefore $S^i_j=T^i_j$. Therefore $S=T$ and $T$ is the only cmx 
representation of $\hat{T}$.\,\,
\epf 
\bb\ni  
Similar considerations apply to the chaos decompositions of $\FH_{t]}$ 
and $\FH_{(t}$ and the corresponding exponential domains 
$\CE_{t]}$ and $\CE_{(t}$. 

\section{\label{BP}Basic Processes} 
In quantum stochastic analysis the basic integrators 
have particularly simple cmx 
representations: diagonal, subdiagonal and superdiagonal. 
This and Lemma \ref{L 73.01} make feasible a chaos matrix based theory of 
quantum stochastic integration compatible with the Hudson--Parthasarathy 
theory.
\bb\ni
If ${J}\in\CB(L^2(\BBR_+))$ define $\GL(J)^j_j$ in $\CB(\FH^j)$ by 
putting $\GL(J)^0_0=0$ and 
$$ 
({\GL}({J})^j_j\psi)(x_1,\ldots,x_j)=\dyl\sum_{{\nu}=1}^j 
({J}_{\nu}\psi)(x_1,\ldots,x_j), 
$$ 
whenever $\psi\in \FH^j$, 
$j=1,2,\ldots$. In this formula $J_{\nu}$ acts on $\psi$ by 
fixing all variables other than $x_{\nu}$ and then applying $J$. 
The theory of Lebesgue integration shows 
that ${J}_{\nu}$ is in $\CB(\FH)$ with $\nm{{J}_{\nu}}=\nm{{J}}$ for 
$j=1,2,\ldots$ and ${\nu}=1,\ldots, j$. Therefore ${\GL}({J})^j_j$ is well 
defined with $\|{\GL}({J})^j_j\|=j\nm{J}$.
\bb\ni
The second quantisation of $J$ \cite[\S 5.2.1]{BR2} 
is represented by the chaos matrix 
$$
\GL(J)= {\rm diag}[0,\GL(J)^1_1,\GL(J)^2_2,\GL(J)^3_3,\ldots].
$$
If $J$ is self-adjoint then 
$\GL(J)^j_j$ is self-adjoint for all $j\in \BBN$ and $\GL(J)\st=\GL(J)$.
\bb\ni
If $f\in L^2[0,1]$ define the adjoint pair, 
$a(f)^j_{j+1}$ in $\CB(\FH^{j+1},\FH^j)$ and $a\d(f)^{j+1}_j$
in $\CB(\FH^{j},\FH^{j+1})$,  by putting
\begin{eqnarray*}
({a}(f)^j_{j+1}\psi)(x_1,\ldots,x_j)&=&
(j+1)^{1/2}\dyl\int_0^{1} \ol{f(s)}\psi(x_1,\ldots,x_j,s)\,ds,\\
({a}\d(f)^{j+1}_j\Gv)(x_1,\ldots,x_{j+1})&=&
(j+1)^{-1/2}\dyl\sum_{{\nu}=1}^{j+1} 
f(x_{\nu})\Gv(x_1,\ldots,\hat{x}_{\nu},\ldots,x_{j+1}),
\end{eqnarray*}
whenever, $j=1,2,\ldots$, $\psi\in \FH^{j+1}$ and $\Gv\in \FH^j$. The `hat' 
over ${x}_{\nu}$ in the above formula suppresses $x_{\nu}$. Then 
\beqn
\nm{a(f)^j_{j+1}}=(j+1)^{1/2}\nm{f}, \qquad 
\nm{a\d(f)^{j+1}_{j}}= (j+1)^{1/2}\nm{f},
\eeqn
and the chaos matrices 
\begin{eqnarray*}
\mbox{\footnotesize$ \left[\ba{ccccc}
0&a(f)^0_1&0&0&\ldots \\
0&0&a(f)^1_2&0&\ldots \\
0&0&0&a(f)^2_3&\ldots \\
\vdots&\vdots&\vdots&\vdots&\ddots 
\ea \right]$},
&\quad&
\mbox{\footnotesize$ \left[\ba{ccccc}
0&0&0&0&\ldots\\
a\d(f)^1_0&0&0&0&\ldots\\
0&a\d(f)^2_1&0&0&\ldots\\
0&0&a\d(f)^3_2&0&\ldots\\
\vdots&\vdots&\vdots&\vdots&\ddots
\ea \right]$}
\end{eqnarray*}
represent, respectively, the annihilation and creation operators 
${a}(f)$ and ${a}\d(f)$ on $\CD(N^{1/2})$ \cite[\S 5.2.1]{BR2}. 
The same symbols ${a}(f)$ and ${a}\d(f)$ denote the matrices 
above and the operators they represent.
\bb\ni
If $f$ is in $L\iiy[0,1]$ and $j\in \BBN\cup\{0\}$ then,
for $\psi\in \FH^{j+1}$ and $\Gv\in L^2([0,1],\FH^j)$,
the formulae
\begin{eqnarray*}
\nabla(f)^j_{j+1}(\psi)(s)(x_1,\ldots,x_j) 
&=& (j+1)^{\frac12}\psi(x_1,\ldots,x_j,s)\ol{f(s)},
\\
\CS(f)^{j+1}_j(\Gv)(x_1,\ldots,x_{j+1}) 
&=& (j+1)^{-\frac12}
\sum_{k=1}^{j+1}f(x_k)\Gv(x_k)(x_1,\ldots,\hat{x}_k,\ldots,x_{j+1}),
\end{eqnarray*}
define an adjoint pair of bounded linear transformations
$$
\nabla(f)^j_{j+1}:\FH^{j+1}\row L^2([0,1], \FH^j),\qquad 
\CS(f)^{j+1}_j:L^2([0,1], \FH^j)\row \FH^{j+1}.
$$
A simple calculation shows that 
\beqn
\label{E 73.211}
\hspace{-0.2in}\nm{\nabla(f)^j_{j+1}}=(j+1)^{1/2}\nm{f}\iy,& \qquad &
\nm{\CS(f)^{j+1}_{j}}= (j+1)^{1/2}\nm{f}\iy.
\eeqn
The chaos matrices 
\beqn\label{M 73.01}
\hspace{-0.2in}\mbox{\footnotesize $\left[\ba{ccccc}
0&\nabla(f)^0_1&0&0&\ldots \\
0&0&\nabla(f)^1_2&0&\ldots \\
0&0&0&\nabla(f)^2_3&\ldots \\
\vdots&\vdots&\vdots&\vdots&\ddots 
\ea \right],$}
&\,\,&
\mbox{\footnotesize$\left[\ba{ccccc}
0&0&0&0&\ldots\\
\CS(f)^1_0&0&0&0&\ldots\\
0&\CS(f)^2_1&0&0&\ldots\\
0&0&\CS(f)^3_2&0&\ldots\\
\vdots&\vdots&\vdots&\vdots&\ddots
\ea \right]$,}
\eeqn 
represent ${\nabla}(f)$, the {\it gradient 
operator}, and ${\CS}(f)$, the {\it Hitsuda--Skorohod} operator, on 
$\CD(N^{1/2})$ \cite[\S 1.]{Lin1}.
We dispense with the " $\hat{}$ " denoting operator 
and let ${\nabla}(f)$ and ${\CS}(f)$ stand for the matrices and for the 
operators they represent.
Similarly we write $N$ for the number operator.
\bb\ni
It follows from (\ref{E 73.211}) that 
$( {N}+1)^{-1/2}  {\nabla}(f)$ and 
$ {\CS}(f)( {N}+1)^{-1/2}$ are bounded in 
norm by $\nm{f}\iy$.
\bb\ni
If $t\in [0,1]$ then 
$$
\nabla_t\defeq \nabla(\chi_{[0,t]}),
\qquad\CS_t\defeq \CS(\chi_{[0,t]}
$$
We will sometimes write $\nabla(f)(s)\Gv$ 
for $(\nabla(f)\Gv)(s)$ and $\nabla_t(s)\Gv$ for $\nabla_t(\Gv)(s)$. 
\bb\ni
The operators $ {\nabla}(f)$ and $ {\CS}(f)$ are formally 
adjoint \cite{Lin1}. This is reflected in (\ref{E 73.101}) 
below. The identity (\ref{E 73.103}) is the key to the 
quantum Ito formula.
\begin{Propn}\label{P 73.211} Let $f,g\in L^2[0,1]$. 
If $\psi,\Gth\in L^2([0,1],\FH^j)$, $\phi\in L^2
([0,1],\FH^{j+1})$  and $\Gv\in \FH^{j+1}$ then
\begin{eqnarray}
\label{E 73.101}
\ip{\CS(f)\psi}{\Gv}&=& \int_0^{1} 
\ip{\psi(s)}{\nabla(f)(s)\Gv}\,ds=\ip{\psi}{\nabla(f)\Gv};
\\
\label{E 73.102}
\int_0^{1} \ip{\CS(f)\psi}{\phi(u)}\,du&=& \int_0^{1} \int_0^{1} 
\ip{\psi(s)}{\nabla(f)(s)\phi(u)}\,du\,ds;
\\
\nonumber
\ip{\CS(f)\psi}{\CS(g)\Gth} 
&=& \int_0^{1} \int_0^{1} 
\ip{\nabla(\ol{f})(u)\psi(s)}{\nabla(\ol{g})(s)\Gth(u)}\,du\,ds
\\
\label{E 73.103}
&&\hspace{1in}+ \int_0^{1}\ip{f(s)\psi(s)}{g(s)\Gth(s)}.
\end{eqnarray}
\end{Propn}
\pf For each $s\in [0,1]$ the functions $\psi(s)$ and $\Gv$ 
are completely symmetric. Therefore $\ip{\CS(f)\psi}{\Gv}$ equals
\begin{eqnarray*}
&& \hspace{-.5in}(j+1)^{-\frac12}\sum_{k=1}^{j+1}\int_0^{1} 
f(x_k)\psi(x_k)(x_1,\ldots,\hat{x}_k,\ldots,x_{j+1})
\ol{\Gv(x_1,\ldots,x_{j+1})}
\,dx
\\
&=& (j+1)^{\frac12}\int_0^{1}\int_{[0,1]^{j}}\psi(s)(y_1,\ldots,y_{j})
\ol{\Gv(y_1,\ldots,y_{j},s)\ol{f(s)}}\,dy\,ds
\\
&=& \int_0^{1} \ip{\psi(s)}{\nabla(f)(s)\Gv}\,ds,
\end{eqnarray*}
where $dx=dx_1\ldots dx_{j+1}$ and $dy=dy_1\dots dy_j$, proving (\ref{E 73.101}). This leads immediately to (\ref{E 73.102}). To prove 
(\ref{E 73.103}) note that $(j+1)\ip{\CS(f)\psi}{\CS(g)\Gth}$ is equal to
{\small
\begin{eqnarray*}
&&
\hspace{-0.5in}\sum_{k,l=1}^{j+1}\int_{[0,{1}]^{j+1}}
\!\!\!\!\!\!f(x_k)\psi(x_k)(x_1,\ldots,\hat{x}_k,\ldots,x_{j+1})
\ol{g(x_l)\Gth(x_l)(x_1,\ldots,\hat{x}_l,\ldots,x_{j+1})}
\,dx
\\
&=&
\sum_{k\neq l=1}^{j+1}\int_{[0,{1}]^{j+1}}
\!\!\!\!\!\!f(x_k)\psi(x_k)(x_1,\ldots,\hat{x}_k,\ldots,x_{j+1})
\ol{g(x_l)\Gth(x_l)(x_1,\ldots,\hat{x}_l,\ldots,x_{j+1})}\,dx
\\
&&+ 
\sum_{k=1}^{j+1}\int_{[0,{1}]^{j+1}}
\!\!f(x_k)\psi(x_k)(x_1,\ldots,\hat{x}_k,\ldots,x_{j+1})
\ol{g(x_k)\Gth(x_k)(x_1,\ldots,\hat{x}_k,\ldots,x_{j+1})}\,dx
\\
&=& (j+1)\!\!\int_0^{1}\!\! \int_0^{1}\!\! 
\ip{\nabla(\ol{f})(u)\psi(s)}{\nabla(\ol{g})(s)\Gth(u)}\,du\,ds+
(j+1)\!\!\int_0^{1}\!\!\ip{f(s)\psi(s)}{g(s)\Gth(s)}.\,\,
\epf
\end{eqnarray*} }
The proof of (\ref{E 73.102}) yields the following corollary.
\begin{Cor}
If $f\in L\iiy[0,1]$ and $j\in \BBN\cup \{0\}$ then
$$
\left(\nabla(f)^j_{j+1}\right)\st=\CS(f)^{j+1}_j.
$$
\end{Cor}                                                  
If $\Psi$ is one of the basic matrices $\GL(f)$, $a(f)$, 
$a\d(f)$, 
$\nabla(f)$ and $\CS(f)$ then $\Psi$ is tridiagonal. Therefore, 
the product chaos matrices $\Psi T$ and $T\Psi$ exist for each 
chaos matrix $T$. For example the formulae 
\begin{eqnarray*}
((\nabla(f)T)^{i-1}_j \psi)(s)(x_1,\ldots,x_{i-1})&=&
(\nabla(f)^{i-1}_i(s)T^i_j\psi)(x_1,\ldots,x_{i-1}); 
\\
((T\nabla(f))^i_j \psi)(s)(x_1,\ldots,x_i)&=& (T^i_{j-1}
(\nabla(f)^{j-1}_{j}\psi)(s))(x_1,\ldots,x_i),
\end{eqnarray*}
whenever $\psi\in \FH^j$ 
define the products (composition) of the chaos matrices $\nabla(f) T$ and 
$T\nabla(f)$. The other products are similarly defined. 

\section{\label{S 3} Adaptedness}
Let $\hat{T}=\{\hat{T}_t:t\in [0,1]\}$ be a process of 
operators in $\CL(\CH)$ which is adapted in the sense of 
Hudson and Parthasarathy \cite[Definition 3.1]{HP1}.
Then  $\hat{T}_t$ is the ampliation of 
an operator $\hat{T}_{t]}$ in $\CL(\FH_{t]})$ for all $t$.
A possible definition of adaptedness for a cmx process $T$ 
would be that the corresponding operator-valued process $\tilde{T}$ 
be adapted. Our proof of the quantum Ito formula for cmx 
semimartingales requires a more intrinsic definition of adaptedness 
and, in particular, a description of ampliation in terms of 
chaos matrices. This intrinsic definition is shown to be compatible with 
the Hudson-Parthasarathy definition.
\bb\ni
For each $t\in [0,1]$ there is a unique unitary
$\FU=\FU_t:\FH_{t]}\otimes \FH_{(t}\hraa \FH$ satisfying 
\beqn\label{E 65.01}
\FU (e(g_{t]})\otimes e(g_{(t}))=e(g_{t]}+g_{(t})   
\qquad\forall \quad g\in L^2[0,1] 
\eeqn
where $g_{t]}=\chi_{[0,t]}g$ and $g_{(t}=\chi_{(t,1]}g$. The 
map $\FU_t$ is natural and defined at the algebraic level:
$$
\CE=\FU_t(\CE_{t]}\ul{\otimes}\CE_{(t})
$$
where $\CE$, 
$\CE_{t]}$ and $\CE_{(t}$ are the exponential domains of $\FH$, $\FH_{t]}$ 
and $\FH_{(t}$ respectively and `$\ul{\otimes}$' denotes 
{\em algebraic} tensor product.
\bb\ni
For the rest of this section fix $t\in [0,1]$ and put 
$\FU=\FU_t$. Let $\hat{T}_{t]}$ be a linear transformation in $\FH_{t]}$ 
with domain $\CD(\hat{T}_{t]})$. 
{\it The ampliation} of $\hat{T}_{t]}$ to $\FU[\CD(T_{t]})\ul{\otimes} 
\FH_{(t}]$ is the linear transformation
$\hat{T}_t\defeq\FU\circ(\hat{T}_{t]}\ul{\otimes} 
\hat{I}_{(t})\circ \FU\inv$. Thus for $\Ga\in \CD(\hat{T}_{t]})$ 
and $\Gb\in \FH_{(t}$
$$
\hat{T}\FU[\Ga\otimes \Gb]=\FU[(\hat{T}_{t]}\Ga)\otimes\Gb].
$$
Any extension of $\hat{T}_t$  is said to be 
{\it an ampliation} of $\hat{T}_{t]}$.
\bb\ni
If $\hat{T}_{t]}$ is bounded with cmx representation $T_{t]}$ 
there is a unique bounded chaos matrix $T_t$ such that 
$\tilde{T}_t$ is an ampliation of $\hat{T}_{t]}$ to 
$\FH$. Therefore the entries of 
$T_t$ may be represented in terms of the entries of 
$T_{t]}$. To find this representation we use the deconstruction 
of $\FU$ described in Proposition \ref{PP 73.01}.
\bb\ni
It follows from  Proposition \ref{PP 73.01}
that if $\Ga\in \FH_{t]}^{j-r}$ and $\Gb\in \FH_{(t}^r$ then 
\beqn\label{E 75.879}
\FU(\Ga\otimes \Gb)=
\left(\ba{c} j\\ r\ea\right)^{\frac12}\Ga\otsym \Gb.
\eeqn
Moreover $\FH^j$ is the closed 
linear span of all vectors of the form 
\beqn\label{E 65.0001}
\psi^j&=&\Ga^j\otsym  \Gb^0\oplus\cdots\oplus 
\Ga^{j-r}\otsym  \Gb^r\oplus\cdots\oplus  \Ga^0\otsym  \Gb^j,
\eeqn    
where $\Ga^{j-r}\in \FH\lrb{t}^{j-r}$ and $\Gb^r \in \FH^{r}_{(t}$. 
\bb\ni
If $\Ga,\Gg\in \FH_{t]}^{i-r}$ and $\Gb,\Gd \in \FH^r_{(t}$ then 
$$
\nm{\Ga\otimes\Gb}=\nm{\FU(\Ga\otimes\Gb)}
=\choose{i}{r}^{\frac12}\nm{\Ga\otsym  \Gb},
$$
and, by polarisation,
$$
\choose{i}{r}\ip{\Ga\otsym \Gb}{\Gg\otsym \Gd}=
\ip{\Ga}{\Gg}\ip{\Gb}{\Gd}.
$$
This implies that if $S\in \CB(\FH_{t]}^{j-r},\FH_{t]}^{i-r})$ then 
\beqn\label{E 75.875}
\choose{i}{r}\ip{(S\Ga)\otsym \Gb}{\Gg\otsym \Gd}&=&
\choose{j}{r}\ip{\Ga\otsym \Gb}{(S\st\Gg)\otsym \Gd},
\eeqn
whenever $\Ga\in \FH_{t]}^{j-r}$, $\Gg\in \FH_{t]}^{i-r}$ and 
$\Gb,\Gd \in \FH^r_{(t}$. 
\bb\ni
Define $S\otsym  I_{(t}^r$ 
in $\CB(\FH_{t]}^{j-r}\otsym \FH_{(t}^r,\FH_{t]}^{i-r}\otsym 
\FH_{(t}^r)$, where $I_{(t}^r$ is the identity in $\FH_{(t}^r$, by putting 
$$
(S\otsym  I_{(t}^r)(\Ga \otsym  \Gb)=(S\Ga)\otsym  \Gb
$$
whenever $\Ga\in \FH_{t]}^{j-r}$ and $\Gb \in \FH^r_{(t}$. It follows from
(\ref{E 75.879}) that
\beqn\label{E 75.464}
S\otsym  I_{(t}^r&=&\choose{i}{r}^{-\frac12} \choose{j}{r}^{\frac12}
\FU(S\otimes I_{(t}^r)\FU\inv.
\eeqn
We obtain following lemma directly from (\ref{E 75.875}). 
\begin{Lemma}
If $S\in \CB(\FH_{t]}^{j-r},\FH_{t]}^{i-r})$ then 
\beqn
\label{E 75.215}
\choose{i}{r}(S\otsym I_{(t}^r)\st=\choose{j}{r}(S\st\otsym  I_{(t}^r).
\eeqn
\end{Lemma}
\begin{Thm}\label{P 1.001} Let $T_{t]}$ be a chaos matrix  for 
$\FH_{t]}$  and let $T_t$ be the chaos matrix whose entries are 
defined by the formula
\beqn\label{E 1.1}
T_t{\,}^i_j=\bigoplus_{r=0}^{i\wedge j} 
\choose{i}{r}^{\frac12} \choose{j}{r}^{-\frac12}
(T_{t]}{\,}^{i-r}_{j-r} \otsym  I_{(t}^r) 
=\bigoplus_{r=0}^{i\wedge j} 
\FU\circ(T_{t]}{\,}^{i-r}_{j-r} \otimes I_{(t}^r)\circ\FU\inv.
\eeqn
Then 
\bb\ni
{\rm (i)} $\FU[\FD(T_{t]})\ul{\otimes} \FH_{(t}]\subset \FD(T_t)$;
\bb\ni
{\rm (ii)}
$\tilde{T}_t$ is an ampliation of $\tilde{T}_{t]}$ on  
$\FU[\FD(T_{t]})\ul{\otimes} \FH_{(t}]$.
\end{Thm}
\pf 
If $\Ga\in \FD(T_{t]})$ and $\Gb\in \FH_{(t}$ then 
$\eta^i=(\tilde{T}_{t]}\Ga)^i=\sum_{j=0}\iiy 
T_{t]}{\,}^i_j\Ga^j$ where the sum is convergent in 
$\FH^i$. Moreover $\eta=(\eta^0,\eta^1,\ldots,\eta^i,\ldots)^t\in 
\FH$. Therefore
\begin{eqnarray*}
\left((\tilde{T}_{t]}\Ga)\otimes\Gb\right)^i 
&=& \bigoplus_{r=0}^i(T_{t]}\Ga)^{i-r}\otimes \Gb^r  \\
&=& \bigoplus_{r=0}^i\sum_{j=0}\iiy(T_{t]}{\,}^{i-r}_j\Ga^j)\otimes \Gb^r  \\
&=& \bigoplus_{r=0}^i\sum_{j=r}\iiy(T_{t]}{\,}^{i-r}_{j-r}\Ga^{j-
r})\otimes \Gb^r  \\
&=& \sum_{j=0}\iiy\bigoplus_{r=0}^{i\wedge j}(T_{t]}{\,}^{i-r}_{j-r}\Ga^{j-
r})\otimes \Gb^r  \\
&=& \sum_{j=0}\iiy\bigoplus_{r=0}^{i\wedge j}(T_{t]}{\,}^{i-
r}_{j-r}\otimes I_{(t}^r)(\Ga\otimes \Gb).
\end{eqnarray*}
These sums are all convergent in $\FH_{t]}\otimes \FH_{(t}$. Therefore 
\begin{eqnarray}
\nonumber \left(\FU[(\tilde{T}_{t]}\Ga)\otimes \Gb]\right)^i
&=& \FU([(\tilde{T}_{t]}\Ga)\otimes \Gb]^i) \\
\nonumber &=& \sum_{j=0}\iiy\bigoplus_{r=0}^{i\wedge j}
\FU\left[(T_{t]}{\,}^{i-r}_{j-r}\otimes I_{(t}^r)
(\Ga\otimes \Gb)\right] \\
\label{E 75.22} &=& \sum_{j=0}\iiy\bigoplus_{r=0}^{i\wedge 
j}\choose{i}{r}^{\frac12}\!\choose{j}{r}^{-\frac12}\!\!(T_{t]}{\,}^{i-
r}_{j-r}\otsym\!  I_{(t}^r)\FU[(\Ga\otimes \Gb)] \\
\label{E 75.23}&=& \sum_{j=0}\iiy T_t{\,}^i_j 
(\FU[\Ga\otimes\Gb])^j.
\end{eqnarray}
Equation (\ref{E 75.22}) follows from (\ref{E 75.464}) and 
the sums are all convergent in $\FH^i$. Now
\begin{eqnarray*}
\bigoplus_{i=0}\iiy \left((\tilde{T}_{t]}\Ga)\otimes \Gb\right)^i
&=&
\bigoplus_{i=0}\iiy\bigoplus_{r=0}^i (\tilde{T}_{t]}\Ga)^{i-r}\otimes 
\Gb^r \\
&=& 
\bigoplus_{r=0}\iiy\bigoplus_{i=r}\iiy (\tilde{T}_{t]}\Ga)^{i-r}\otimes 
\Gb^r \\
&=&
\bigoplus_{r=0}\iiy\bigoplus_{i=0}\iiy (\tilde{T}_{t]}\Ga)^i\otimes 
\Gb^r
\end{eqnarray*}
which is convergent in $\FH_{t]}\ul{\otimes} \FH_{(t}$ to 
$(\tilde{T}_{t]}\Ga)\otimes \Gb$. Therefore 
$$
\bigoplus_{i=0}\iiy \FU\left[\left((\tilde{T}_{t]}\Ga)\otimes 
\Gb\right)^i\right]
$$
is convergent in $\FH$. It follows from (\ref{E 75.23}) that 
$\FU[\Ga\otimes \Gb]$ belongs to $\FD(T_t)$, proving (i), and also that 
$$
\FU\left[(\tilde{T}_{t]}\Ga)\otimes \Gb\right]= 
\tilde{T}_t\FU[\Ga\otimes \Gb].
$$
Therefore $\tilde{T}_t$ agrees with  the ampliation of $T_{t]}$ 
on $\FD(T_{t]})\ul{\otimes} \FH_{(t}$, proving (ii). \,\,\epf
\bb\ni
We define the {\it ampliation} of $T_{t]}$ to be 
chaos matrix $T_t$ defined by formula (\ref{E 1.1}). 
In this case we write
$$
T_t=T_{t]}\otsym  I_{(t}.
$$
Note that $T_t$ is $(2k+1)$-diagonal whenever $T_{t]}$ is 
$(2k+1)$-diagonal.
\bb\ni
A cmx process $T$ is {\it adapted} if for each $t\in [0,1]$ there exists a 
chaos matrix $T_{t]}$ such that $T_t$ is the ampliation of 
$T_{t]}$ for each $t\in [0,1]$.
\bb\ni
It was shown in Theorem \ref{T 73.001} (ii) that the map $T\mapsto 
\tilde{T}$ is an isometric isomorphism from $\CB_{\rm 
cmx}(\FH_{t]})$ onto $\CB(\FH_{t]})$ for each $t\in [0,1]$. This 
leads to the following corollaries.
\begin{Cor} Let $T_t$ be a bounded chaos matrix. Then $T_t$ is 
the ampliation of $T_{t]}\in \CB_{\rm cmx}(\FH_{t]})$  if and 
only if $\tilde{T}_t$ is the ampliation of $\tilde{T}_{t]}$.
\end{Cor}
\begin{Cor}\label{C 074.01}
$T=\{T_t:t\in [0,1]\}$ is a regular adapted cmx process 
if and only $\tilde{T}=\{\tilde{T}_t:t\in [0,1]\}$ is a regular adapted 
operator process.
\end{Cor}
\bb\ni
Using the gradient operator $\nabla\defeq \nabla(\chi_{[0,1]})$ 
we give a characterisation of ampliation 
which is significantly easier to use in computation than (\ref{E 1.1}).
\bb\ni
In what follows $\nabla(s)(\psi)$ denotes $(\nabla(\psi))(s)$ 
whenever $\psi\in \CD(N^{1/2})$.
\bb\ni
Let $\CM_{j-r}^j$, $j\geq r$ be the set of strictly monotone 
functions $\pi:\{1,2,\ldots,j-r\}\row \{1,2,\ldots,j\}$ and if 
$\pi\in \CM_{j-r}^j$ let $\check{\pi}$ be the unique function in 
$C_r^j$ whose image is disjoint from that of $\pi$. Note that 
$|\CM^j_r|=|\CM_{j-r}^j|=\choose{j}{r}$.
\begin{Lemma}
If $\Ga\in \FH_{t]}^{j-r}$, $\Gb\in \FH^r_{(t}$ and $s_1,\ldots,s_r\in 
(t,1]$ then
\begin{eqnarray}
\label{E 75.621}\nabla(s_1)(\Ga\otsym \Gb)
&=& \frac{r^{\frac12}}{j^{\frac12}}\Ga\otsym \nabla(s_1)\Gb;
\\
\label{E 75.622}
\nabla(s_1)\cdots\nabla(s_r)(\Ga\otsym \Gb)&=&
\choose{j}{r}^{-\frac12}(\Ga\otsym 1)\Gb(s_1,\ldots,s_r),
\end{eqnarray}
\bb\ni
where $1$ is the unit in $\FH_{(t}^0=\FH^0=\BBC$.
\end{Lemma}
\pf
\begin{eqnarray*}
\nabla(s_1) (\Ga\otsym  
\Gb)(x_1,\ldots,x_{j-1})&&
\\
&&\hspace{-1.2in}=
j^{\frac12}\left(\!\ba{c} j\\r\ea\!\right)\inv \!
 \sum_{\stackrel{\pi\in \CM^j_{j-r}}{\check{\pi}_r=j}}
\Ga(x_{\pi_1}\cdots x_{\pi_{j-r}}) 
\Gb(x_{\check{\pi}_1}\cdots x_{\check{\pi}_{r-1}},s_1)
\\
&&\hspace{-1.2in} =
j^{\frac12}\left(\!\ba{c} j\\r\ea\!\right)\inv \!
 \sum_{\pi\in \CM^{j-1}_{j-r}}
\Ga(x_{\pi_1}\cdots x_{\pi_{j-r}}) 
\Gb(x_{\check{\pi}_1}\cdots x_{\check{\pi}_{r-1}},s_1)
\\
&&\hspace{-1.2in} =
\frac{j^{\frac12}}{r^{\frac12}}\left(\!\ba{c} j\\r\ea\!\right)\inv \!
 \sum_{\pi\in \CM^{j-1}_{j-r}}
\Ga(x_{\pi_1}\cdots x_{\pi_{j-r}}) 
(\nabla(s_1)\Gb)(x_{\check{\pi}_1}\cdots x_{\check{\pi}_{r-1}})
\\
&&\hspace{-1.2in} =
\frac{j^{\frac12}}{r^{\frac12}}\left(\!\ba{c} j\\r\ea\!\right)\inv \!
\choose{j-1}{r-1}(\Ga\otsym \nabla(s_1)\Gb)(x_1,\ldots,x_{j-1})
\\
&&\hspace{-1.2in} = \frac{r^{\frac12}}{j^{\frac12}}
(\Ga\otsym \nabla(s_1)\Gb)(x_1,\ldots,x_{j-1}),
\end{eqnarray*}
proving (\ref{E 75.621}). Repeated applications of (\ref{E 75.621}) give
(\ref{E 75.622}).\,\,\epf
\begin{Thm}\label{T 75.221} A chaos matrix $T_t$ for $\FH$ is the 
ampliation of a chaos matrix $T_{t]}$ for 
$\FH_{t]}$ if and only if, for each $s\in (t,1]$,
\beqn\label{E 65.004}
{\rm (i)}\,\,\,T_t\nabla(s)=\nabla(s) T_t\quad\mbox{ and }\quad 
{\rm (ii)}\,\,\, T_t\st\nabla(s)=\nabla(s) T_t\st.
\eeqn
In this case, if $i,j\in \BBN$, then 
\beqn\label{E 75.213}
T_{t]}{}^i_j=\chi_{[0,t]^j}T_t{\,}^i_j\chi_{[0,t]^i}.
\eeqn
\end{Thm}
\pf If $T_t$ is the ampliation of $T_{t]}$ and 
$\Ga^{j-r}\in \FH_{t]}^{j-r}$ and $\Gb^r \in \FH_{(t}^r$ then
\begin{eqnarray*}
\nabla(s) T_t{\,}^i_j (\Ga^{j-r}\otsym \Gb^r)
&=&
\choose{i}{r}^{\frac12}\choose{j}{r}^{-\frac12}\nabla(s)(T_{t]}{\,}^{i-r}_{j-
r}\Ga^{j-r}\otsym \Gb^r)
\\
&&\hspace{-0.5in} =
\left(\frac{r}{i}\right)^{\frac12}\choose{i}{r}^{\frac12}\choose{j}{r}^{-\frac12}
(T_{t]}{\,}^{i-r}_{j-r}\Ga^{j-r}\otsym \nabla(s)\Gb^r)
\\
&&\hspace{-0.5in} =
\left(\frac{r}{j}\right)^{\frac12}
\choose{i-1}{r-1}^{\frac12}\choose{j-1}{r-1}^{-\frac12}
(T_{t]}{\,}^{i-r}_{j-r}\Ga^{j-r}\otsym \nabla(s)\Gb^r)
\\
&&\hspace{-0.5in} =
\left(\frac{r}{j}\right)^{\frac12}
T_{t}{\,}^{i-1}_{j-1}(\Ga^{j-r}\otsym \nabla(s)\Gb^r)
\\
&&\hspace{-0.5in} =
T_{t}{\,}^{i-1}_{j-1}\nabla(s)(\Ga^{j-r}\otsym \Gb^r).
\end{eqnarray*}
It follows by linearity that $ T_t{\,}^{i-1}_{j-1}\nabla(s) \psi^j
=\nabla(s) T_t{\,}^i_j\psi^j $ almost everywhere for all $\psi^j$ 
of the form (\ref{E 65.0001}). Since such $\psi^j$ are total in 
$\FH^j$ the identity (\ref{E 65.004})(i) is valid.
\bb\ni
It follows from (\ref{E 75.464}) and definition (\ref{E 1.1}) 
that $T_t\st$ is the ampliation of $T_{t]}\st$. Therefore  
the identity (\ref{E 65.004})(ii) 
is also valid.
\bb\ni
Suppose conversely that $T_t$ satisfies the identities 
(\ref{E 65.004}) (i) and (ii) and define $T_{t]}$ by (\ref{E 75.213}).
If $\Ga^{j-r}\in \FH_{t]}^{j-r}$ and $\Gb^r \in \FH_{(t}^r$ and $i\geq j$ 
then
\begin{eqnarray*}
\nabla(s_1)\cdots\nabla(s_r)T_t{\,}^i_j(\Ga\otsym \Gb)
&=&
T_t{\,}^{i-r}_{j-r}\nabla(s_1)\cdots\nabla(s_r)(\Ga\otsym \Gb)
\\
&&\hspace{-0.5in} =
\choose{j}{r}^{-\frac12}T_t{\,}^{i-r}_{j-r}(\Ga\otsym 1)\Gb(s_1,\ldots,s_r),
\\
&&\hspace{-0.5in} =
\choose{j}{r}^{-\frac12}(T_{t]}{}^{i-r}_{j-r}\Ga\otsym 1)
\Gb(s_1,\ldots,s_r),    
\\
&&\hspace{-0.5in} = \choose{i}{r}^{\frac12}\choose{j}{r}^{-\frac12}
\nabla(s_1)\cdots\nabla(s_r)
(T_{t]}{}^{i-r}_{j-r}\Ga\otsym \Gb).
\end{eqnarray*}
The third equality is valid since 
$\nabla(s)T_t{\,}^{i-r}_{j-r}(\Ga\otsym  1)=0$ whenever $s\in (t,1]$. 
Therefore \beqn\label{E 75.214}
T_t{\,}^i_j=\sum_{r=0}^j \choose{i}{r}^{\frac12}\choose{j}{r}^{-\frac12}
(T_{t]}{}^{i-r}_{j-r}\otsym  I_{(t}^r)
\eeqn
whenever $i\geq j$. Now put $S_t=T_t{\,}\st$, apply (\ref{E 75.214}) 
to $S_t$, take adjoints and use (\ref{E 75.215}) to give
\begin{eqnarray*}
T_t{\,}^j_i=(S^i_j)\st  &=& \sum_{r=0}^{i}\choose{i}{r}^{\frac12}
\choose{j}{r}^{-\frac12}(S_{t]}{}^{i-r}_{j-r}\otsym  I_{(t}^r)\st
\\
&=& \sum_{r=0}^{i}\choose{j}{r}^{\frac12}
\choose{i}{r}^{-\frac12}((S_{t]}{}^{i-r}_{j-r})\st\otsym  I_{(t}^r)
\\
&=& \sum_{r=0}^{i}\choose{j}{r}^{\frac12}
\choose{i}{r}^{-\frac12}(S\st_{t]}{}^{j-r}_{i-r}\otsym  I_{(t}^r)
\\
&=& \sum_{r=0}^{i}\choose{j}{r}^{\frac12}
\choose{i}{r}^{-\frac12}(T_{t]}{}^{j-r}_{i-r}\otsym  I_{(t}^r).
\end{eqnarray*}
Thus Equation (\ref{E 75.214}) is valid for all $i,j$ and $T_t{\,}$ 
is the ampliation of $T_{t]}$ defined by Formula (\ref{E 1.1}).\,\,\epf
\bb\ni
The characterisation (\ref{E 65.004}) does not transfer directly to
operators $\tilde{T}$ since there is no guarantee 
that the product operators $\nabla(s) \tilde{T}_t$ exist on the domain 
of $\nabla(s)$ even when $\tilde{T}_t$ is bounded. At the formal 
level it is closely related to the extended definition of adaptedness 
given by Attal \cite[\S I.1.2]{Att1}.
\bb\ni
It is also related to \cite[Definition 4.3.]{OB2}
which appears to be an integrated version of (\ref{E 65.004})
\begin{Cor}\label{C 74.01}
If $S$ and $T$ are adapted cmx processes then 
\bb\ni
{\rm(i)} $T\st$ is adapted; 
\bb\ni
{\rm (ii)} the product process $ST$ is adapted whenever it exists.
\end{Cor}
\pf (i) It is immediate from (\ref{E 65.004}) that $T$ is adapted if and only 
if $T\st$ is adapted.
\bb\ni
(ii) The product process $ST$ exists if and only if $(ST)\st$ exists. 
Since $S$ and $T$ are adapted
$$
S_t T_t\nabla(s)=S_t\nabla(s) T_t=\nabla(s)S_t T_t,
$$
and 
$$ 
T_t\st S_t\st \nabla(s)=T_t\st \nabla(s)S_t\st = \nabla(s)T_t\st S_t\st, 
$$
so that $(ST)_t \nabla(s)=\nabla(s)(ST)_t$ and 
$(ST)\st_t \nabla(s)=\nabla(s)(ST)\st_t$ whenever $s>t$. Therefore
$(ST)_t$ satisfies (\ref{E 65.004}) and $ST$ is an adapted process.\,\,\epf

\section{\label{PDK} Processes defined by Kernels}
A large class of adapted processes in $\CL^p\cmx(\FH)$ 
may be defined 
using $L^2$ kernels. Such kernel processes may be used to 
provide non-trivial examples of quantum semimartingales 
which satisfy the quantum Duhamel and quantum Ito formulae. 
They should not be confused with the  ``integral kernel 
operators'' studied by Obata[\S 4.3]{OB1}
For $i,j\in \BBN\cup  \{0\}$ let
$$
k_t{\,}^i_j(x;y)\equiv k^i_j(x;y;t),\quad  
(x;y;t)\in [0,1]^i\times[0,1]^j\times[0,1],
$$ 
where $x=(x_1,\ldots,x_i)$, $y=(y_1,\ldots,y_j)$, 
be a Lebesgue-measurable kernels completely symmetric in both 
$x_1,\ldots,x_i$ and $y_1,\ldots,y_j$. 
Suppose that $k_t{\,}^i_j$ is bounded in the $L^2$-norm, 
for each $t\in [0,1]$:
$$
\nm{k_t{\,}^i_j}=\left(\int_{[0,1]^i\times[0,1]^j} 
|k^i_j (x;y;t)|^2\,dx\,dy \right)^{1/2}< \infty,
$$
and let $\LGa_t$ be the scalar matrix with entries 
$\LGa_t{\,}^i_j=\|k_t{\, }^i_j\|$.
If $t\in [0,1]$ let $k^{\phantom{t]}\,i}_{t]\,j}$ be the 
restriction of $k^i_j$ to $[0,t]^i\times [0,t]^j\times[0,t]$.  
The formula
\beqn\label{E 74.02} 
(K^{\phantom{t]}\,i}_{t]\,j}\Gth)(x)
&=&\int_{[0,1]^j}
k^{\phantom{t]}\,i}_{t]\,j}(x;y;t)
\Gth(y)\,dy,\quad \Gth\in \FH^i_j,
\eeqn 
defines $K_{t]}{\,}^i_j$ in 
$\CB(\FH^{j}\lrb{t},\FH^{i}\lrb{t})\so$ 
with $\nm{K_{t]}{\,}^i_j}\leq \LGa_t{\,}^i_j$. 
\begin{Thm} 
Let $K_{t]}$ be the chaos matrix for $\FH_{t]}$ 
with entries $K_{t]}{\,}^i_j$ defined by {\rm (\ref{E 074.02})}.
\bb\ni
{\rm(i)} Then $K=\{K_t:t\in [0,1]\}$ is an adapted cmx process
where  $K_t$ is the ampliation of $K_{t]}$ defined by
(\ref{E 1.1}). 
bb\ni
{\rm (ii)} If $\LGa^i_j\in L^p[0,1]$ for all $i,j$ then 
$K\in \CL^p_{\rm cmx}(\FH)\so$.
\bb\ni
{\rm (iii)} If the matrix $[k_t{\,}^i_j]$ is $(2k+1)$-diagonal 
then so is $K_t$.
\end{Thm}
The adapted cmx process $K$ represents the adapted operator valued 
process $\tilde{K}=\{\tilde{K}_t:t\in [0,1]\}$. 
However, without extra conditions on $K$, the process $\tilde{K}$ 
need not be regular. Nor need the operators $\tilde{K}_t$ 
have a common dense domain. 
\bb\ni
Suppose $\LGa_t$ belongs to $\CB\cmx(\ell^2)$ for each $t\in [0,1]$.  
If $\psi\in \FH_{t]}$ then 
$\nu(\psi)=(\|\psi^0\|,\|\psi^1\|,\|\psi^2\|,\ldots)^t$ belongs to 
$\ell^2$ and the series $\sum_i\left(\sum_j\LGa^i_j(t)\nm{\psi^{j}}\right)^2$ 
is convergent. For each  $\psi\in \FH\lrb{t}$ define
\begin{equation} \label{E 074.01}
{\as{1.5}\ba{lcl}
(\tilde{K}\lrb{t}\psi)^{i}& =&\dyl\sum_{j=0}\iiy 
K^{\phantom{t]}\,i}_{t]\,j}\psi^{j}\\
\tilde{K}\lrb{t}\psi &=& \dyl\bigoplus_{i=0}\iiy(K\lrb{t}\psi)^{i}
\ea }
\end{equation}
These series converge since $\nu(\psi)\in \ell^2$ and
$$
\sum_{i=0}\iiy \nm{(\tilde{K}\lrb{t}\psi)^{i}}^2\leq
\sum_{i=0}\iiy\left(\sum_{j=0}\iiy 
\nm{K^{\phantom{t]}\,i}_{t]\,j}\psi^{j}}\right)^2\leq
\sum_{i=0}\iiy\left(\sum_{j=0}\iiy \LGa_t{\,}^i_j\nm{\psi^{j}}\right)^2
\leq \nm{\LGa_t}^2\nm{\psi}^2
$$
The operator $\tilde{K}_t$ is bounded with 
$\nm{\tilde{K}_t}\leq \nm{\LGa_t}$. We immediately have the 
following proposition.
\begin{Propn} The formulae (\ref{E 074.01}) define an adapted 
regular operator process  $\tilde{K}=\{\tilde{K}_t:t\in [0,1]\}$ with  
$\tilde{K}_t=\tilde{K}\lrb{t}\otsym I_{(t}$ 
such that $K\lrb{t}$ is a cmx representation of $\tilde{K}\lrb{t}$ 
and $K_t$ is a cmx representation of $\tilde{K}_t$. 
\bb\ni
If $\LGa$ belongs to the Bochner--Lebesgue space 
$L^p([0,1],\CB\cmx(\ell^2))$ then $\tilde{K}\in L^p([0,1],\CB(\FH))$.
\end{Propn}
Given a $(2k+1)$-diagonal matrix $\hat{k}^i_j$ of 
$L^2$-kernels with $\nm{k_t{\,}^i_j}<\leq C$ it is possible to 
construct a variety of $(2k+1)$-diagonal processes in $L^p([0,1],\CB\cmx(\ell^2))$. Choose functions 
$t\mapsto c^i_j(t)=c^j_i(t)$ with $|c^i_j|<\Ga\in L^p[0,1]$ and 
put ${k}_t{\,}^i_j\defeq c^i_j(t)\hat{k}_t{\,}^i_j$. Then the 
process $K$ defined by (\ref{E 74.02}) belongs to 
$L^p([0,1],\CB\cmx(\ell^2))$.   
\bb\ni 
If the matrices $\{\LGa_t:t\in[0,1]\}$ have common domain 
$\FD\subset \ell^2$ then $\tilde{K}$ is an operator process 
the operators $\{\tilde{K}_t:t\in [0,1]\}$ have common  domain $\CD$ 
such that $\{\psi\in \FH:\nu(\psi)\in
\FD\}\subset \CD$.
If $\FD$ contains a dense subset of the positive cone in $\ell^2$ then 
$\CD$ is dense in $\FH$.

\section{\label{QSI} Quantum Stochastic Integrals}
Lindsay \cite[Definition 2.1, Proposition 2.4]{Lin1} extends 
Hudson and Parthasarathy's definition of quantum stochastic 
integrals \cite{HP1} to include non-adapted integrands.
Lindsay's integrals are defined in terms of the Hitsuda-Skorohod 
process $t\mapsto \CS_t$ and the 
gradient process $t\mapsto\nabla_t$. The cmx 
representations of these two processes are given by (\ref{M 73.01}). 
Their simple form facilitates a  formal calculation
of the cmx representations of Lindsay's quantum stochastic 
integrals. This leads us to
definition of the cmx quantum stochastic integrals of cmx processes 
compatible with Lindsay's definition for operator processes.
\bb\ni
The integrators are the quadruple $(\GL,A,A\d,t)$ of 
basic processes: 
\bb\ni
the {\it gauge} process ${\GL}:t\row 
{\GL}(\chi_{[0,t]})$; 
\bb\ni
the {\it annihilation } process, ${A}:t\row  {a}(\chi_{[0,t]})$;
\bb\ni
the {\it creation } process, ${A}\d:t\row  a\d(\chi_{[0,t]})$; 
\bb\ni
the {\it time} process $tI=t$. 
\bb\ni
In the definition of the gauge process  
$\chi_{[0,t]}\in \CB(L^2[0,1])$ is the operator of 
multiplication by $\chi_{[0,t]}$;
\bb\ni
A quadruple $(E,F,G,H)$ of cmx processes is {\it integrable} 
if $E\in \CL\iiy_{\rm cmx}(\FH)\so$, 
$F,G\in \CL^2_{\rm cmx}(\FH)\so$, and $H\in \CL^1_{\rm 
cmx}(\FH)\so$ and {\it symmetric} if $E=E\st$, $G=F\st$ and 
$H=h\st$. The processes $E$, $F$, $G$, and $H$ are  
integrands for $d\GL$, $dA$, $dA\d$ and $dt$ respectively. This is consistent 
with, but wider than,  the definition, (\ref{E 075.04}) of an integrable quadruple of 
operator processes introduced in \cite{Att}.
\bb\ni
For the sake of clarity and to prevent a `debauch of indices' in the proofs 
we have kept to this old-fashioned notation rather than using 
Evans' notation \cite{Ev1}, \cite[V \S 2.1]{MEY}. However Mr A 
Belton's extension of the Evans notation allows 
a very efficient proof of the local Ito Product formula, Lemma 
\ref{L 76.211}. In the higher dimensional case the benefits of 
the Evans-Belton notations are even greater.
\bb\ni
Let $(E,F,G,H)$ be an integrable quadruple of cmx processes. The formulae
\begin{eqnarray}
\label{E 75.101} M(E:d\GL)= \CS E\nabla, 
&& M_t(E:d\GL)^{i+1}_{j+1}= 
\CS_t{}^{i+1}_iE^i_j\nabla_{t}{}^j_{j+1},
\\
\label{E 75.102} M(F:dA)= \int_0^{1}\! F_{s}\nabla(s)\,ds, && 
M_t(F:dA)^{i}_{j+1}= \int_0^{1}\! F_{s}{}^i_j
\nabla_{t}{}^j_{j+1}(s)\,ds,
\\
\label{E 75.103} M(dA\d:G)= \CS G, && 
M_t(dA\d:G)^{i+1}_{j}= \CS_{t}{}^{i+1}_i G^i_j,
\\
\label{E 75.104} M_t(H:\,dt)=\int_0^t\! H_{s}\,ds, && 
M_t(H:\,dt)^i_j=\int_0^t \!H_{s}{}^{i}_{j}\,ds,
\end{eqnarray}
$$
M_t(E:d\GL)^0_j=M_t(E:d\GL)^i_0 = M_t(F:dA)^i_0 =M_t(dA\d:G)^0_j =0
$$
where $i,j=0,1,2,\ldots$, define cmx processes
$$
M(E:d\GL),\quad M(F:dA), \qquad M(dA\d:G), \qquad M(H,dt). 
$$
It may be helpful to think of these processes as
cmx representations of the 
Hudson--Parthasarathy quantum stochastic integrals 
$$
\int_0^t \tilde{E_s}\,d\GL_s,\quad\int_0^t \tilde{F}_s\,dA_s,\quad
\int_0^t dA\d_s \tilde{G}_s,\quad \int_0^t \tilde{H}_s\,ds.          
$$
We examine the first three of these definitions in more detail. 
In (\ref{E 75.101}) $E^i_j$  denotes the mapping 
$E^i_j:L^2([0,1],\FH^j) \row L^2([0,1],\FH^i)$ defined by the 
formula
$$
(E^i_j\Gv)(s)=E_s{}^i_j\Gv(s),\qquad \Gv\in L^2([0,1],\FH^j).
$$
Consider the diagram
$$
\FH^{j+1}\stackrel{\nabla_t{}^j_{j+1}}{\longrightarrow} L^2([0,1],\FH^j)
\stackrel{E^i_j}{\longrightarrow} L^2([0,1],\FH^i)
\stackrel{\CS_t{}^{i+1}_i}{\longrightarrow} \FH^{i+1}
$$
The transformations are all bounded so that 
$M_t(E:d\GL)^{i+1}_{j+1}= \CS_t{}^{i+1}_i E^i_j 
\nabla_t{}^j_{j+1}$ is a bounded linear 
transformation and (\ref{E 75.101}) defines a 
cmx process $M(E:d\GL)$.
\bb\ni
Since $\|\nabla_t{}^j_{j+1}\|\leq (j+1)^{1/2}$ it follows that $s\mapsto 
F_s{}^i_j\nabla_t{}^j_{j+1}$ is in the Lebesgue space
$L^2([0,1],\CB(\FH^{j+1},\FH^i)\so)$ and (\ref{E 75.102}) defines a 
cmx process $M(F:dA)$.
\bb\ni
If $\psi\in \FH^j$ then $s\mapsto G_s{}^i_j\psi$ is a process $G^i_j\psi$ in 
$L^2([0,1],\FH^i)$ and $\psi\mapsto \CS_t{}^{i+1}_i G^i_j\psi$ is a 
linear transformation $\CS_t{}^{i+1}_i G$ in 
$\CB(\FH^j,\FH^{j+1})$. Therefore (\ref{E 75.103}) defines a 
cmx process $M(dA\d:G)$.
\bb\ni
To avoid reindexing we define $E^i_j=F^i_j=G^i_j=0$ whenever 
$i<0$ or $j<0$. 
\bb\ni
The following theorem shows that  the integrals defined above 
belong to $\CL\iiy_{\rm cmx}(\FH)\so$ and that it is natural to 
choose the respective integrands in 
the spaces $\CL\iiy_{\rm cmx}(\FH)\so$, 
$\CL^2_{\rm cmx}(\FH)\so$, $\CL^2_{\rm cmx}(\FH)\so$ 
and $\CL^1_{\rm cmx}(\FH)\so$. 
\begin{Thm}\label{T 74.801}
The formulae {\rm(\ref{E 75.101}), (\ref{E 75.102}), (\ref{E 75.103})} and 
{\rm (\ref{E 75.104})} define processes in $\CL\iiy_{\rm cmx}(\FH)\so$. 
Moreover
\begin{eqnarray}
\label{E 75.201}(M_t(E:d\GL)^{i+1}_{j+1})\st&=&M_t(E\st:d\GL)^{j+1}_{i+1}
\\
\label{E 75.202}(M_t(F:dA)^i_{j+1})\st&=& M_t(dA\d:F\st)^{j+1}_i, 
\\
\label{E 75.203}(M_t(H:dt)^i_{j})\st&=& M_t(H\st:dt)^{j}_i,
\\ 
\label{E 75.204}\nm{M(E:d\GL)^{i+1}_{j+1}}\iy
&\leq& (i+1)^{\frac12}(j+1)^{\frac12}\nm{E^i_j}\iy,
\\
\label{E 75.205}\nm{M(F:dA)^{i}_{j+1}}\iy&\leq&(j+1)^{\frac12}\nm{F^i_j}_2,
\\
\label{E 75.206}\nm{M(H:dt)^{i}_{j}}\iy&\leq&\nm{H^i_j}_1.
\end{eqnarray}
\end{Thm}
\pf The measurability of the processes follows from the 
measurability of 
$t\mapsto \nabla(t)$ and $t\mapsto \CS_t$ and an extension of Corollary 
\ref{CC 73.01}. 
It follows from (\ref{E 73.211}) that
$$      
\nm{M(E:d\GL)^{i+1}_{j+1}}\iy\leq\sup_{t\in [0,1]}
\nm{\CS_{t}{}_i^{i+1}}\nm{E_{t}{}^i_j}\nm{\nabla_{t}{}^j_{j+1}}=
(i+1)^{\frac12}(j+1)^{\frac12}\nm{E^i_j}\iy,
$$
proving (\ref{E 75.204}). If $\psi\in \FH^{j+1}$ then 
\begin{eqnarray*}
\nm{M_t(F:dA)^{i}_{j+1}\psi}&\leq& 
\int_0^t\nm{F_{s}{}^i_j}\nm{\nabla_{t}{}^j_{j+1}(s)\psi}\,ds
\\
&\leq& \left(\int_0^{1}\nm{F_{s}{}^i_j}^2\,ds\right)^{\frac12}  
\left(\int_0^{1}\nm{\nabla_t{\,}^j_{j+1}(s)\psi}^2\,ds\right)^{\frac12}  
\\
&\leq& (j+1)^{\frac12}\nm{F^i_j}_2 \nm{\psi},
\end{eqnarray*}
proving (\ref{E 75.205}). A similar argument gives (\ref{E 75.206}).
\bb\ni
Suppose $f,g\in L^2(0,1)$. It follows from (\ref{E 73.101}) that
\begin{eqnarray}
\nonumber 
\ip{M_t(E:d\GL)^{i+1}_{j+1}e^{j+1}(f)}{e^{i+1}(g)} &=& 
\ip{\CS_t{\,}^{i+1}_i E^i_j\nabla_t{\,}^j_{j+1}e^{j+1}(f)}{e^{i+1}(g)}
\\
\nonumber 
&&\hspace{-1.2in} = 
\int_0^{1}\ip{ E_{s}{\,}^i_j\nabla_t{\,}^j_{j+1}(s)e^{j+1}(f)}
{\nabla_t{\,}^i_{i+1}(s)e^{i+1}(g)}\,ds
\\
&&\hspace{-1.2in} = 
\label{E 74.067}
\int_0^{t}f(s)\ol{g(s)}\ip{ E_{s}{\,}^i_j e^{j}(f)}
{e^{i}(g)}\,ds
\\
\nonumber 
&&\hspace{-1.2in} = 
\int_0^{t}f(s)\ol{g(s)}\ip{ e^{j}(f)}{(E\st)_{s}{\,}^j_i e^{i}(g)}\,ds
\\
\nonumber 
&&\hspace{-1.2in} = 
\int_0^{1}\ip{ \nabla_t{\,}^j_{j+1}(s)e^{j+1}(f)}
{(E\st)_{s}{\,}^j_i\nabla_t{\,}^i_{i+1}(s)e^{i+1}(g)}\,ds
\\
\nonumber 
&&\hspace{-1.2in} = 
\ip{e^{j+1}(f)}
{\CS_t{\,}^{j+1}_j(E\st)^j_i\nabla_t{\,}^i_{i+1}e^{i+1}(g)}\,ds
\\
\nonumber 
&&\hspace{-1.2in} =
\ip{e^{j+1}(f)}{M_t(E\st:d\GL)^{j+1}_{i+1}e^{i+1}(g)}.
\end{eqnarray}
Since $\{e^{n}(f):f\in L^2[0,1]\}$ is total in $\FH^n$ and 
$M_t(E:d\GL)^{i+1}_{j+1}$ is bounded this proves (\ref{E 75.201}).  
Put $G=F\st$. It follows from (\ref{E 73.101}) that 
\begin{eqnarray}
\nonumber 
\ip{M_t(dA\d:G)^{i+1}_{j+1}e^{j+1}(f)}{e^{i+1}(g)}
&=& \ip{\CS_t{\,}^{i+1}_i G^i_{j+1}e^{j+1}(f)}{e^{i+1}(g)}
\\
\nonumber 
&&\hspace{-0.8in} = \int_0^1\ip{G_s{\,}^i_{j+1}e^{j+1}(f)}{\nabla_t{\,}^i_{i+1}(s)e^{i+1}(g)}\,ds
\\
\label{E 74.068}
&&\hspace{-0.8in} = \int_0^t\ol{g(s)}\ip{G_s{\,}^i_{j+1}e^{j+1}(f)}{e^{i}(g)}\,ds
\\
\label{E 74.0068}
&&\hspace{-0.8in} = \int_0^t\ol{g(s)}\ip{e^{j+1}(f)}{F_s{\,}^{j+1}_i e^{i}(g)}\,ds
\\
\nonumber 
&&\hspace{-0.8in} = \int_0^1\ip{e^{j+1}(f)}{F_s{\,}^{j+1}_i\nabla_t{\,}^i_{i+1}(s)e^{i+1}(g)}\,ds
\\
\nonumber 
&&\hspace{-0.8in} = \ip{e^{j+1}(f)}{M_t(F:dA)^{j+1}_{i+1}e^{i+1}(g)},
\end{eqnarray}
and (\ref{E 75.202}) follows as above. A similar argument gives 
(\ref{E 75.203}).  \,\,\epf 
\bb\ni
The {\it cmx quantum stochastic integral}, $M=M(E,F,G,H)$, 
of an integrable quadruple $(E,F,G,H)$ of chaos matrices is the process
\beqn\label{E 73.11}
M_t=  M_{t}(E:\,d\GL)+ M_{t}(F:\,dA)+M_{t}(\,dA\d:G)+M_t(H:\,dt).
\eeqn
The term cmx quantum stochastic integral may also be used for 
the definite integral.
\bb\ni
$M$ is {\it regular} if $M_t\in \CB\cmx(\FH)$ for $t\in [0,1]$.
We shall us the notation
\begin{eqnarray*}
M_t &=& \int_0^t (E_s\,d\GL_s+F_s\,dA_s+\,dA\d_s G_s+H_s\,ds)
\end{eqnarray*}
although the subscript $_s$ in the integrand is sometimes suppressed. 
If $0\leq t_1, t_2\leq 1$ define, 
$$
\int_{t_1}^{t_2} (E\,d\GL+F\,dA+\,dA\d G+H\,ds)
= M_{t_2}(E,F,G,H)-M_{t_1}(E,F,G,H).
$$
The integral behaves in behaves in the usual way with respect to $t_1$ and 
$t_2$. 
\bb\ni
The {\it control matrix}  of $M$ (of the integrable quadruple 
$(E,F,G,H)$) is the scalar matrix $\LGk=\LGk(M)$ with entries
\beqn\label{E 64.011}
\LGk^i_j=i^{1/2}\nm{E^{i-1}_{j-1}}\iy j^{1/2}+
\nm{F^i_{j-1}}_2 j^{1/2}+i^{1/2}\nm{G^{i-1}_j}_2+\nm{H^i_j}_1,
\eeqn
with the convention that $E^i_j=F^{\nu}_j=G^i_{\nu}=0$ if $i<0$ or 
$j<0$. It follows from (\ref{E 75.204}), (\ref{E 75.205}) and 
(\ref{E 75.206}) that 
$$
\nm{M^i_j}\iy \leq \LGk^i_j \qquad i,j=0,1,2,\ldots.
$$
\begin{Cor} If $\LGk$ is the cmx representation of a bounded operator 
in $\ell^2$ then $M=\{M_t:t\in [0,1]\}$ is a regular cmx 
process.
\end{Cor}
The following corollary to Theorem \ref{T 74.801} is 
a chaos matrix version of \cite[Theorems 4.1, 4.4]{HP1} and is 
used to prove the uniqueness of cmx quantum stochastic integrals. 
It may be shown that, for adapted integrands, 
the definition of cmx quantum stochastic integral is compatible 
with Hudson and Parthasarathy's definition of quantum stochastic 
integral  \cite{HP1}. 
\begin{Cor}
If $M=M(E,F,G,H)$ is the cmx quantum stochastic integral of an integrable 
quadruple of cmx processes and $f,g\in L^2[0,1]$ then
\begin{eqnarray}\label{EE 77.011} 
\ip{M_t{\,}^{i+1}_{j+1} e^{j+1}(f)}{e^{i+1}(g)} &=& \\
\nonumber&& \hspace{-1.25in}\int_0^t
\left(f(s)\ol{g(s)}\ip{E_s{\,}^i_j e^j( f)}{e^i(g)}
+ f(s)  \ip{F_s{\,}^{i+1}_j e^j( f)}{e^{i+1}(g)}
\right.
\\
\nonumber&& \hspace{-1in}
+\left. \ol{g(s)}\ip{G_s{\,}^i_{j+1} e^{j+1}( f)}{e^i(g)}+ 
\ip{H_s{\,}^{i+1}_{j+1} e^{j+1}(f)}{e^{i+1}(g)}\right)\,ds    
\end{eqnarray}
\end{Cor}
\pf. The terms on the right involving $E$, $F$ and $G$ follow from 
(\ref{E 74.067}), (\ref{E 74.068}) and (\ref{E 74.0068}). 
The remaining term is immediate from the definition of $M(H:dt)$.\,\,\epf
\bb\ni
\begin{Cor}\label{CC 75.01}  The mapping $t\mapsto M_t{\,}^{i+1}_{j+1}$ is 
strongly continuous and strongly measurable.
\end{Cor}
\pf Since  $\nm{M_t{\,}^{i+1}_{j+1}}\leq \LGk^{i+1}_{j+1}$ 
it follows from Lebesgue's dominated convergence 
theorem that 
$t\mapsto \ip{M_t{\,}^{i+1}_{j+1} e^{j+1}(f)}{e^{i+1}(g)}$ is 
continuous whenever  $f,g\in L^2[0,1]$.  The conclusion follows 
since $\{e^k(f):f\in L^2[0,1]\}$ is total in $\FH^k$ for $k=0,1,2,\ldots$.
\,\,\epf
\bb\ni
The following Proposition extends the uniqueness theorems of 
Parthasarathy \cite{P2} Vincent--Smith \cite{Vin} and Lindsay \cite{Lin2}
to cmx quantum stochastic integrals.
\begin{Propn}\label{P 74.01} If $(E,F,G,H)$ is an integrable 
quadruple and  
$$
\int_0^t(E\,d\GL+F\,dA+G\,dA\d+H\,ds)=0
$$
then $E\equiv F\equiv G\equiv H\equiv 0$, the zero process.
\end{Propn}
\pf We follow Lindsay's line of argument.
\bb\ni
Let $L\step_{\BBQ}=L\step_{\BBQ}[0,1]$ be the vector space, over 
$\BBQ$, of 
$\BBQ$-valued step functions with discontinuity points in $\BBQ$ and, for 
$s\in [0,1]$, let $L\step_{\BBQ}(\hat{s})
=\{f\in L\step_{\BBQ}:f=0\mbox{ in a neighbourhood of } s\}$. Then both 
$\{e^j(f):f\in L\step_{\BBQ}(\hat{s})\}$ and 
$\{e^j(f):f\in L\step_{\BBQ}\setminus L\step_{\BBQ}(\hat{s})\}$ 
are total subsets of $\FH^j$. 
\bb\ni
The integrand in (\ref{EE 77.011}) is zero almost everywhere for each $f, g$ in 
the countable set $L\step_{\BBQ}$. Therefore there exists a 
null set $\CN$ such that
\begin{eqnarray*}
f(s)\ol{g(s)}\ip{E_s{\,}^i_j e^j( f)}{e^i ( g)}
+ f(s)  \ip{F_s{\,}^{i+1}_j e^j( f)}{e^{i+1} ( g)} \hspace{0.5in} &&\\
\hspace{0.5in}+\ol{g(s)}\ip{G_s{\,}^i_{j+1} e^{j+1}( f)}{e^i ( g)} 
+\ip{H_s{\,}^{i+1}_{j+1} e^{j+1}( f)}{e^{i+1} ( g)}&=&0    
\end{eqnarray*}
whenever $s\in [0,1]\setminus \CN$, for all $ f,g\in L\step_{\BBQ}$ and 
$i,j=0,1,2,\ldots$.
\bb\ni
Fix $s\in [0,1]\setminus \CN$. Then 
$\ip{H_s{\,}^{i+1}_{j+1} e^{j+1}( f)}{e^{i+1} ( g)}=0 $
whenever $f,g\in L\step_{\BBQ}(\hat{s})$. Since $H_s{\,}^{i+1}_{j+1}$ is bounded it 
follows that $H_s{\,}^{i+1}_{j+1}=0$ for all $i,j$. Therefore $H\equiv 0$.
\bb\ni
Similarly $f(s)\ip{F_s{\,}^{i+1}_j e^j( f)}{e^{i+1} ( g)}=0$ whenever 
$f\in L\step_{\BBQ}$ and $g\in L\step_{\BBQ}(\hat{s})$. Since 
$\{e^i(g):g\in L\step_{\BBQ}(\hat{s})$ is total in $\FH^j$ it follows that
$F_s{\,}^i_j e^j( f)=0$ whenever 
$f\in L\step_{\BBQ}\setminus L\step_{\BBQ}(\hat{s})$.
These $e^j(f)$ are total in $\FH^j$ and $F_s=0$. Therefore $F\equiv 0$ and, 
taking adjoints, $G\equiv 0$.
\bb\ni
Finally  $\ip{E_s{\,}^i_j e^j( f)}{e^i ( g)}=0$ whenever 
$f,g\in L\step_{\BBQ}\setminus L\step_{\BBQ}(\hat{s})$ and, 
as in the other cases, $E\equiv 0$.\,\,\epf 
\bb\ni
We shall also need the following convergence theorem.
\begin{Thm}\label{T 74.802}
Let $(E_n,F_n,G_n,H_n)$ be a sequence of integrable quadruples of 
cmx processes and suppose $E_n$, 
$F_n$, $G_n$ and $H_n$ converge to $E$, $F$, $G$, $H$ respectively. Then 
$M(E_n,F_n,G_n,H_n)$ converges to $M(E,F,G,H)$
\end{Thm}
\pf This follows directly from Theorem \ref{T 74.801}. For example, for all 
$i,j$, 
$$
\lim_{n\roy}\nm{M(E_n-E:d\GL)^{i}_{j}}\iy 
\leq \lim_{n\roy}i^{1/2} j^{1/2}\nm{E^{i-1}_{j-1}-E_n{\,}^{i-
1}_{j-1}}\iy=0,
$$
by the inequality (\ref{E 75.204}). The inequalities (\ref{E 75.205}) and 
(\ref{E 75.206}) may be similarly used to prove the convergence 
in $\CL\iiy_{\rm cmx}(\FH)\so$ of the 
sequences $M(F_n-F:dA)$, $M(dA\d:G_n-G)$ and $M(H_n-H:dt)$ 
to the zero process.\,\,\epf 
\bb\ni
\begin{Thm}\label{T 74.803} If $(E,F,G,H)$ is an adapted 
integrable cmx quadruple
then the cmx quantum stochastic integral $M(E,F,G,H)$ 
is an adapted cmx process.
\end{Thm}
\pf The cmx quantum stochastic integrals are given by (\ref{E 75.101}), 
(\ref{E 75.102}), (\ref{E 75.103}), (\ref{E 75.104}) and (\ref{E 73.11}). 
Consider first $M=M(F:dA)$ given by (\ref{E 75.102}). If $t<s$ then 
$$ M_t\nabla(s)=\int_0^t F_u\nabla(u)\nabla(s)\,du=
\int_0^t F_u\nabla(s)\nabla(u)\,du= \int_0^t \nabla(s)F_u\nabla(u)\,du
$$
$$
= \nabla(s)\int_0^t F_u\nabla(u)\,du=\nabla(s) M_t.
$$
These identities may be verified by using the entry-wise definition of $M$ in 
(\ref{E 75.102}) and the norm bounds (\ref{E 75.205}) and computing 
entry-wise. In particular the norm bounds on the entries allow the interchange 
of $\nabla(s)$ and the integral over $[0,t]$.
\bb\ni
If $\psi\in L^2([0,1],\FH^j)$ and $t<s$ then
\begin{eqnarray*}
\nabla(s)(\CS_t\psi)(x_1,\ldots,x_{j+1}) &=&
(j+1)^{\frac12}(\CS_t\psi)(x_1,\ldots,x_j,s)
\\
&=& 
\sum_{k=0}^j\chi_{[0,t]}(x_k)\psi(x_k)(x_1,\ldots,\hat{x}_k,\ldots,x_{j},s)
\\
&& \hspace{0.5in}
+ \chi_{[0,t]}(s)\psi(s)(x_1,\ldots,x_k,\ldots,x_{j})
\\
&=& j^{-\frac12}\sum_{k=0}^j\chi_{[0,t]}(x_k)(\nabla(s) 
\psi)(x_k)(x_1,\ldots,\hat{x}_k,\ldots,x_j)
\\
&=& (\CS_t\nabla(s)\psi)(x_1,\ldots,x_j).
\end{eqnarray*}
Therefore $\nabla(s)\CS_t=\CS_t\nabla (s)$ whenever $t<s$.
\bb\ni
Now $M\st =M(dA\d:G)$ where $G=F\st$. Since $F$ is adapted so is $G$.
$$
M\st_t\nabla(s)=\CS_t G_t\nabla(s)=\CS_t \nabla(s)G_t=
\nabla(s)\CS_t G_t=\nabla(s) M\st_t
$$ 
whenever $t<s$ and $M$ is adapted. Since $M(dA\d:G)\st=M(G\st:dA)$ it follows 
that $M(dA\d:G)$ is adapted.
\bb\ni
If $M= M(E:d\GL)$ then, since $E$ is adapted 
$$
M_t\nabla(s)=\CS_t E_t \nabla_t \nabla(s)=\CS_t E_t \nabla(s)\nabla_t =
\nabla(s)\CS_t E_t \nabla_t=\nabla(s)M_t
$$
whenever $t<s$. Since $E=E\st$ it follows that $M=M\st$ so that $M$ is 
adapted. We omit the proof that $M(H:dt)$ is adapted. \,\,\epf 
\bb\ni
Although the natural integrands for quantum stochastic 
integrals are integrable quadruples of 
operator processes the processes are  more amenable when the 
operators are bounded. For example the quantum 
semimartingales of Attal and Meyer \cite{AM} have integrands in 
the Lebesgue spaces $L^p([0,1],\CB(\FH)\so)$ and the 
quantum Ito formula \cite{Vin2} is valid for regular quantum 
semimartingales. 
\bb\ni
A {\it quantum semimartingale} (non-commutative 
semimartingale \cite[\S II]{Att}) is 
an adapted quantum stochastic process $\hat{M}=\{\hat{M}_t:t\in [0,1]\}$ with 
representation as a Hudson--Parthasarthy quantum stochastic 
integral
\beqn\label{E 75.0004}
\hat{M}_t
=M_0+\int_0^t (\hat{ E}_s\,d{\GL}_s+\hat{ F}\,d{ A}+
d{ A}\d{\,}\hat{G}+\hat{H}_s\,d{s}),\qquad t\in [0,1],
\eeqn
where the adapted processes 
$(\hat{E},\hat{F},\hat{G},\hat{H})$  with  
$$
\hat{E}\in L\iiy([0,1],\CB(\FH)\so),\quad \hat{F}, 
\hat{G}\in L^2([0,1],\CB(\FH)\so), 
$$
\beqn
\label{E 075.04}\quad \hat{H}\in L^1([0,1],\CB(\FH)\so),
\eeqn
form an {\it integrable quadruple}.
\bb\ni
If there is no initial space $\hat{M}_0$ is a multiple of the 
identity. Otherwise $\hat{M}_0$ is the 
ampliation of an operator in the initial space. It is shown in 
\cite[\S 6]{MEY}, \cite[\S 4]{HP1} that $\hat{M}_t$ 
is defined as an operator in $\CH$ with 
domain $\CE$ for each $t\in [0,1]$.  Sometimes $\hat{M}_t$ is extended to its 
maximal domain \cite{AM}. The representation (\ref{E 75.0004}) 
of $\hat{M}$ is unique.
\bb\ni
If $\hat{E}_t=(\hat{E}_t)\st$, $\hat{G}_t=(\hat{F}_t)\st$ and 
$\hat{H}_t=(\hat{H}_t)\st$ for almost all $t\in [0,1]$ then 
$\hat{M}$ is said to be {\em symmetric}. In this case 
$\hat{M}_t$, with domain $\CE$, is a symmetric operator for each 
$t\in [0,1]$.
\bb\ni
If $\hat{M}$ is symmetric and the operators $\hat{M}_t$ are 
essentially self-adjoint on a common core for almost all $t\in [0,1]$ 
then $\hat{M}$ is {\it essentially self-adjoint}. 
\bb\ni
If $\hat{M}_t$ is a bounded operator  
then $\hat{M}_t\in \CB(\FH)$ will always denote the unique 
extension of $\hat{M}_t$ to 
$\FH$. If $\hat{M}\in L\iiy([0,1],\CB(\FH)\so)$
then the quantum semimartingale $\hat{M}$ is said 
to be {\it regular}. If $\hat{M}$ is regular and symmetric then 
$\hat{M}$ is said to be {\it self-adjoint}.
\bb\ni
The vector space of quantum semimartingales is denoted $\CS\pr$ 
and that of regular quantum semimartingales is denoted 
$\CS$.
\bb\ni
A {\it cmx semimartingale} is a cmx quantum stochastic integral
\beqn\label{E 75.0005}
M_t=M_0+\int_0^t (E_s\,d\GL_s+ F_s\,d A_s+ G_s\,d A\d_s+ 
H_s\,d s), \qquad t\in [0,1],
\eeqn
where $(E,F,G,H)$ is an adapted integrable quadruple with .
$$
E\in L\iiy\cmx([0,1],\CB\cmx(\FH)\so),\quad
F,G\in L^2\cmx([0,1],\CB\cmx(\FH)\so),
$$
$$
H\in L^1\cmx([0,1],\CB\cmx(\FH)\so)
$$
It follows from Proposition \ref{P 74.01} and Theorem 
\ref{T 74.803} 
that $M$ is adapted and well defined by (\ref{E 75.0005}).
\bb\ni
If $(E,F,G,H)$ is symmetric ${M}$ is said to be {\it symmetric}. 
\bb\ni
If $M$ is regular it follows 
from Corollary \ref{CC 75.01} (ii) that 
$M\in L\iiy\cmx([0,1],\CB\cmx(\FH)\so)$. If $M$ is regular 
and symmetric then $M$ is said to be {\it self-adjoint}.
\bb\ni
The vector spaces of quantum semimartingales and  cmx semimartingales 
are denoted $\CS\pr$ and $\CS_{\rm cmx}\pr$. The vector spaces of 
regular quantum semimartingales and regular 
cmx semimartingales are denoted $\CS$ and $\CS_{\rm cmx}$.
\begin{Thm}\label{T 74.444} 
Let $\hat{M}$ be a quantum semimartingale
\beqn\label{E 74.73} 
\hat{M}_t=\int_0^t (\hat{E}\,d\GL+\hat{F}\,d 
A+\hat{G}\,d A\d+\hat{H}\,d s)
\eeqn
and let $M$ be the cmx semimartingale, 
$$
M_t=\int_0^t (E\,d\GL+F\,d A+G\,d A\d+H\,d s),
$$
where $E,F,G,H$ are the cmx representations of 
$\hat{E},\hat{F},\hat{G},\hat{H}$ respectively.
Then
\bb\ni
{\rm(i)} $M$ is a cmx representation of $\hat{M}$ with 
$\CE\subset\FD(M_t)$ for all $t\in [0,1]$. 
\bb\ni
{\rm(ii)} The mapping $M\mapsto \hat{M}$ is a bijection 
from $\CS\pr_{\rm cmx}$ onto $\CS\pr$.
\bb\ni
{\rm(iii)} The mapping $M\mapsto \hat{M}$ is an isometric 
isomorphism from $\CS_{\rm cmx}$ onto $\CS$.
\end{Thm}
\pf It follows from Proposition \ref{PP 73.101} and Corollary \ref{C 074.01} 
that the formula (\ref{E 74.73}) defines a quantum 
semimartingale $\hat{M}$. If $t\in [0,1]$ then
$\hat{M}_t$ is characterised in \cite{HP1} and \cite{MEY} as 
the unique operator in $\CL(\FH)$ with domain $\CE$ such that, 
whenever $f,g\in L^2[0,1]$,
\begin{eqnarray*}
\ip{\hat{M}_t e(f)}{e(g)}
&=&\int_0^t\ip{(f\bar{g}\hat{ E}+f\hat{ F}
+\bar{g}\hat{G}+\hat{H})(s)e(f)}{e(g)}\,ds.
\end{eqnarray*}
We shall only consider the case $\hat{F}=\hat{G}=\hat{H}=0$. The proofs of the  
three other cases with a single non-zero integrand are simple variants.
\bb\ni
If $g,h\in L^2[0,1]$ then
\begin{eqnarray}
\nonumber\sum_{i,j=0}\iiy\ip{g(s) E_s{\,}^i_j e^j(g)}{h(s) e^i(h)}
&\leq & |g(s)h(s)| \sum_{i,j=0}\iiy\nm{E^i_j}\iy|\ip{e^j(g)}{e^i(h)}|
\\
{\label{S 64.01}} &\leq& |g(s)h(s)| \nm{ \tilde{E}}\iy \ip{e(|g|)}{e(|h|)},
\end{eqnarray}
which is an integrable function. 
\begin{eqnarray}
\nonumber \ip{{ \hat{M}}_te(g)}{e(h)}
&=&
\int_0^t \ip{g(s) \hat{E}_s e(g)}{ h(s) e( h)}\,ds
\\
\label{E 74.70} &=& \int_0^t \ip{g(s) \sum_{j=0}\iiy E_s{\,}^i_j e^j(g)}{ 
h(s)\sum_{i=0}\iiy e^i( h)}\,ds
\\
\label{E 74.71} &=& \sum_{i,j=0}\iiy \int_0^t\ip{g(s)  E_s{\,}^i_j e^j(g)}{ 
h(s)e^i( h)}\,ds
\\
\nonumber&=& \sum_{i,j=0}\iiy \int_0^t\ip{E_s{\,}^i_j \nabla_s e^{j+1}(g)}{ 
h(s)\nabla_s e^{i+1}(h)}\,ds
\\
\label{E 74.72} &=& \sum_{i,j=0}\iiy \ip{M_t{\,}^i_j e^j(g)}{e^i(h)}.
\end{eqnarray}
(\ref{E 74.70}) follows from the absolute convergence of the 
left hand series in (\ref{S 64.01}). The third equality is a 
consequence of Lebesgue's dominated convergence theorem gives 
(\ref{E 74.71}) and (\ref{E 74.72}) follows from the definition 
(\ref{E 75.101}) and (\ref{E 73.101}). Now
$$
\nm{M_t{\,}^i_j e^j(g)}\leq 
i^{\frac12}j^{\frac12}\nm{E}\iy \frac{\nm{g}^j}{(j!)^{\frac12}}\leq 
i^{\frac12}\nm{E}\iy \frac{\nm{g}^j}{((j-1)!)^{\frac12}}
$$
and $\sum_{j=0}\iiy M_t{\,}^i_j e^j(g)$ converges to $\eta^i$ in $\FH^i$. 
Therefore, if $z\in \BBC$, 
\begin{eqnarray*}
\ip{\hat{M}_te(g)}{e(z  h)}&=&
\sum_{i,j=0}\iiy z ^i \ip{M_t{\,}^i_j e^j(g)}{e^i(h)}
\\
&=&
\sum_{i=0}\iiy z ^i\ip{\eta^i}{e^i(h)};
\\
\ip{\hat{M}_t e(g)}{e(z  h)}&=&\sum_{i=0}\iiy
z ^i\ip{(\hat{M}_t e(g))^i}{e^i(h)},
\end{eqnarray*}
and
$$
\sum_{i=0}\iiy z ^i\ip{(\hat{M}_t e(g))^i
-\eta^i}{e^i(h)} =0.
$$
Therefore $\ip{(\hat{M}_t) e(g))^i-\eta^i}{e^i(h)} =0$ 
for all $h\in 
L^2[0,1]$. Since the exponential domain is dense in $\FH$ it follows that
$(\hat{M}_t e(g))^i=\eta^i$ for all $i$ and 
$\sum_{i=0}\iiy \eta^i$ is convergent 
in $\FH$ to  $\hat{M}_t e(g)$. Therefore $e(g)\in 
\FD(M_t)$  and 
$$
\tilde{M}_t e(g)=M_t e(g)=\hat{M}_t e(g),\qquad t\in 
[0,1],\,\,g\in L^2[0,1],
$$
Therefore $M$ is a cmx representation of $\hat{M}$. 
This proves (i) in the case $F=G=H=0$. 
\bb\ni
By linearity $M_t \psi=\tilde{M}_t \psi$ for all 
$\psi\in \CE$ and $M_t$ is a cmx representation of 
$\tilde{M}_t$.
\bb\ni
It follows from Corollary \ref{C 074.01} that $(E,F,G,H)\mapsto 
(\hat{E}, \hat{F}, \hat{G},\hat{H})$, is a one to one 
correspondence between the integrands of $\CS_{\rm cmx}$ and 
the integrands of $\CS\pr$. The one to one correspondence 
between $M_t$  and $\hat{M}_t$ now follows from Proposition \ref{P 74.01} 
and the uniqueness theorem for 
quantum stochastic integrals \cite{Lin2}, \cite[Theorem 
4.7]{Vin}\,\,\epf
\bb\ni
The quantum Ito product formula \cite[Theorem 4.5]{HP1}, 
its extension \cite[Theorem 4]{AM}, and the functional quantum 
Ito formulae \cite[Theorem 4.2, Theorem 6.2.]{Vin2} are true for regular 
quantum semimartingales. Using the isometric isomorphism 
$M\mapsto \tilde{M}$ of $\CS_{\rm cmx}$ and $\CS$ 
in Theorem \ref{T 74.444} the corresponding Ito formulae 
are also valid for regular cmx semimartingales. 
\bb\ni
We will reverse this procedure. 
The functional Ito formula will be proved for a class of process 
$M\in \CS_{\rm cmx}$. The corresponding formula for the 
Hudson--Parthasarathy process $\tilde{M}$ will be deduced 
from the correspondence between quantum 
semimartingales and cmx semimartingales in Theorem \ref{T 74.444}.
\bb\ni
We conclude this section with a sufficient condition that $M$ be 
a regular cmx process.
\begin{Thm}\label{T 74.821} If $E,F,G,H$ are cmx processes such 
that $(N+I)^{1/2}E(N+I)^{1/2}$ is in $L\iiy_{\rm cmx}(\FH)\so$, $F(N+I)^{1/2}$ and  
$(N+I)^{1/2} G$ are in $L^2_{\rm cmx}(\FH)\so$ and $H$ is in $
L^1_{\rm cmx}(\FH)\so$ then  $M(E,F,G,H)$ is in $L\iiy_{\rm cmx}(\FH)\so$.
\end{Thm}
\pf The conditions automatically imply that the cmx quantum 
stochastic integral 
$$
M_t=\int_0^t (E\,d\GL+F\,dA +G\,dA\d+H\,ds)
$$
exists.
It follows from (\ref{E 73.211}) that
$\CS_t(N+I)^{-1/2}$ and $(N+I)^{-1/2}\nabla_t$ 
represent contractions. It follows from (\ref{E 75.204}) 
that the operator
$$
M_t(E:d\GL)=\CS_t E\nabla_t
=\CS_t(N+I)^{-1/2}(N+I)^{1/2} E(N+I)^{1/2}(N+I)^{-1/2}\nabla_t
$$
is bounded in norm by $\leq\nm{(N+I)^{1/2} E(N+I)^{1/2}}\iy$.
\bb\ni
Similarly, if $\psi\in \FH^0\oplus \FH^1\oplus\cdots\oplus\FH^j$ 
for some $j$ then
\begin{eqnarray*}
\nm{M_t(F:dA)\psi}&=&\nm{\int_0^{1} F_s(N+I)^{1/2}(N+I)^{-1/2}
\nabla_t(s)\psi}\,ds
\\ 
&\leq &
\int_0^{1} \nm{F_s(N+I)^{1/2}}\nm{(N+I)^{-1/2}\nabla_t(s)\psi}\,ds 
\\
&\leq &
\nm{F(N+I)^{1/2}}_2\nm{\psi}.
\end{eqnarray*}
Since such $\psi$ are dense in $\FH$ it follows that $M(F:dA)$ is in 
$L\iiy_{\rm cmx}(\FH)\so$.
\bb\ni
If $\psi=\psi^0+\psi^1+\ldots+\psi^k$ with $\psi^j\in \FH^j$ 
then
\begin{eqnarray*}
\nm{M_t(dA\d:G)\psi}^2&=&\nm{\CS_t (N+I)^{1/2}}\nm{(N+I)^{-1/2} G\psi}^2
\\
&=& \sum_{j=0}^k \int_0^t\left(G_s\psi^j(x_1,\ldots,x_j)\right)^2
 \,ds\, dx_1\ldots dx_j
\\
&\leq & \sum_{j=0}^k \int_0^T\nm{G_s}^2|\psi^j(x_1,\ldots,x_j)|^2 
\,ds\, dx_1\ldots dx_j
\\
&\leq & \nm{G}_2^2 \nm{\psi}^2,
\end{eqnarray*}
and $M(d\GL:G)$ is in $L\iiy_{\rm cmx}(\FH)\so$. The remaining case is 
straightforward. \,\,\epf
\bb\ni
It is now easy to find bounded cmx quantum stochastic 
integrals. If $E,F,G,H$ are as in Theorem \ref{T 74.821} and 
$$
\breve{E}=(N+1)^{-\frac12}E(N+1)^{-\frac12},\quad
\breve{F}=F(N+1)^{-\frac12},\quad \breve{G}=(N+1)^{-\frac12}G
$$
then $\breve{E}\in L\iiy_{\rm cmx}(\FH)\so$ and 
$\breve{F},\,\, \breve{G}\in L^2_{\rm cmx}(\FH)\so$. It follows from  
Theorem \ref{T 74.821} that $M(\breve{E},\breve{F},\breve{G},H)$ belongs to 
$L\iiy_{\rm cmx}(\FH)\so$.
\bb\ni
The use of Theorem \ref{T 74.821} in the hunt for 
sufficient conditions on its integrands for a quantum 
semimartingale to be regular may be a chimera.  
This is because the number 
operator is not adapted so that the quadruple 
$(\breve{E},\breve{F},\breve{G}, H)$ will not, normally, be adapted.
Similar difficulties are encountered if the number operator 
$N$ is replaced by the adapted number process 
$\{N_{t]}\otimes I_{(t}$, 
where $N_{t]}$ is the number operator for $\FH_{t]}$.

\section{\label{SMX}Scalar Matrices}
We recall Nelson's analytic vector theorem \cite{Nel}. A proof of 
the theorem may be found in \cite[\S 8.5]{Wei}. If $\hat{T}$ is a 
self-adjoint operator in a Hilbert space $\CH$ then $C\iiy 
(\hat{T})=\bigcap_{k=0}\iiy \CD(\hat{T}^k)$. An element 
$\Gv\in C\iiy(\hat{T})$ is 
an analytic vector of $\hat{T}$ with radius of convergence 
$r(\Gv)$ if $r$ is the radius  of convergence of the complex 
power series
$$ 
\sum_{k=0}\iiy \frac{\|\hat{T}^k \Gv\|}{k!} {z^k}
$$
The set of analytic vectors of $\hat{T}$ is denoted $\CA(\hat{T})$ and 
$\CA_r(\hat{T})=\{\Gv\in \CA(\hat{T}):r(\Gv)\geq r\}$.
\begin{Thm}[Nelson]\label{T 76.01} Let $\hat{T}$ be a symmetric operator on a 
Hilbert space $\CH$ whose domain $\CD(\hat{T})$ 
contains a dense subspace $\CD$ of analytic vectors.  Then 
\bb\ni
{\rm (i)} $\hat{T}$ is essentially self-adjoint.
\bb\ni
{\rm(ii)} $\CD$ is a core for $\hat{T}$.
\bb\ni
{\rm(iii)} If $\Gv$ is an analytic vector of $\hat{T}$ and $p\in \BBC$ 
then
$$
e^{ip\ol{T}}\Gv=\sum_{k=0}\iiy \frac{p^k}{k!}\hat{T}^k \Gv\qquad 
\mbox{ whenever } |p|<r(\Gv)
$$
where $\ol{T}$ is the (self-adjoint) closure of $\hat{T}$.
\end{Thm} 
\begin{Lemma}\label{L 76.01} Let $\Lnu\succ 0$ be a scalar matrix. Then 
$\ell_{00}\subset C\iiy(\tilde{\Lnu})$ if and only if 
$\ell_{00}\subset C\iiy(\Lnu)$.
\bb\ni
In this case $\Lnu^k x=(\tilde{\Lnu})^k x$ whenever $x\in 
\ell_{00}$.
\end{Lemma}
\pf It is enough to consider $x=e^j=(0,0,\ldots,0,1,0,\ldots)^t$, 
$j=0,1,\ldots$, which span $\ell_{00}$.
\bb\ni
Suppose $e^j\in C\iiy(\tilde{\Lnu})$ for all $j$. A simple 
induction shows that, for each $k$,
\beqn\label{E 76.51}
\eta^{(k)}{\,}^i_j=\sum\{\Lnu^{i}_{j_k}\Lnu^{j_k}_{j_{k-1}}\cdots
\Lnu^{j_k}_{j}:(j,j_2,\ldots,j_k)\in \BBN^k\}
\eeqn
is convergent for all $i$ and $j$, and that 
$$ 
(\tilde{\Lnu})^k e^j =
(\eta^{(k)}{\,}^0_j,\eta^{(k)}{\,}^1_j,\eta^{(k)}{\,}^i_j,\ldots)^t\in 
\ell^2.
$$
This implies that $\Lnu^k$ exists for all $k$ and has entries
$\eta^{(k)}{\,}^i_j$.
\bb\ni
Since 
$\eta^{(k+1)}{\,}^i_j=\sum\{\Lnu^i_{j_{k+1}}\eta^{(k)}{\,}^{j_{k+1}}_j:
j_{k+1}=0,1,2,\ldots\}$ it follows that $e^j\in \CD(\Lnu^k)$ for 
all $k$.
\bb\ni
Conversely, if $\Lnu^k$ exists for all $k$ the formula (\ref{E 
76.51}) defines $\eta^{(k)}{\,}^i_j$ for all $i,j$ and $k$. If 
$e^j\in C\iiy(\Lnu)$ then, for all $k$, 
$$ 
{\Lnu}^k e^j =
(\eta^{(k)}{\,}^0_j,\eta^{(k)}{\,}^1_j,\eta^{(k)}{\,}^i_j,\ldots)^t\in 
\ell^2.
$$
Since 
$\eta^{(k+1)}{\,}^i_j=\sum\{\Lnu^i_{j_{k+1}}\eta^{(k)}{\,}^{j_{k+1}}_j:
j_{k+1}=0,1,2,\ldots\}$ it follows that $\Lnu^{k+1} 
e^j=\tilde{\Lnu}(\Lnu^k e^j)$ and $e^j\in C\iiy(\tilde{\Lnu})$ 
for all $j$. \,\,\epf
\begin{Propn}\label{P 76.02} Let $\Xi\succ 0$ be a (2k+1)-diagonal scalar 
matrix, let $\xi>0$ and let $\Ge=\xi\inv(4k^2+2k)\inv$. Then
\bb\ni
{\rm (i)} If $\Xi^i_j\leq \xi(i+j)$ for all $i,j$ then 
$\ell_{00}\subset \CA_{\Ge}(\Xi)$;
\bb\ni
{\rm (ii)} If $\Xi^i_j\leq \xi(i^{1/2}+j^{1/2})$ for all $i,j$ then 
$\ell_{00}\subset \CA\iy(\Xi)$
\end{Propn}
\pf (i) It is enough to consider the case $\xi=1$. 
A simple induction shows that $x\ub{n}=\Xi^n \Ge_j\in \ell_{00}$ 
with $x\ub{n}_r=0$ for $r>j+nk$. Assume, inductively, that 
\beqn\label{E 79.131}
x\ub{n-1}_r &\leq& (2j+k)(2j+3k)\cdots(2j+(2n+1)k)(2k+1)^{n-1}
\eeqn
for all $r$. Now $x\ub{n}=\Xi x\ub{n-1}$ and, since $\Xi$ is $2k+1$ 
diagonal, the 
largest entry in $\Xi$ which can multiply a non-zero element of 
$x\ub{n-1}$ is $\Xi^{j+nk}_{j+(n-1)k}\leq 2j+(2n+3)k$. 
Since $\Xi$ is $2k+1$ diagonal at most $2k+1$ of the entries of 
$x\ub{n-1}$ contribute to each entry of $x\ub{n}$. Therefore
$$
x\ub{n}_r \leq (2j+k)(2j+3k)\cdots(2j+(2n+3)k)(2k+1)^{n}
$$
$$
\frac{\nm{\Xi^n \Ge_j}}{n!}=\frac{\|x\ub{n}\|}{n!}\leq 
\frac{(2j+k)(2j+3k)\cdots(2j+(2n+3)k)(2k+1)^n(j+nk)^{\frac12}}{n!}=c_n.
$$
Now 
$$
\lim_{n\roy}\frac{c_{n+1}}{c_n}= 
\lim_{n\roy}\frac{(2j+(2n+5)k)(2k+1)(j+(n+1)k)^{\frac12}}{(n+1)(j+nk)^{\frac12}}
=2k(2k+1).
$$
By the limit ratio test $\Ge_j$ is an analytic vector for $\Xi$ with 
$r(\Ge_j) =1/(4k^2+2k)=\Ge$ for all $j$. Therefore $\ell_{00}\subset 
\CA_{\Ge}(\Xi)$.
\bb\ni
(ii) Using the same argument as in (i) with the bound
$\Xi^{j+nk}_{j+(n-1)k}\leq 2(j+(n+1)k)^{1/2}$ gives that
$$
x\ub{n-1}_r \leq 2^{n-
1}[(j+k)(j+2k)\cdots(j+nk)]^{\frac12}(2k+1)^{n-1}.
$$
The limit ratio test then shows that the power series 
$\sum_{n=0}\iiy z^n\nm{\Xi^n \Ge_j}/n!$ has infinite 
radius of convergence and $\Ge_j\in \CA\iy(\Xi)$. Therefore 
$\ell_{00}\subset \CA\iy(\Xi)$.\,\,\epf
\bb\ni
The following proposition follows directly from the definition of 
$\LGk$. 
\begin{Propn}\label{PP 76.01}
Let $M$ be a cmx semimartingale, with control matrix $\LGk$, 
whose integrands $E,F,G,H$ are $(2k+1)$-diagonal matrices. 
Suppose there exists a constant $\xi$ such that
\beqn\label{E 78.978}
\nm{E^i_j}\iy,\,\,\,\,\nm{F^i_j}_2,\,\,\,\,\nm{G^i_j}_2,
\,\,\,\,\nm{H^i_j}_1\leq \xi
\eeqn
for all $i,j=0,1,2,\ldots$. If $\Ge=\xi\inv(4(k+1)^2+2(k+1))\inv$ then 
\bb\ni
(i) $\ell_{00}\subset \CA_{\Ge}(\LGk)$.
\bb\ni
(ii) If $E=0$ then  $\ell_{00}\subset \CA\iy(\LGk)$.
\end{Propn}
If $K$ is a $(2k+1)$-diagonal process defined by (\ref{E 074.02})
and with $\nm{k_t{}\,^i_j}\leq \xi$ for all $t\in [0,1]$ then 
$K$ satisfies (\ref{E 78.978}). Thus there is a large class of cmx semimartingales whose control matrices   
satisfy the conclusions of Proposition \ref{P 76.02}. These control 
matrices represent essentially self-adjoint operators.
\begin{Propn}\label{P 76.01} Let $\Lnu$ be a scalar matrix with 
$\Lnu\succ 0$, $\Lnu=\Lnu\st$. 
\bb\ni
{\rm (i)} If $\ell_{00}\subset \CA(\Lnu)$ then $\tilde{\Lnu}$ is 
essentially self-adjoint, with closure $\ol{\Lnu}$.
\bb\ni 
{\rm (ii)}  If $\ell_{00}\subset \CA_1(\Lnu)$ and $-1\leq \rho\leq 1$ then 
$\sum_{k=0}\iiy (i\rho \Lnu)^k/k!$ converges entrywise to the 
unique bounded scalar matrix $e^{i\rho \Lnu}$ representing 
$e^{i\rho \ol{\Lnu}}$.
\end{Propn}
\pf (i) By Lemma \ref{L 76.01} if $x\in \ell_{00}$ then $x\in 
\CA(\tilde{\Lnu})$ and the two series
\beqn\label{E 76.02}
\sum_{k=0}\iiy \frac{i\rho^k}{k!}(\tilde{\Lnu})^k x &=&
\sum_{k=0}\iiy \frac{i\rho^k}{k!} \Lnu^k x
\eeqn
are identical and absolutely convergent for $|\rho|\leq r(x)$. 
Therefore $\ell_{00}\subset \CA(\tilde{\Lnu})$. Since 
$\ell_{00}$ is dense in $\ell^2$ it follows from Nelson's 
analytic vector theorem that
$\tilde{\Lnu}$ is essentially self-adjoint. It 
also follows that the sum of the series (\ref{E 76.02}) is 
$e^{i\rho \ol{\Lnu}}x$ for each $x\in \ell_{00}$
\bb\ni
(ii) $\Lnu^k e^j$ is the $j$th column vector in $\Lnu^k$. Taking 
$x=e^j$, $j=0,1,2,\ldots$ in (\ref{E 76.02}) shows that the 
columns of $\sum_{k=0}\iiy (i\rho \Lnu)^k/k!$ sum in $\ell^2$, 
and {\em a fortiori} coordinatewise to the corresponding column
of the matrix representing the bounded linear transformation 
$e^{i\rho \ol{\Lnu}}$.\,\,\epf
\bb\ni
The {\it scalar matrix} of $X\in \CL^p_{\rm cmx}(\FH)\so$ is the real 
matrix
$\Lnu(X)$ with entries
$$
\Lnu^i_j(X)=\nm{ X^i_j}_p,\qquad i,j=0,1,2,\ldots.
$$
Processes in $\CL^p_{\rm cmx}(\FH)\so$ are controlled by their scalar matrices. The 
proof of the following lemma is an immediate consequence of the definitions.
\begin{Lemma}\label{L 78.019} Let $S_n$ be a sequence in $\CL^p_{\rm cmx}(\FH)\so$.
\bb\ni
(i) $S_n$ is convergent to $S$ in 
$\CL^p_{\rm cmx}(\FH)\so$ if and only if $\Lnu(S-S_n)$ converges to the zero matrix.
\bb\ni
(ii) $S_n$ is absolutely summable if and only if $\Lnu(S_n)$ is summable.
\end{Lemma}
\begin{Lemma}\label{L 78.119} 
Let $p,q,r\in \{1,2,\infty\}$ with $1/p+1/q=1/r$ and let $0\prec \Phi$ and 
$0\prec\Psi$ be scalar matrices such that the 
product matrix $\Phi\Psi$ exists.
If $S\in \CL^p_{\rm cmx}(\FH)\so$ and $T\in \CL^q_{\rm cmx}(\FH)\so$ are such 
that $\Lnu(S)\prec \Phi$ and $\Lnu(T)\prec \Phi$ then the product process 
$ST$ exists in $\CL^r_{\rm cmx}(\FH)\so$ and $\Lnu(ST)\prec \Phi\Psi$.
\end{Lemma}
\pf The product $S^i_k T^k_j$ exists in $L^r([0,1],\CB(\FH^j,\FH^i))$ 
and $\nm{S^i_k T^k_j}_r\leq \Lnu(S)^i_k\Lnu(T)^k_j$. Therefore 
$$
\sum_{k=0}\iiy \nm{S^i_k T^k_j}_r
\leq \sum_{k=0}\iiy \Lnu(S)^i_k\Lnu(T)^k_j=(\Lnu(S)\Lnu(T))^i_j\leq
(\Phi\Psi)^i_j<\infty
$$
and the sequence $\{S^i_k T^k_j:k=0,1,\ldots\}$ is absolutely summable and 
therefore summable in $L^r([0,1],\CB(\FH^j,\FH^i))$. Therefore 
the product process $ST$ exists in $\CL^r_{\rm cmx}(\FH)\so$ and is dominated by 
$\Phi\Psi$.\,\,\epf
\bb\ni
If $\Gv\in \FH$ define  $\nu(\Gv)\in \ell^2$ by the formula
$$
\nu(\Gv)=(\|\Gv^0\|,\|\Gv^1\|,\ldots,\|\Gv^k\|,\ldots)^t
$$
\begin{Lemma}\label{L 76.03} If $M\in \CL\iiy_{\rm cmx}(\FH)\so$ and 
$\Lnu(M)^k$ exists 
then $M_t^k$ exists and $\widetilde{M_t^k}=(\tilde{M_t})^k$ for 
almost all $t\in [0,1]$. 
\bb\ni
If $\psi\in \FD(M_t^k)$ then $\nu(M_t^k \psi)\prec 
\Lnu(M)^k\nu(\psi)$ for almost all $t$.
\end{Lemma}
\pf It follows from Lemma \ref{L 78.119} that $M_t^k$ exists. 
Let $\Lnu=\Lnu(M)$. If $\psi\in \FH$ then
$$
\nm{\sum_{j=0}\iiy M_t{\,}^i_j \psi^j} \leq 
\sum_{j=0}\iiy \nm{M_t{\,}^i_j}\nm{\psi^j}
\leq 
\sum_{j=0}\iiy \nm{\Lnu^i_j}\nm{\psi^j},
$$
and $\nu(M_t \psi)\prec \Lnu \nu(\psi)$. This may be iterated to 
give $\nu(M_t^k \psi)\prec \Lnu(M)^k\nu(\psi)$.\,\,\epf
\begin{Lemma}\label{L 76.02} If $X\in \CL\iiy_{\rm cmx}(\FH)\so$ 
and $\psi\in \FH$ with 
$\nu(\psi)\in \CA(\Lnu(X))$ then $\psi\in \CA(X_t)$ for almost 
all $t\in [0,1]$.
\end{Lemma}
\pf Let $\Lnu=\Lnu(X)$. Since $C\iiy(\Lnu)\neq \emptyset$ it 
follows by  Lemma \ref{L 78.119} that 
$\Lnu^k$ and therefore $X^k$ each exist for all $k$. Since 
$\Lnu_t(X^k)\prec \Lnu^k$ for almost all $t$ it follows that $\psi\in 
\CA(X_t)$ for almost all $t$.\,\,\epf
\begin{Thm}\label{T 76.02}
Let $M$ be symmetric cmx process in $\CL\iiy_{\rm cmx}(\FH)\so$ 
with scalar matrix $\Lnu=\Lnu(M)$. If $\ell_{00} \subset 
\CA(\Lnu(M))$ then 
\bb\ni
{\rm(i)} 
$\tilde{M}$ is essentially self-adjoint with core 
$\FH_{00}\subset \CA(M_t)\cap\CA(\tilde{M_t})$ for almost all $t\in [0,1]$
\bb\ni
{\rm(ii)}
If $\ell_{00} \subset \CA_1(\Lnu(M))$ and $\rho\in [-1,1]$ 
then $\sum_{k=0}\iiy (i\rho M)^k/k!$ is convergent in 
$\CL\iiy_{\rm cmx}(\FH)\so$ to a process 
$J$.         
\bb\ni
For $t\in [0,1]$ let $\ol{M}_t$ be  the closure of 
$\tilde{M}_t$.  Then $J_t$ is the cmx representation of 
the unitary operator $e^{i\rho \ol{M}_t}$.
\end{Thm}
\pf By definition $\Lnu^k$ exists and by Lemma \ref{L 76.03} that 
$M_t^k$ exists for all $k$. 
\bb\ni
(i) If $\psi\in \FH_{00}$ then $\nu(\psi)\in \ell_{00}$ and
\beqn\label{E 76.03}
\sum_{n=0}\iiy \frac{\rho^n}{n!}\nm{M_t^n \psi}\leq \sum_{n=0}\iiy 
\frac{\rho^n}{n!}\nm{\Lnu^n \nu(\psi)}<\infty
\eeqn
whenever $|\rho|<r(\nu(\psi))$. By Lemma \ref{L 76.03} this 
remains true if $M_t$ replaced by $\tilde{M_t}$.
Therefore  $\FH_{00}\subset \CA(M_t)\cap\CA(\tilde{M}_t)$ 
for all $t\in [0,1]$ and (i) follows from Theorem \ref{T 76.01}.
\bb\ni
\bb\ni
(ii) It follows from Proposition \ref{P 76.01} that  $\sum_{k=0}\iiy 
(\Lnu)^k/k!$ is entrywise convergent to a scalar matrix. By 
Lemma \ref{L 78.119}  
$\Lnu(M^k)\prec \Lnu(M)^k= \Lnu^k$ and,  by Lemma \ref{L 78.019}, 
$\sum_{k=0}\iiy (i\rho M)^k/k!$ is convergent to a process 
$J$ in $\CL\iiy_{\rm cmx}(\FH)\so$.
\bb\ni
If $\psi\in \FH_{00}$ and $|\rho|\leq 1$ then by Lemma \ref{L 76.03}
$$
J_t\psi
=\sum_{n=0}\iiy \frac{i\rho^n}{n!}{M_t^n \psi}
=\sum_{n=0}\iiy \frac{i\rho^n}{n!}{\widetilde{M_t^n} \psi}
=\sum_{n=0}\iiy \frac{i\rho^n}{n!}{(\tilde{M}_t)^n \psi}
= \sum_{n=0}\iiy \frac{i\rho^n}{n!}{\ol{M_t}^n \psi}
=e^{i\rho \ol{M_t}}\psi.
$$
It follows from Theorem \ref{T 73.001} (iii) that $J_t$ is the unique 
cmx representation of $e^{i\rho \ol{M}_t}$
\,\,\epf
\bb\ni
\begin{Propn}\label{P 78.311}
Let $p,q,r\in \{1,2,\infty\}$ with $1/p+1/q=1/r$ and let $0\prec \Phi(u)$ and 
$0\prec\Psi(u)$ be continuous scalar matrix valued functions 
such that the product matrix $\Phi(u)\Psi(u)$ exists for $u\in [0,1]$.
\bb\ni
Let $u\mapsto S_n(u)$ and $u\mapsto T_n(u)$ be sequences of 
continuous functions $S_n:[0,1]\row \CL^p_{\rm cmx}(\FH)\so$ and 
$T_n:[0,1]\row \CL^q_{\rm cmx}(\FH)\so$ respectively such that, for each $u\in [0,1]$,
\bb\ni
{\rm(a)} $S_n(u)$ converges to $S(u)$ in $\CL^p_{\rm cmx}(\FH)\so$ 
and $\Lnu(S_n)\prec \Phi(u)$ $n=1,2,\ldots$;
\bb\ni
{\rm(b)} $T_n(u)$ converges to $T(u)$ in $\CL^q_{\rm cmx}(\FH)\so$ 
and $\Lnu(T_n)\prec \Psi(u) $ $n=1,2,\ldots$.
\bb\ni
Then $S_n(u)T_n(u)$ and $S(u)T(u)$ both exist and 
$S_n(u)T_n(u)$ converges to $S(u)T(u)$ in $\CL^r_{\rm cmx}(\FH)\so$ 
for each $u\in 
[0,1]$. Moreover and $u\mapsto S(u)T(u)$ is entrywise Bochner--Lebesgue integrable and 
$$
\lim_{n\roy}\int_0^1 S_n(u)T_n(u)\,du=\int_0^1 S(u)T(u)\,du.
$$
\end{Propn}
\bb\ni
\pf It follows from Lemma \ref{L 78.119} that  
$S_n(u)T_n(u)$ and $S(u)T(u)$ both exist for all $u\in  [0,1]$.
\bb\ni
Suppose $\Phi_n$ is a sequence of scalar matrices with $0\prec\Phi_n\prec 
\Phi$ such that $\lim_{n\roy}\Phi_n$is the zero matrix. Then 
$(\Phi_n)^i_k\Psi^k_j\leq (\Phi_n)^i_k\Psi^k_j$ for each $k$ for all $n$. 
Now $\sum_{k=0}\iiy (\Phi_n)^i_k\Psi^k_j=(\Phi\Psi)^i_j<\infty$ and 
by a series version of the dominated convergence theorem
$$
\lim_{n\roy}(\Phi_n \Psi)^i_j=\lim_{n\roy}\sum_{k=0}\iiy 
(\Phi_n)^i_k\Psi^k_j=0
$$
for all $i,j$. Since $\Lnu(S(u)-S_n(u))\prec 2\Phi(u)$ 
it follows that  $\lim_{n\roy}\Lnu(S(u)-S_n(u))\Lnu(T(u))=0$. Similarly, since 
$\Lnu(S_n(u))\Lnu(T(u)-T_n(u))\prec \Phi\Lnu(T(u)-T_n(u))$ it follows that 
$\lim_{n\roy}(S_n(u))\Lnu(T(u)-T_n(u))=0$. Since
$$
\Lnu(S(u)T(u)-S_n(u)T_n(u))\prec  \Lnu(S(u)-S_n(u))\Lnu(T(u))
+\Lnu(S_n(u))\Lnu(T(u)-T_n(u))
$$
it follows that $\lim_{n\roy}\Lnu(S(u)T(u)-S_n(u)T_n(u))=0$. By Lemma
\ref{L 78.019} $S_n(u)T_n(u)$ converges to $S(u)T(u)$ in $\CL^r_{\rm cmx}(\FH)\so$ 
for each $u\in [0,1]$.
\bb\ni
Since  $S_n(u)T_n(u)\prec \Phi(u)\Psi(u)$ for all $n$ and $u\in [0,1]$ 
it follows from the bounded convergence theorem that 
$u\mapsto S(u)T(u)$ is entrywise Bochner--Lebesgue integrable and 
$$
\lim_{n\roy}\int_0^1 S_n(u)T_n(u)\,du=\int_0^1 S(u)T(u)\,du.\,\,\epf
$$

\section{\label{HPP}Hudson--Parthasarathy Processes}
Theorem \ref{T 74.444} exhibits the natural correspondence between 
cmx semimartingales and quantum semimartingales: 
processes with regular integrands.
\bb\ni
In the Hudson--Parthasarathy theory 
there are interesting quantum stochastic integrals
whose integrands are processes of unbounded operators. 
Polynomials in the Brownian semimartingale $(A_A\d)$, for 
example, do not satisfy the conditions of Theorem \ref{T 74.444}. A 
polynomial 
$p(\GL_t,A_t,A_t\d)$, of degree greater than one, whose 
coefficients are continuously differentiable functions of $t$ 
may be written as a quantum stochastic integral whose integrands are 
processes of unbounded operators. For example 
$$
{A}_t^2=2\int_0^t {A}_s\,d{A}_s.
$$
In this section we show that the correspondence in 
Theorem \ref{T 74.444} may be extended to such processes.
\bb\ni
The integrands in the Hudson--Parthasarathy theory are 
adapted processes of operators with common domain 
$\CE_S$ where $S$ is an admissible 
subspace of $L\iiy[0,1]$ \cite[\S 3]{HP1}, \cite[\S 24]{P1}. 
The restriction that $S\subset L\iiy[0,1]$ imposed in \cite{HP1} is 
unnecessary: {\em cf}. \cite[\S 25]{P1} and \cite[VI. 6]{MEY}. 
\bb\ni
Measure theoretic considerations 
require a slight modification of some definitions in \cite{HP1},
as in \cite[\S 25]{P1}. A process 
$\hat{X}=\{\hat{X}_t:t\in [0,1]\}$ is {\it measurable} with 
respect to $S$ if $t\mapsto \hat{X}_t e(g)$ is Borel measurable, for each 
$g\in S$. 
\bb\ni
This definition requires only that $e(g)\in \CD(X_t)$ for almost 
all $t\in [0,1]$.
\bb\ni
$\hat{X}$ is {\it adapted} with respect to $S$ if, 
for each $t\in [0,1]$, there exists a linear 
transformation $\hat{X}_{t]}$ in $\FH_{t]}$ such that
$\hat{X}_t e(g)=\FU_t[(\hat{X}_t e(g_{t]}))\otimes e(g_{(t})]$ 
for almost all $t\in [0,1]$.
\bb\ni
For the rest of this article, unless otherwise stated, the 
integrands of quantum stochastic integrals and cmx quantum 
stochastic integrals will be adapted.
\bb\ni
Let $p=1,2$ or $\infty$ and define the following vector spaces of 
operator processes in $\FH$
\bb\ni
$\FM(\CE_S)=\{\hat{X}: \hat{X} \mbox{ is measurable with respect to $S$}\}$;
\bb\ni
$\FL^p(\CE_S)=\{\hat{X}\in \FM(\CE_S):t\mapsto \hat{X}_t e(g) \mbox{ is in }
L^p([0,1], \FH) \mbox{ for each }g\in S\}$; 
\bb\ni
For simplicity we let $S=L^2[0,1]$ and put $\CE=\CE_S$.
\bb\ni
The usual measure theoretic conventions apply. If $\hat{X}\in 
\FL^p(\CE)$ and $g\in L^2[0,1]$ then $\hat{X}_s e(g)$ need only be defined 
almost everywhere on $[0,1]$. If $\hat{Y}\in \FL^p(\CE)$ and
$\hat{Y}_t e(g)=\hat{X}_t e(g)$ for almost all $t$, whenever 
$g\in S$, then the processes $\hat{Y}$ and $\hat{X}$ are 
equivalent, and will normally be identified. If $g\in L^2[0,1]$ define
$$
\check{X}_t e(g)=\left\{ \ba{rl} 
\hat{X}_t e(g) &\qquad\mbox{ if } e(g)\in \CD(\hat{X}_t)  \\
0 & \qquad \mbox{ otherwise }
\ea \right. \quad \forall \,t\in [0,1].
$$
Then $\check{X}$ is equivalent to $\hat{X}$ and $ 
\CD(\check{X}_t)=\CE$ for all $t\in [0,t]$. If $\hat{X}$ is an 
adapted process in $ \FL^p(\CE)$ 
then $\check{X}$ is a process of the type considered by
Hudson and Parthasarathy \cite[\S 3]{HP1}. 
\bb\ni
If $\hat{E}\in \FL\iiy(\CE)$, $\hat{F}, \hat{G}\in \FL^2(\CE)$
and $\hat{H}\in \FL^1(\CE)$ are all adapted 
it follows from \cite[Theorem 4.4]{HP1} 
that there is a unique process 
$\check{M}=\check{M}( \check{E},\check{F},\check{G},\check{H})$ with 
$\CD(\check{M}_t)=\CE$ for all $t\in [0,1]$:
$$
\check{ M}_t=\int_0^t{\,} \check{ E}\,d{\GL}
+\check{ F}\,d{ A}+\check{ G}\,d{ A}\d+\check{ H}\,d{ s}
$$
such that
\begin{eqnarray}
\nonumber \ip{\check{M}_t e(f)}{e(g)} &=& \int_0^t\la  
(f\ol{g}\check{E}+f\check{F}+\ol{g}\check{G}
+\check{H})(s)e(f),e(g)\ra\,ds\\
\label{E 77.31} &=& \int_0^t\la  
(f\ol{g}\hat{E}+f\hat{F}+\ol{g}\hat{G}
+\hat{H})(s)e(f),e(g)\ra\,ds
\end{eqnarray}
for all $f,g\in L^2[0,1]$. The {\it quantum stochastic 
integral} $\hat{M}=\{\hat{M}_t:t\in [0,1]\}$
of the quadruple $(\hat{E},\hat{F},\hat{G},\hat{H})$ is defined 
to be the process $\check{M}=\{\check{M}_t:t\in [0,1]\}$.
It follows from (\ref{E 77.31}) that this definition is 
consistent with that of Hudson and Parthasarathy \cite[Theorem 
4.4]{HP1} and depends 
only on the equivalence classes of $\hat{E}$, $\hat{F}$, $\hat{G}$ and 
$\hat{H}$.
\bb\ni
A cmx process $X\in \CL^p_{\rm cmx}(\FH)\so$ is a 
{\it cmx representation} of 
$\hat{X}\in \FL^p(\CE)$ if for each $g\in L^2[0,1]$ there 
exists a null set $\CN_g$ such that $e(g)\in \FD(X_t)$
and $X_t e(g)$ is the chaos representation of $\hat{X}_te(g)$ 
whenever $t\in [0,1]\setminus \CN_g$. 
\bb\ni
We give a sufficient condition that this be the case.
\bb\ni
The {\it exponential domain}, $\ell^2_{\rm exp}$, in $\ell^2$ is the 
linear span of vectors of the form
$$
\Gve(x)=(1,x,{x^2}/{(2!)^{1/2}},\ldots,{x^j}/{(j!)^{1/2}},\ldots), \quad x\in 
\BBC.
$$
\begin{Thm} \label{T 77.801} 
Let $p=1,2$ or $\infty$. If $X\in \CL^p_{\rm cmx}(\FH)\so$ and 
$\ell^2_{\rm exp}\subset \CD(\Lnu(T))$ then 
\bb\ni
{\rm (i)} $\tilde{X} \in \FL^p(\CE)$;
\bb\ni
{\rm (ii)} If $X$ is adapted so is $\tilde{X} \in L^p(\CE)$.
\end{Thm}
\def\Lnu{\mbox{\Large$\nu$}}
\pf
(i) $p=\infty$, $\Lnu^i_j=\nm{X^i_j}\iy$. 
$$
\sum_{j=0}\iiy \nm{X_s{\,}^i_j e^j(g)}   \leq   
\sum_{j=0}\iiy \nm{X_s{\,}^i_j}\iy\nm{ e^j(g)}
 \leq \sum_{j=0}\iiy \Lnu^i_j \nm{ e^j(g)}
 = (\Lnu\Gve(\nm{g}))^i,
$$
and $\sum X_s{\,}^i_j e^i(g)$ is absolutely convergent in $\FH^i$ to 
$\eta_s^i(g)$ with 
$\nm{\eta_s^i}\leq \nm{( \Lnu\Gve(\nm{g}))^i}$ for almost all $s\in 
[0,1]$. Therefore 
$$
\sum_{i=0}\iiy \nm{\eta_s^i(g)}^2\leq \sum_{i=0}\iiy 
\nm{( \Lnu\Gve(\nm{g}))^i}^2 = \nm{\Lnu\Gve(\nm{g})}^2,
$$
and $\sum_{i=0}\iiy \eta_s^i(g)$ is convergent, for almost all $s$, to 
$\eta_s(g)$ in $\FH$. For each $s\in [0,1]$ and $g\in L^2[0,1]$ define 
$$
{\hat{X}}_s e(g)=\left\{ \ba{ll} \eta_s(g)&\qquad\mbox{ whenever 
$\sum_{i=0}\iiy \eta_s^i(g)$ is convergent }
\\
0&\qquad\mbox{ otherwise}.
\ea\right.
$$                     
Since the exponential vectors are linearly independent $\hat{X}_s e(g)$ 
is well defined for each $g\in L^2[0,1]$ and $\hat{X}_s$ 
extends uniquely to a linear operator ${\hat{ X}}_s:\CE\row \FH$. 
Since $t\mapsto X_t{\,}^i_j e^i(g)$ is measurable it is clear 
from above that 
$t\mapsto \hat{X}_t e(g)$ is also measurable for each $g\in 
L^2[0,1]$. Moreover $\nm{\hat{X}_s e(g)}\leq \nm{\Lnu\Ge(\nm{g})}$ for 
almost all $s$ in $[0,1]$  and $\hat{X}\in \FL\iiy(\CE)$.
\bb\ni
From the definition of $\tilde{X_t}$ in Section \ref{CMX} 
it follows that $\hat{X}=\tilde{X}$ and 
$\tilde{X} \in \FL\iiy(\CE)$.
\bb\ni
(ii) $p=2$, $\Lnu^i_j=\nm{X^i_j}_2$. 
If $g\in L^2[0,1]$ then
\begin{eqnarray}
\nonumber \nm{\Lnu\Gve(\nm{g})}^2 &\geq&
\sum_{i=0}\iiy \left(\sum_{j=0}\iiy \nm{X^i_j}_2\nm{e^j(g)}\right)^2
\\
\nonumber &=&
\sum_{i,j,k=0}\iiy \nm{X^i_j}_2\nm{X^i_k}_2 \nm{e^j(g)}\nm{e^k(g)}
\\
\nonumber &\geq&
\sum_{i,j,k=0}\iiy \int_0^t\nm{X_s{\,}^i_j}\nm{X_s{\,}^i_k}\,ds 
\nm{e^j(g)}\nm{e^k(g)}
\\
\label{AA 79.01}&=&
\sum_{i=0}\iiy \int_0^t\sum_{j,k=0}\iiy \nm{X_s{\,}^i_j}\nm{X_s{\,}^i_k}\,ds 
\nm{e^j(g)}\nm{e^k(g)}
\\
\nonumber &\geq&
\sum_{i=0}\iiy \int_0^t\left|\ip{\sum_{j=0}\iiy  X_s{\,}^i_j e^j(g)}
{\sum_{k=0}\iiy  X_s{\,}^i_k e^k(g)}\right|\,ds
\\
\nonumber &\geq&
\sum_{i=0}\iiy \int_0^t\nm{\sum_{j=0}\iiy  X_s{\,}^i_j e^j(g)}^2\,ds
\\
\nonumber &=&
\int_0^t\sum_{i=0}\iiy \nm{\sum_{j=0}\iiy  X_s{\,}^i_j e^j(g)}^2\,ds.
\end{eqnarray}
It follows from the Lebesgue series theorem that, for almost all $s\in 
[0,1]$, the sum $\sum_{j=0}\iiy X_s{\,}^i_j e^j(g)$ is absolutely 
convergent in 
$\FH^j$ to $\eta_s^j$ and the direct sum $\bigoplus_{i=0}\iiy \eta_s^j$ is 
convergent in $\FH$ to $\eta_s$. Furthermore the above inequalities show that 
$s\mapsto \nm{\eta_s}$ is in $L^2[0,1]$. 
For each $s\in [0,1]$ and $g\in L^2[0,1]$ define 
$$
{\hat{ X}}_s e(g)=\left\{ \ba{ll} \eta_s(g)&\qquad\mbox{ whenever 
$\sum_{i=0}\iiy \eta_s^i(g)$ is convergent }
\\
0&\qquad\mbox{ otherwise}.
\ea\right.
$$
As in (i) ${\hat{ X}}_s$ extends uniquely to a linear 
operator ${\hat{ X}}_s:\CE\row \FH$ such that 
$s\mapsto \hat{X}e(g)$ is measurable for all $g\in L^2[0,1]$.
\bb\ni
For almost all $s$ 
$$
\nm{\hat{X}_s e(g)}^2=\sum_{i=0}\iiy \nm{\sum_{j=0}\iiy X_s{\,}^i_j 
e^j(g)}^2.
$$
Since $\Lnu$ is defined on the exponential domain, 
it follows from the inequalities (\ref{AA 79.01}) that
$$
\int_0^t\nm{\hat{ X}_s e(g)}^2\,ds=
\int_0^t\sum_{i=0}\iiy \nm{\sum_{j=0}\iiy  X_s{\,}^i_j e^j(g)}^2\,ds
\leq \nm{\Lnu \Gve(\nm{g})}^2<\infty
$$
and $\hat{X}$ belongs to $\FL^2(\CE)$. As in (i) $\hat{X}=\tilde{X}$ 
and $\tilde{X} \in \FL^2(\CE)$.
\bb\ni
(iii) $p=1$, $\Lnu^i_j=\nm{X^i_j}_1$. If $g\in 
L^2[0,1]$ then
\beqn\label{E 74.05}
\sum_{j=0}\iiy \int_0^{\tau}\nm{X_s{\,}^i_j e^j(g)}\,ds
\leq
\sum_{j=0}\iiy \Lnu^i_j \nm{e^j(g)}<\infty
\eeqn
and $\sum_{j=0}\iiy X_s{\,}^i_j e^j(g)$ is almost everywhere absolutely 
convergent to 
$\eta^i_s(g)$ in $\FH^i$. If $\psi\in \FH$ then
$$
\sum_{i=0}\iiy \nm{\Gv^i}\int_0^{1} \nm{\eta^i_s(g)}\,ds
\leq 
\sum_{i,j=0}\iiy \nm{\Gv^i}\Lnu^i_j \nm{e^j(g)}=
\sum_{i=0}\iiy \nm{\Gv^i}\sum_{i=0}\iiy\Lnu^i_j \nm{e^j(g)}
<\infty
$$
since $(\nm{\Gv^0},\nm{\Gv^1},\nm{\Gv^2},\ldots)^t\in \ell^2$. 
Therefore, for each $x\in \ell^2$, the sum 
$\sum_{i=0}\iiy x^i \int_0^{1} \nm{\eta^i_s(g)}\,ds$ is convergent. It 
follows from the uniform boundedness theorem and Lebesgue's dominated 
convergence theorem that
$$
\sum_{i=0}\iiy \int_0^{1} \nm{\eta^i_s(g)}\,ds<\infty\quad\mbox{ and }\quad
\sum_{i=0}\iiy \nm{\eta^i_s(g)}<\infty
$$
almost everywhere. Therefore $\sum_{i=0}\iiy \nm{\eta^i_s(g)}^2<\infty$ almost 
everywhere. 
For each $s\in [0,1]$ and $g\in L^2[0,1]$ define 
$$
\hat{ X}_s e(g)=\left\{ \ba{ll} \eta_s(g)&\qquad\mbox{ whenever 
$\sum_{i=0}\iiy \eta_s^i(g)$ is convergent }
\\
0&\qquad\mbox{ otherwise}.
\ea\right.
$$
As in (ii) this extends uniquely to a linear 
operator $\hat{ X}_s:\CE\row \FH$.
\bb\ni
It follows from (\ref{E 74.05}) that
$$
\int_0^{1} \nm{\tilde{ X}_s e(g)}\,ds =
\int_0^{1}\nm{\bigoplus_{i=0}\iiy \sum_{j=0}\iiy X_s{\,}^i_j e^j(g)}\,ds
\leq
\int_0^{1}\sum_{i=0}\iiy \nm{\sum_{j=0}\iiy X_s{\,}^i_j e^j(g)}\,ds
<\infty
$$
so that the process $\tilde{ X}$ is in 
$\FL^1(\CE)$.  From the definition of $\tilde{X_t}$ in Section \ref{CMX} 
it follows that $\hat{X}=\tilde{X}$ and 
$\tilde{X} \in \FL^1(\CE)$.
\bb\ni
(ii) This follows from Theorem \ref{P 1.001} (ii).\,\,\epf 
\begin{Thm}\label{T 74.001} Let $(E,F,G,H)$ be an adapted 
integrable cmx quadruple such that $\ell^2_{\rm exp}\subset 
\CD(\Lnu)$ whenever $\Lnu$ is one of the scalar matrices 
$\Lnu(E)$, $\Lnu(F)$, $\Lnu(G)$, $\Lnu(H)$.
\bb\ni
If $0\leq t\leq 1$ the quantum stochastic integral 
\beqn \label{E 77.800}
\hat{M}_t =\int_0^t{\,} (\tilde{ E}\,d{\GL} +\tilde{ F}\,d{ A}
+\tilde{ G}\,d{ A}\d+\tilde{ H}\,d{s})
\eeqn
exists with $\CD(\hat{M}_t)=\CE$. Moreover, if 
$$
M_t =\int_0^t E_s\,d\GL_s+F_s\,dA_s+G_s\,dA\d_s+H_s\,ds
$$
then $\CE\subset \FD(M_t)$ and $M$ is a cmx representation of 
$\hat{M}$.
\end{Thm}
\pf Let $S=L^2[0,1]$. 
From the remarks preceding Theorem \ref{T 77.801} 
and its Corollary it follows that that the quantum stochastic integral 
(\ref{E 77.800}) exists and is the unique process such that
\beqn
\ip{\hat{M}_t e(f)}{e(g)} &=& \int_0^t\la  
(f\ol{g}\tilde{E}+f\tilde{F}+\ol{g}\tilde{G}
+\tilde{H})(s)e(f),e(g)\ra\,ds
\eeqn
for all $f,g\in L^2[0,1]$.
\bb\ni
We treat separately the individual constituents of the 
quantum stochastic integrals.
\bb\ni
(i) Let $z\in \BBC$, $g, h\in L^2[0,1]$ and let 
$M_t=M_t(E:\,d\GL)$. If $\Gv\in \FH^i$ then
\begin{eqnarray*}
\left|\ip{M_t{\,}^i_j e^j(g)}{\Gv }\right| &=&
\left|\int_0^t\ip{ g(s) E_s{\,}^{i-1}_{j-1} e^{j-
1}(g)}{\nabla_s\Gv }\,ds \right|
\\
&\leq & 
\int_0^t |g(s)|\nm{ E_s{\,}^{i-1}_{j-1}}\frac{\nm{g}^{j-1}}{((j-
1)!)^{\frac12}} \nm{\nabla_s \Gv }\,ds 
\\
&\leq &
\Lnu^{i-1}_{j-1}\frac{\nm{g}^{j-1}}{((j-1)!)^{\frac12}} \int_0^t 
|g(s)| \nm{\nabla_s \Gv }\,ds 
\\
&=&
\Lnu^{i-1}_{j-1}\frac{\nm{g}^{j-1}}{((j-1)!)^{\frac12}} i^{\frac12}\int_0^t 
|g(s)| \left(\int_{[0,1]^{i-1}}|\Gv (x;s)|^2\,dx\right)^{\frac12}\,ds
\\
&=& \Lnu^{i-1}_{j-1}\frac{\nm{g}^{j-1}}{((j-1)!)^{\frac12}} i^{\frac12}
\nm{g}\nm{\Gv }.
\end{eqnarray*}
Therefore $\nm{M_t{\,}^i_j e^j(g)}\leq 
\dyl i^{\frac12} \nm{g}\Lnu^{i-1}_{j-1}
\frac{\nm{g}^{j-1}}{((j-1)!)^{\frac12}}$
and,  since $\Gve(\nm{g})\in \FD(\Lnu)$, it follows that $\sum_{j=0}\iiy 
M_t{\,}^i_j e^j(g)$ is absolutely convergent to $\eta^i_t$ in $\FH^i$. If 
$z\in \BBC$ then  
\begin{eqnarray*}
\ip{\hat{ M}_t e( g)}{e(z h)} 
&=& \int_0^t\ip{  g(s) \tilde{ E}_s e(  g)}{z  h(s) e(z  h)}\,ds
\\
&=& \int_0^t\ip{  g(s) \sum_{j=0}\iiy E_s{\,}^i_j e^j(  g)}
{z  h(s) \sum_{i=0}\iiy e^i(z h)}\,ds
\\
&=& \sum_{i,j=0}\iiy\int_0^t\ip{  g(s)  E_s{\,}^i_j e^j(  g)}
{z  h(s) e^i(zh)}\,ds
\\
&=& \sum_{i,j=0}\iiy \ol{z} ^{i+1} \int_0^t\ip{g(s)  E_s{\,}^i_j e^j(g)}
{h(s) e^i(h)}\,ds
\\
&=& \sum_{i,j=1}\iiy\ol{z}^{i}
\int_0^t\ip{g(s)E_s{\,}^{i-1}_{j-1} e^{j-1}(g)}
{h(s) e^{i-1}(h)}\,ds
\\
&=& \sum_{i,j=1}\iiy \ol{z}^{i}\ip{M_t{\,}^i_j e^j(g)}{e^i(h)}
\\
&=&
\sum_{i=1}\iiy \ol{z}^{i}\ip{\eta^i_t}{e^i(h)}.
\end{eqnarray*}
The interchange of sum and integral in the third inequality follows from the 
Lebesgue series theorem since 
$$
\sum_{j=0}\iiy\nm{E_s{\,}^i_j e^j(g)}\leq \Lnu^i_j \Gve(\nm{g})<\infty 
$$
for almost all $s\in [0,1]$. Now
$$
\ip{\hat{M}_te(g)}{e( z h)}=\sum_{i=0}\iiy \ol{z}^i\ip{\left(\hat{M}_t 
e(g)\right)^i}{e^i(h)}
$$
for all $z\in \BBC$ so that 
$\ip{\eta ^i}{e^i(h)}=\ip{(\hat{ M}_t e(g))^i}{e^i(h)}$ whenever 
$h\in L^2[0,1]$. Therefore $(\hat{ M}_t e(g))^i=\eta^i_t$ and 
$$
\eta_t=(\eta^0_t,\eta^1_t,\ldots,\eta^i_t,\ldots)^t
$$
belongs to $\FH_{t]}$ and $e(g)\in \FD(M_t)$ with 
$M_t e(g)=\eta_t=\hat{M}_t e(g)$. Thus  $M_t$ is a 
cmx representation of $\hat{ M}_t$ for all $t$ in $[0,1]$.
\bb\ni
(ii) Let $z\in \BBC$, $g, h\in L^2[0,1]$ and let $M_t=M_t(F:\,dA)$. If 
$\Gv \in \FH^i$ then                                            
\begin{eqnarray*}
\left|\ip{M_t{\,}^i_j e^j(g)}{\Gv } \right|&=&
\left|\int_0^t \ip{F_s{\,}^i_{j-1} \nabla_s 
e^j(g)}{\Gv }\,ds\right|
\\
&=&
\left|\int_0^t \ip{g(s)F_s{\,}^i_{j-1} e^{j-1}(g)}{\Gv }\,ds\right|
\\
&\leq &
\int_0^t |g(s)|\nm{F_s{\,}^i_{j-1}} \nm{e^{j-1}(g)}\nm{\Gv }\,ds
\\
&\leq &
\nm{g}\nm{F^i_{j-1}}_2 \nm{e^{j-1}(g)}\nm{\Gv }
\\
&= &
\Lnu^i_{j-1} \nm{e^{j-1}(g)}\nm{g}\,\nm{\Gv },
\end{eqnarray*}
and $\nm{M_t{\,}^i_j e^j(g)}\leq \nm{g}\Lnu^i_j \nm{e^{j-1}(g)}$. Since 
$\Gve(\nm{g})$ belongs to $\CD(\Lnu)$ it follows that 
$\sum_{j=0}\iiy M_t{\,}^i_j e^j(g)$ is absolutely convergent to $\eta^i_t$ in 
$\FH^i$. If $z\in \BBC$ then 
\begin{eqnarray*}
\ip{ \hat{ M}_t e(g)}{e(z h)} &=& 
\int_0^t \ip{g(s) F_s e(g)}{e(z h)}\,ds
\\
&=&
\int_0^t \ip{g(s) \sum_{j=0}\iiy F_s{\,}^i_j  e^j(g)}{\sum_{i=0}\iiy z^i 
e^i(h)}\,ds
\\
&=&
\sum_{i,j=0}\iiy \ol{z}^i \int_0^t \ip{g(s)  F_s{\,}^i_j  e^j(g)}{
e^i(h)}\,ds
\\
&=&
\sum_{i,j=0}\iiy \ol{z}^i (j+1)^{\frac12}\int_0^t \ip{  F_s{\,}^i_j 
\nabla_s e^{j+1}(g)}{e^i(h)}\,ds
\\
&=&
\sum_{i,j=0}\iiy \ol{z}^i \ip{ M_t{\,}^i_{j+1} 
e^{j+1}(g)}{e^i(h)}
\\
&=& \sum_{i=0}\iiy \ol{z}^i\ip{\eta^i_t}{e^i(h)}.
\end{eqnarray*}
Since $\ip{\hat{ M}_t e(g)}{e(z h)}=\sum_{i=0}\iiy \ol{z}^i\ip
{(\hat{ M}_t e(g))^i}{e^i(h)} $ it follows as in (i) 
that $\eta^i_t=(\hat{ M}_t e(g))^i$ 
for all $i$ and that $M_t$ is a cmx representation of $\hat{ M}_t$ for 
all $t\in [0,1])$.
\bb\ni
(iii) Let $M_t=M_t(\,dA\d:G)$. 
Arguing as in (ii), if $\Gv \in \FH^i$ then
$$
|\ip{M_t{\,}^i_j e^j(g)}{\Gv }|\leq i^{\frac12} \nm{G^i_j}_2 
\nm{e^j(g)}\,\nm{\Gv }, 
$$
and $M_t$ is a cmx representation of $\hat{ M}_t$ for 
all $t\in [0,1]$.
\bb\ni
(iv) The case $M_t=M_t(H:\,ds)$ may be dealt with similarly. \,\,\epf
\bb\ni
Let $k$, $m$ and $n$ be non-negative integers and consider the 
cmx 
process $X_t=X_t(k,m,n)$ where $X_t(k,m,n)$ equals (the cmx representation 
of) $\GL^k_t A^m_t (A\d_t)^n$. A simple matrix calculation shows 
that $X_t{\,}^i_j=0$ unless $i=j+n-m$. 
\bb\ni
Let $\CP$ be the space of polynomials $P$ in $(\GL_t,A_t,A\d_t)$ with 
continuously differentiable coefficients. Each polynomial $P$ has a unique 
representation 
$$
P_t=\sum_{k,m,n=0}^K c_{k,m,n}(t)X_t(k,m,n).
$$
If $(E,F,G,H)$ is a quadruple in $\CP$ then  $(E,F,G,H)$
satisfies the conditions of Theorem \ref{T 74.001}.
\bb\ni
For simplicity we consider only $X_t(k,m,n)$ where $k=1$.
If $\psi^j\in \FH^j$ and $t\in [0,1]$ then
\begin{eqnarray*}
\nm{(A\d_t)^n\psi^j} &\leq & 
\left((j\!+\!1)\ldots(j\!+\!n)\right)^{\frac12}\nm{\psi^j}
\\
\nm{A_t^m(A\d_t)^n\psi^j} &\leq & 
\left((j\!+\!1)\ldots(j\!+\!n)(j\!+\!n\!-\!1)
\cdots(j\!+\!n\!-\!m)\right)^{\frac12}\nm{\psi^j}
\\
\nm{\GL_t A_t^m(A\d_t)^n\psi^j} &\leq & 
\left((j\!+\!1)\ldots(j\!+\!n)(j\!+\!n\!-\!1)
\cdots(j\!+\!n\!-\!m)\right)^{\frac12}(j\!+\!n\!-\!
m)^k\nm{\psi^j}.
\end{eqnarray*}
and 
$$
c_j=\nm{X^{j\!+\!n\!-\!m}_j}\leq 
\left((j\!+\!1)\ldots(j\!+\!n)(j\!+\!n\!-\!1)
\cdots(j\!+\!n\!-\!m)\right)^{1/2}(j\!+\!n\!-\!m)^k.
$$
It follows that $X$ is in $\CL^p_{\rm cmx}(\FH)\so$ with 
$\Lnu^i_j=0$ if  $i\neq j+n-m$ and $\Lnu^{j+n-m}_j\leq 2c_j$ 
for $p=1,2,\infty$. By the limit ratio test 
$$
\sum_{j=0}\iiy c_j^2 \frac{x^j}{j!}<\infty
$$
for all $x\in \BBR$  and the domain of $\Lnu(X)$ contains $\ell^2_{\rm exp}$.
\bb\ni
Since $\tilde{X}_t(k,m,n)$ equals ${\GL}_t^k A_t^m (A_t\d)^n$ on 
$\CE$ we have the following corollary.
\begin{Cor} 
Suppose the quadruple $(E,F,G,H)$ consists of polynomials in the 
basic cmx processes $\GL,A,A\d,t$. Then the quadruple
$(\tilde{E},\tilde{F},\tilde{G},\tilde{H})$ consists of the 
corresponding polynomials in the basic operator processes all defined 
on $\CE$. Let $M$ be the cmx quantum stochastic integral
$$
M_t=\int_0^t 
E_s\,d\GL_s+F_s\,dA_s+G_s\,dA\d_s+H_s\,ds.
$$
Then $\CE\subset \FD(M_t)$ for each $t\in [0,1]$ and, 
for each $g\in L^2[0,1]$,
$$
\tilde{M}_t e(g)
=\int_0^t{\,} (\tilde{ E}_s\,d{\GL}_s
+\tilde{ F}_s\,d{ A}_s+\tilde{ G}_s\,d{ A}\d_s+\tilde{ 
H}_s\,d{s}) e(g).
$$
\end{Cor}

\section{\label{IPF}The Ito Product Formula}
Our original proof of the Local Ito Formula was in the spirit of 
\cite{HP1} and used approximation by step processes. It was 
excessively long and we thank Mr A Belton for his very concise 
proof.
\bb\ni
He uses the Evans notation \cite{Ev1}.
$A_t{\,}^1_1=\GL_t$, $A_t{\,}^0_1=A_t$, $A_t{\,}^1_0=A\d_t$, 
$A_t{\,}^0_0=tI$.
If $(E,F,G,H)$ is an adapted integrable cmx quadruple put
\begin{eqnarray*}
M_t{\,}^1_1(E)=\int_0^t E_s\,dA_s{\,}^1_1,&\quad&
M_t{\,}^0_1(F)=\int_0^t F_s\,dA_s{\,}^0_1,
\\
M_t{\,}^1_0(G)=\int_0^t dA_s{\,}^1_0 G_s ,&\quad&
M_t{\,}^0_0(H)=\int_0^t H_s\,dA_s{\,}^0_0,
\end{eqnarray*}
We thank Mr A Belton for his assistance in proving the following lemma.
\begin{Lemma}[Local Ito Formula] \label{L 76.211} 
If $\Ga,\Gb,\Gg,\Gd\in \{0,1\}$ let 
$X\in \CL^p_{\rm cmx}(\FH)\so$ and $Y\in \CL^q_{\rm cmx}(\FH)\so$ 
be adapted processes, 
where $p=2(2-\Ga-\Gb)\inv$ and $q=2(2-\Gg-\Gd)\inv$. 
If $i,j,\mu,\nu\in \BBN$ and $\nu+\Gb=\mu+\Gg$ then
\begin{eqnarray}
\nonumber
M_t{\,}^{\Ga}_{\Gb}(X)^{i+\Ga}_{\nu+\Gb}
M_t{\,}^{\Gg}_{\Gd}(Y)^{\mu+\Gg}_{j+\Gd}
&=& \int_0^t M_v{\,}^{\Ga}_{\Gb}(X)^{i+\Ga-\Gg}_{\mu} 
Y_{v}{\,}^{\mu}_j\,dA_v{\,}^{\Gg}_{\Gd}
\\
\label{E 75.301}
&&
\hspace{-1.5in}+\int_0^t X_u{\,}^i_{\nu} M_u{\,}^{\Gg}_{\Gd}(Y)^{\nu}_{j+\Gd-
\Gb}\,dA_u{\,}^{\Ga}_{\Gb}
+\Gd^{\nu}_{\mu}\Gd^1_{\Gb}\Gd^{\Gg}_1\int_0^t(X_s{\,}^i_{\nu} 
Y_s{\,}^{\mu}_j)\,dA_s{\,}^{\Ga}_{\Gd}
\end{eqnarray}
\end{Lemma}
\pf In the following formulae we use two exponents of $\nabla$: 
$\nabla^1=\nabla$ and $\nabla^0=I$. If 
$\psi\in \FH^{j+\Gg}$ and $\Gv\in \FH^{i+\Ga}$ then
\begin{eqnarray}
\nonumber\ip{M_t{\,}^{\Gg}_{\Gd}(Y)^{\mu+\Gg}_{j+\Gd}\psi}
{M_t{\,}_{\Ga}^{\Gb}(X\st)_{i+\Ga}^{\nu+\Gb}\Gv}
&=& \\
\label{E 76.211}
&&\hspace{-1.8in}
\int_0^t\int_0^t\ip{\nabla^{\Gb}(u)Y^{\mu}_j\nabla^{\Gd}(v)\psi}
{\nabla^{\Gg}(v)(X_u\st)^{\nu}_i\nabla^{\Ga}(u) \Gv}\,dudv
\\
&& \label{E 76.212}
\hspace{-1.5in}+ \Gd^{\mu}_{\nu}\Gd^{\Gg}_1\Gd^1_{\Gb}\int_0^t
\ip{Y_s{\,}^{\mu}_j\nabla^{\Gd}(s)\psi}{(X_s\st)^{\nu}_i\nabla^{\Ga}(s)\Gv}\,ds.
\end{eqnarray}
This identity may be checked for each quadruple $\Ga,\Gb,\Gg,\Gd\in \{0,1\}$.
The entries in the left-hand side are given by 
the Definitions (\ref{E 75.101}), (\ref{E 75.102}), 
(\ref{E 75.103}) and (\ref{E 75.104}). The terms $\CS$ are then replaced by
terms in $\nabla$ using Proposition \ref{P 73.211} which shows that 
$\nabla=\CS\st$ and gives commutation relations for $\CS$ and $\nabla$.
\bb\ni
When $\nu=\mu$ the integral in (\ref{E 76.212}) is 
\beqn
\label{E 75.303}
\int_0^t\ip{X_s{\,}^i_{\nu} 
Y_s{\,}^{\nu}_j\nabla^{\Gd}(s)}{\nabla^{\Ga}(s)\Gv}\,ds=
\ip{\left(\int_0^t 
X_s{\,}^i_{\nu}Y_s{\,}^{\nu}_j\,dA_s{\,}^{\Ga}_{\Gd}\right)\psi}{\Gv}.
\eeqn
The double integral over $\{0\leq u,\,v\leq t\}$ in (\ref{E 76.212}) 
is the sum, $I_1+I_2$, of the integrals over 
$\{0\leq u\leq v\leq t\}$ and $\{0\leq v< u\leq t\}$. 
Since $X$ is adapted it follows from Theorem \ref{T 75.221} that 
$\nabla^{\Gg}(v)(X_u\st)^{\nu}_i\nabla^{\Ga}(u) \Gv=
(X_u\st)^{\nu-\Gg}_{i-\Gg}\nabla^{\Ga}(u)\nabla^{\Gg}(v) \Gv$
whenever $u\leq v$. Therefore 
\begin{eqnarray}
\nonumber
I_1 &=&
\int_0^t\int_0^v\ip{\nabla^{\Gb}(u)Y_v{\,}^{\mu}_j\nabla^{\Gd}(v)\psi}
{(X_u\st)^{\nu-\Gg}_{i-\Gg}\nabla^{\Ga}(u)\nabla^{\Gg}(v) \Gv}\,dudv
\\
\nonumber &=&
\int_0^t\int_0^v
\ip{X_u{\,}_{\nu-\Gg}^{i-\Gg}\nabla^{\Gb}(u)Y_v{\,}^{\mu}_j\nabla^{\Gd}(v)\psi}
{\nabla^{\Ga}(u)\nabla^{\Gg}(v) \Gv}\,dudv
\\
\nonumber &=&
\int_0^t\ip{ M_v{\,}^{\Ga}_{\Gb}(X)^{i+\Ga-\Gg}_{\nu}Y_v{\,}^{\nu}_j 
\nabla^{\Gd}(v)\psi}{\nabla^{\Gg}(v) \Gv}\,dv
\\
&=&
\label{E 75.304}
\ip{\left(\int_0^t M_v{\,}^{\Ga}_{\Gb}(X)^{i+\Ga-\Gg}_{\nu}Y_v{\,}^{\nu}_j 
\,dA_v{\,}^{\Gg}_{\Gd}\right)\psi}{\Gv}
\end{eqnarray}
Similarly 
\beqn
\label{E 75.305}
I_2= \ip{\left(\int_0^t X_u{\,}^i_{\nu} M_u{\,}^{\Gg}_{\Gd}(Y)^{\nu}_{j+\Gd-
\Gb}\,dA_u{\,}^{\Ga}_{\Gb}\right)\psi}{\Gv}.
\eeqn
The identity (\ref{E 75.301}) follows immediately from 
(\ref{E 75.303}), (\ref{E 75.304}) and (\ref{E 75.305}).\,\,\epf 
\begin{Cor}[Quantum Ito Product Formula] \label{C 75.801} If 
the product processes 
$$
M^{\Ga}_{\Gb}(X)Y,\qquad X M^{\Gg}_{\Gd}(Y),\qquad \Gd^1_{\Gb}\Gd^{\Gg}_1 XY
$$
exist in $\CL^q_{\rm cmx}(\FH)\so$, $\CL^p_{\rm cmx}(\FH)\so$ and 
$\CL^r_{\rm cmx}(\FH)\so$, where 
$1/p+1/q=1/r$ then  the product process
$M^{\Ga}_{\Gb}(X)M^{\Gg}_{\Gd}(Y)$ exists in $\CL\iiy_{\rm cmx}(\FH)\so$ and
\begin{eqnarray}\nonumber
M_t{\,}^{\Ga}_{\Gb}(X)M_t{\,}^{\Gg}_{\Gd}(Y)&=&
\int_0^t M_s{\,}^{\Ga}_{\Gb}(X) Y_s\,dA_s{\,}^{\Gg}_{\Gd}+
\int_0^t X_sM_s{\,}^{\Gg}_{\Gd}(Y)\,dA_s{\,}^{\Ga}_{\Gb}\\
\label{E 75.800}
&&+\Gd^1_{\Gb}\Gd^{\Gg}_1 \int_0^t X_s Y_s \,dA_s{\,}^{\Ga}_{\Gd}.
\end{eqnarray}
\end{Cor}
\pf This follows directly from the identity (\ref{E 75.301}) together with  
the norm bounds (\ref{E 75.204}), (\ref{E 75.205}) and (\ref{E 
75.206}) in Theorem \ref{T 74.801}.\,\,\epf
\bb\ni
We give sufficient conditions that the conditions of Corollary \ref{C 75.801} 
are satisfied.
\begin{Cor}\label{C 75.810} Suppose $\Ga,\Gb,\Gg,\Gd, p,q,X,Y$ 
satisfy the conditions of Lemma
{\rm \ref{L 76.211}} and define scalar matrices $\LGk(X), 
\LGk(Y)$ by the formulae
$$
\LGk(X)^i_j= i^{\frac{\Ga}{2}}\Lnu(X)^{i-\Ga}_{j-\Gb} j^{\frac{\Gb}{2}},\quad
\LGk(Y)^i_j= i^{\frac{\Gg}{2}}\Lnu(X)^{i-\Gg}_{j-\Gd} j^{\frac{\Gd}{2}}.
$$
If the scalar matrix 
$\LGk(X)\Lnu(Y)+\Lnu(X)\LGk(Y)+\Lnu(X)\Lnu(Y)$ has finite entries then $X$ 
and $Y$ satisfy the conditions of Corollary {\rm \ref{C 75.801}} and the 
quantum Ito  product formula {\rm (\ref{E 75.800})} is valid.
\end{Cor}
\pf  It follows from the norm bounds (\ref{E 75.204}), (\ref{E 75.205}) and  
(\ref{E 75.206}), that the r.h.s of (\ref{E 
75.301}) is bounded by
\begin{eqnarray*}
&& (i+\Ga)^{\frac{\Gg}{2}}(j+\Gd)^{\frac{\Gd}{2}}(i+\Ga-\Gg)^{\frac{\Ga}{2}}
\mu^{\frac{\Gb}{2}} \nm{X^{i-\Gg}_{\mu-\Gb}}_p\nm{Y^{\mu}_{\Gg}}_q 
\\
&& \hspace{.25in}+
(i+\Ga)^{\frac{\Ga}{2}}(j+\Gd)^{\frac{\Gb}{2}}(j+\Gd-\Gb)^{\frac{\Gd}{2}}
\nu^{\frac{\Gg}{2}} \nm{X^{i}_{\nu}}_p\nm{Y^{\nu-\Gg}_{j-\Gb}}_q 
\\
&&\hspace{.5in}+ \Gd^1_{\Gb}\Gd^{\Gg}_1(i+\Ga)^{\frac{\Ga}{2}}(j+\Gd)^{\frac{\Gd}{2}}
\nm{X^i_{\nu}}p\nm{Y^{\nu}_j}_q
\\
&& \leq 
(j+\Gd)^{\frac{\Gd}{2}}(i+\Ga-\Gg)^{\frac{\Ga}{2}}
\LGk(X)^i_{\mu}\Lnu(Y)^{\mu}_j
+
(i+\Ga)^{\frac{\Ga}{2}}(j+\Gd-\Gb)^{\frac{\Gd}{2}}
\Lnu(X)^i_{\nu}\LGk(Y)^{\nu}_j 
\\ 
&& \hspace{0.5in}
+ \Gd^1_{\Gb}\Gd^{\Gg}_1(i+\Ga)^{\frac{\Ga}{2}}(j+\Gd)^{\frac{\Gd}{2}}
\Lnu(X)^i_{\nu}\Lnu(Y)^{\nu}_j.
\end{eqnarray*}
Therefore
\begin{eqnarray*}
\sum_{\nu=0}\iiy\nm{M_t{\,}^{\Ga}_{\Gb}(X)^{i+\Ga}_{\nu+\Gb}
M_t{\,}^{\Gg}_{\Gd}(Y)^{\mu+\Gg}_{j+\Gd}}
&\leq &
(j+\Gd)^{\frac{\Gd}{2}}(i+\Ga-\Gg)^{\frac{\Ga}{2}}
(\LGk(X)\Lnu(Y))^i_j
\\
&&
\hspace{-180pt}+
(i+\Ga)^{\frac{\Ga}{2}}(j+\Gd-\Gb)^{\frac{\Gd}{2}}
(\Lnu(X)\LGk(Y))^{i}_j 
+
\Gd^1_{\Gb}\Gd^{\Gg}_1(i+\Ga)^{\frac{\Ga}{2}}(j+\Gd)^{\frac{\Gd}{2}}
(\Lnu(X)\Lnu(Y))^i_j.
\end{eqnarray*}
which is finite for all $i,j$. It follows that the 
the product process $M_t{\,}^{\Ga}_{\Gb}(X)M_t{\,}^{\Gg}_{\Gd}(Y)$ exists. 
\bb\ni
Now
\begin{eqnarray*}
\nm{M_v{\,}^{\Ga}_{\Gb}(X)^{i+\Ga-\Gg}_{\mu} Y_{v}{\,}^{\mu}_j }
&\leq &
(i+\Ga-\Gg)^{\frac{\Ga}{2}}\mu^{\frac{\Gb}{2}}\nm{X^{i-\Gg}_{\mu-
\Gb}}_p\nm{Y_v{\,}^{\mu}_j}
\\
&\leq &
(i+\Ga-\Gg)^{\frac{\Ga}{2}}\LGk(X)^i_{\mu}\nm{Y_v{\,}^{\mu}_j},
\end{eqnarray*}
and
$$
\sum_{\mu=0}\iiy 
\nm{M^{\Ga}_{\Gb}(X)^{i+\Ga-\Gg}_{\mu} Y^{\mu}_j }_q\leq
(i+\Ga-
\Gg)^{\frac{\Ga}{2}}\sum_{\mu=0}\iiy\LGk(X)^i_{\mu}\Lnu(Y)^{\mu}_j<\infty.
$$
Therefore the product chaos matrix $M_s{\,}^{\Ga}_{\Gb}(X) Y_s$ exists for 
almost all $s\in [0,1]$ and defines a cmx process in $\CL^p_{\rm cmx}(\FH)\so$.
Similar arguments show that the product chaos matrix
$X_sM_s{\,}^{\Gg}_{\Gd}(Y)$ exists almost everywhere and defines a process in 
$\CL^p_{cmx}(\FH)\so$ and that, for $\Gb=\Gd=1$, $s\mapsto X_s Y_s$ defines a process 
in the appropriate $\CL^r_{\rm cmx}(\FH)\so$. 
\bb\ni
Taking weak sums on both sides and using the DCT yields 
(\ref{E 75.800}).\,\,\epf
\bb\ni
In particular if $X$ and $Y$ are $
(2k+1)$-diagonal matrices for some $k\in \BBN$
then the conditions of Corollary \ref{C 75.810} are satisfied. Since the 
basic processes are tri-diagonal any polynomial in them is $(2k+1)$-diagonal for 
some $k\in \BBN$.
\begin{Cor} If $X$ and $Y$ are polynomials in the basic processes then the 
quantum Ito product formula {\rm (\ref{E 75.800})} is valid.
\end{Cor}
\section{\label{QDF} Quantum Duhamel Formula} In this section we show that 
the quantum Duhamel formula \cite[Proposition 5.2]{Vin2} for 
regular quantum semimartingales remains true for some quantum 
stochastic integrals and some essentially self-adjoint
quantum semimartingales. The proofs, which use series expansions 
and Nelson's analytic vector theorem, 
follow the same lines as the proof for regular quantum 
semimartingales.
\bb\ni
\begin{Thm}\label{T 79.797} Let $(E,F,F\st,H)$ be an adapted 
symmetric integrable cmx quadruple and let $\LGk$ be the control 
matrix of the process
\beqn\label{E 79.770}
M_t=\int_0^t (E \,d\GL+ F \,dA+F\st\,dA\d+H\,ds).
\eeqn
If $\Ge>0$ and $\ell_{00}\subset \CA_{\Ge}(\LGk)$ then, for  
$-\Ge <p<\Ge$,
\bb\ni 
{\rm (i)} the product cmx processes 
\beqn\label{E 75.307}
M^{\Ga} F^{\mu}(M+E)^{\Gb} G^{\nu} M^{\Gg} \quad
\mbox{ exists in } \CL^p_{\rm cmx}(\FH)\so,\quad p=\frac{2}{\mu+\nu}
\eeqn
whenever $\Ga,\Gb,\Gg\in \BBN$ and $\mu,\nu\in \{0,1\}$.
\bb\ni 
{\rm (ii)} for $n=1,2,\ldots$ the 
formulae 
\begin{eqnarray*} E_n&=& (M+E)^n-M^n \\ F_n &=& 
\sum_{\Ga+\Gb=n-1} M^{\Ga}F(M+E)^{\Gb} \\ (F_n)\st &=& 
\sum_{\Ga+\Gb=n-1} (M+E)^{\Ga}F\st M^{\Gb} \\ H_n &=& 
\sum_{\Ga+\Gb=n-1} M^{\Ga}H M^{\Gb}+ \sum_{\Ga+\Gb+\Gg=n-2} 
M^{\Ga}F(M+E)^{\Gb}F\st M^{\Gg} 
\end{eqnarray*} 
define an adapted symmetric 
integrable cmx quadruple $(E_n, F_n, (F_n)\st, H_n)$.
\bb\ni
{\rm (iii)} $M^n=\{(M_t)^n:t\in [0,1]\}$ is a symmetric cmx 
quantum stochastic integral with
\beqn\label{E 79.771} 
\qquad(M_t)^n=\int_0^t (E_n 
\,d\GL+F_n\,dA+(F_n)\st\,dA\d+H_n\,ds),\quad n=1,2,\ldots.
\eeqn         
{\rm (iv)} $\FH_{00}\subset \CA_{\Ge}(M_t)\cap 
\CA_{\Ge}(\tilde{M}_t)$ and $\FH_{00}\subset \CA_{\Ge}((M+E)_t)\cap 
\CA_{\Ge}(\widetilde{(M+E)}_t)$ for $t\in [0,1]$.
\bb\ni 
{\rm (v)} $\tilde{M}$ and $\widetilde{M\!+\!E}$ are 
essentially self-adjoint with core $\FH_{00}$.
\bb\ni
{\rm (vi)} The series $I+\sum_{k=1}\iiy(ip M)^k/k!$ is convergent in 
$\CL\iiy_{\rm cmx}(\FH)\so$ to a process $J(p)$ and 
$$
\tilde{J}_t(p)=e^{ip\ol{M}_t}
$$
{\rm (vii)} The series $I+\sum_{k=1}\iiy(ip (M+E))^k/k!$ is convergent in 
$\CL\iiy_{\rm cmx}(\FH)\so$ to a process $K(p)$ and 
$$
\tilde{K}_t(p)=e^{ip\ol{M+E}_t}
$$
\end{Thm} 
\pf 
The proofs for $M+E$ are similar to the proofs for $M$ and will 
be omitted. 
\bb\ni
Since $\CA(\LGk)\neq \emptyset $ the  
product matrix $\LGk^k$ exists for all $k\in\BBN$ by definition. 
If $\Lnu=\Lnu(M)$ is the scalar 
matrix of $M$ then $\LGk+\Lnu\prec 2\LGk$ and $(\LGk+\Lnu)^k$ 
exists for all $k\in \BBN$. Since $\LGk$ and $\Lnu$ have 
non-negative entries the 
product scalar matrix $\LGk^{\Ga_1}\Lnu^{\Gb_1}\cdots 
\LGk^{\Ga_k}\Lnu^{\Gb_k}$ exists whenever 
$\Ga_1,\ldots,\Ga_k,\Gb_1,\ldots,\Gb_k$ are non-negative integers 
with $\Ga_1+\cdots+\Ga_k+\Gb_1+\cdots+\Gb_k=n$. Moreover $(\LGk+\Lnu)^n$ is a 
finite sum of such monomial products. 
\bb\ni 
(i) If $1,\Ga,\Gb,\Gg\leq k$ then $\Lnu(M^{\Ga}), \Lnu(F), 
\Lnu((M+E)^{\Gb}), \Lnu(G), \Lnu(M^{\Gg}) \prec (\LGk+\Lnu)^k$. 
Since $M,E\in \CL\iiy_{\rm cmx}(\FH)\so$, $F,G\in \CL^2_{\rm cmx}(\FH)\so$, 
and $H\in \CL^1_{\rm cmx}(\FH)\so$ and the product scalar matrix 
$(\LGk+\Lnu)^{5k}$ exists for all $k$ it follows from Lemma \ref{L 78.119} 
that $M^{\Ga}F^{\mu}(M+E)^{\Gb}G^{\nu} M^{\Gg}$ 
exists whenever $\Ga,\Gb,\Gg\in \BBN$ and $\mu,\nu\in \{0,1\}$. 
\bb\ni
(ii) It follows from (i) and Corollary \ref{C 75.801} 
that the formulae in (ii)  define 
cmx processes $E_n$, $F_n$, $(F_n)\st$ and $ H_n$. 
The following algebraic identities follow directly from those formulae.
\begin{eqnarray*}
E_{n+1}&=&M(M+E)^n-M^{n+1}+EM^n+E(M+E)^n-EM^n
\\
&=& ME_n+EM_n+EE_n;
\\
F_{n+1}&=& \!\!\!\!\sum_{\Ga+\Gb=n-1}\!\!\!\!M^{\Ga+1}F(M+E)^{\Gb}+FM^n+F(M+E)^n-FM^n
\\
&=& MF_n+FM^n+FE_n;
\\
(F_{n+1})\st &=&\!\!\!\! \sum_{\Ga+\Gb=n-1}\!\!\!\!M(M+E)^{\Ga}M^{\Gb}+F\st M^n
+ \!\!\!\!\sum_{\Ga+\Gb=n-1}\!\!\!\!E(M+E)^{\Ga}F\st M^{\Gb}
\\
&=& M(F_n)\st+F\st M^n+E(F_n)\st;
\\
H_{n+1}&=& HM^n+\!\!\!\!\!\sum_{\Ga+\Gb=n-1}\!\!\!\!\!M^{\Ga+1}HM^{\Gb}
+\!\!\!\!\!\sum_{\Ga+\Gb+\Gg=n-2}\!\!\!\!\!M^{\Ga+1} F(M+E)^{\Gb} F\st M^{\Gg}
\\
&& + \sum_{\Gb+\Gg=n-1}\!\!\!\!\!F(M+E)^{\Gb} F\st M^{\Gg}
\\
&=& MH_n+HM^n+F(F_n)\st.
\end{eqnarray*}
The quadruple $(E_n,F_n,G_n,H_n)$ is integrable by (i) and adapted 
by Corollary \ref{C 74.01}. It follows from (i) and 
Lemma \ref{L 73.00} than $G_n=(F_n)\st$ for all $n\in \BBN$.
\bb\ni
(iii) Applying the Quantum Ito product formula (\ref{E 75.800}) to $M_t 
M^n_t$, where $M$ and $M^n$ are given by (\ref{E 79.770})and 
(\ref{E 79.771}), gives
\begin{eqnarray*}
M_t^{n+1} &=& \int_0^t (M\,dM^n+dM M^n+dM\,dM^n)
\\
&=& \int_0^t \left((M E_n+EM^n+EE_n)\,d\GL+(MF_n+FM^n+FE_n)\,dA
\phantom{A\d}\right.
\\
&& \left.
\hspace{-0.5in}+(M(F_n)\st+F\st M^n+E(F_n)\st)dA\d+(MH_n+HM^n+F(F_n)\st)ds\right)
\\
&=& \int_0^t (E_{n+1}d\GL+F_{n+1}dA+ F\st_{n+1}dA\d+H_{n+1}ds).
\end{eqnarray*}
(iv), (v), (vi) Since $\Lnu(M)\prec \LGk$ and $\ell_{00}\subset 
\CA_{\Ge}(\LGk)$ implies that $\ell_{00}\subset 
\CA_{\Ge}(\Lnu(M))$. Parts (iv), (v) and (vi) 
then follow from Theorem \ref{T 76.02}. 
\,\,\epf
\bb\ni
Under the above conditions we define $e^{ipM_t}$ and 
$e^{ip(M+E)_t}$ to be the chaos matrices $ J_t(p)$ and $K_t(p)$ 
respectively. We now have the following chaos matrix version of 
the quantum Duhamel formula.
\begin{Cor}[Cmx Duhamel Formula]\label{C 79.01} The chaos matrix process 
$\{e^{ipM_t}:t\in [0,1]\}$ has a representation as a cmx quantum 
stochastic integral
\beqn\label{EE 79.01}
e^{ipM_t}= I+\int_0^t (E_{{\rm exp}(ipM)}\,d\GL+F_{{\rm exp}(ipM)}\,dA
+G_{{\rm exp}(ipM)}\,dA\d+H_{{\rm exp}(ipM)}\,ds).
\eeqn
The integrands form an adapted integrable quadruple and are 
defined by the formulae
\begin{eqnarray}
\nonumber E_{{\rm exp}(ipM)}(t)&=& ip\int_0^1 e^{i(1-u)pM_t} E_t e^{iupM_t}\,du 
=e^{ip(M_t+E_t)}-e^{ipM_t}
\\
\nonumber F_{{\rm exp}(ipM)}(t)&=& ip\int_0^1 e^{i(1-u)pM_t} F_t e^{iup(M_t+E_t)}\,du
\\
\label{E 79.977}G_{{\rm exp}(ipM)}(t)&=& ip\int_0^1 e^{i(1-u)p(M_t+E_t)} F\st_t e^{iu pM_t}\,du
\\
\nonumber H_{{\rm exp}(ipM)}(t)&=& ip\int_0^1 e^{i(1-u)pM_t} H_t e^{iu pM_t}\,du
\\
\nonumber &&\hspace{0.5in} -p^2\int_0^1\int_0^1
ue^{i(1-u)pM_t}F_te^{iu(1-v)p(M_t+E_t)} F\st_t e^{iuv 
pM_t}\,du\,dv.
\end{eqnarray} 
They have power series representations
$$
E_{{\rm exp}(ipM)}  =  \sum_{n=0}\iiy \frac{(ip)^n E_n}{n!},
 \quad
F_{{\rm exp}(ipM)}  =  \sum_{n=0}\iiy \frac{(ip)^n F_n}{n!},
 \quad
G_{{\rm exp}(ipM)}  =  \sum_{n=0}\iiy \frac{(ip)^n F\st_n}{n!},
$$
$$
H_{{\rm exp}(ipM)}  =  \sum_{n=0}\iiy \frac{(ip)^n H_n}{n!},
$$
the series converging in $\CL\iiy_{\rm cmx}(\FH)\so$, 
$\CL^2_{\rm cmx}(\FH)\so$, 
$\CL^2_{\rm cmx}(\FH)\so$ and $\CL^1_{\rm cmx}(\FH)\so$ respectively.
\end{Cor}
\pf There is no loss of generality in taking $\Ge>1$ and $p=1$.
Both $M$ and $M+E$ satisfy the conditions of 
Theorem \ref{T 76.02}.
\bb\ni
If $0\leq u\leq 1$ define $S_k(u)$ in $\CL\iiy_{\rm cmx}(\FH)\so$ and 
$T_k(u)$ in $\CL^2_{\rm cmx}(\FH)\so$ by
$$
S_k(u)=\sum_{\Ga=0}^k\frac{(i(1-u)M)^{\Ga}}{\Ga!},\qquad
T_k(u)=\sum_{\Gb=0}^k (iF)\frac{(iu(M+E))^{\Gb}}{\Gb!}.
$$
By Theorem \ref{T 79.797}, 
the product processes $M^\Ga F(M+E)^\Gb$ exist in 
$\CL^2_{\rm cmx}(\FH)\so$ for all $\Ga,\Gb\in \BBN$. Therefore product process $S_k(u) T_k(u)$ 
exists for all $u\in [0,1]$ and 
$$
S_k(u) T_k(u)=\sum_{\Ga,\Gb=0}^k 
\frac{(i(1-u)M)^{\Ga}}{\Ga!} (iF)\frac{(iu(M+E))^{\Gb}}{\Gb!}
$$
Now $u\mapsto S_k(u)$ and $u\mapsto T_k(u)$ are continuous functions into 
$\CL\iiy_{\rm cmx}(\FH)\so$ and  $\CL^2_{\rm cmx}(\FH)\so$ respectively and:
\bb\ni
(a) $S_k(u)$ converges to $e^{i(1-u)M}$ in $\CL\iiy_{\rm cmx}(\FH)\so$ and 
$\Lnu(S_k)\prec e^{(1-u)\LGk}$ $k=1,2,\ldots$;
\bb\ni
(b) $T_k(u)$ converges to $(iF)e^{iuM}$ in $\CL^2_{\rm cmx}(\FH)\so$ and $
\Lnu(T_k)\prec \LGk e^{u\LGk}$ $k=1,2,\ldots$.
\bb\ni
Therefore $\Lnu(S_k(u))\Lnu(T_k(u)) \prec \LGk e^{\LGk}$ for all $u\in [0,1]$ so 
by Proposition  \ref{P 78.311} it follows that 
$S_k(u)T_k(u)$ converges in $\CL^2_{\rm cmx}(\FH)\so$ 
to $e^{i(1-u)M}(iF) e^{iu(M+E)}$  for 
each $u\in [0,1]$ and 
\begin{eqnarray*}
\lim_{k\roy}\int_0^1 S_k(u)T_k(u)\,du
&=&\lim_{k\roy} \sum_{\Ga,\Gb=0}^k \frac{(iM)^{\Ga} (iF)(i(M+E))^{\Gb}}
{(\Ga+\Gb+1)!}
\\
&=&\int_0^1 e^{i(1-u)M}(iF) e^{iu(M+E)}\,du=F_{{\rm exp}(iM)}.
\end{eqnarray*}
If $\CR_k$ in $\CL^2_{\rm cmx}(\FH)\so$ is defined by
$$
\CR_k\defeq \left(\int_0^1 S_k(u)T_k(u)\,du\right)-\left(\sum_{n=0}^k \frac{i^n 
F_n}{n!}\right)
=\!\!\!\!\!\!\sum_{0\leq \Ga,\Gb\leq k<\Ga+\Gb}
\!\!\!\!\!\!\frac{(iM)^{\Ga} (iF)(i(M+E))^{\Gb}}
{(\Ga+\Gb+1)!}
$$
then $\lim_{k\roy}\CR_k=0$ in $\CL^2_{\rm cmx}(\FH)\so$. 
This follows from Lemma \ref{L 78.019} since
$$
\Lnu(\CR_k)\prec 
\sum_{\Ga+\Gb>k}\frac{\LGk^{\Ga+\Gb+1}}{(\Ga+\Gb+1)!}\prec
\sum_{n=k}\iiy \frac{\LGk^n}{n!}
$$
which tends to the zero matrix as $k\roy$. Therefore
$$
F_{{\rm exp}(iM)}=\sum_{n=0}\iiy \frac{i^n F_n}{n!},
$$
the series converging in $\CL^2_{\rm cmx}(\FH)\so$.
\bb\ni
The proof of the power series representations 
are similar and will be omitted.
\bb\ni
It follows from (\ref{E 79.771}) that
\begin{eqnarray*}
e^{ipM_t}&=&I+\sum_{n=1}\iiy(ip)^n\int_0^t(E_n 
\,d\GL+F_n\,dA+(F_n)\st\,dA\d+H_n\,ds) \\
&=& I+\int_0^t (E_{{\rm exp}(ipM)}\,d\GL+F_{{\rm exp}(ipM)}\,dA
+G_{{\rm exp}(ipM)}\,dA\d+H_{{\rm exp}(ipM)}\,ds).
\end{eqnarray*}
The interchange of the integrals with the sums is justified by 
Theorem \ref{T 74.802}.
\,\,\epf
\bb\ni
The next problem is to deduce from (\ref{EE 79.01}) 
a quantum Duhamel formula for quantum stochastic integrals. When the 
integrands $(\hat{E},\hat{F},\hat{F}\st,\hat{H})$ are not 
required to be bounded we are unable to find reasonable 
sufficient conditions to ensure that integrands can be defined on 
$\CE$ by the formulae  (\ref{E 79.999}). We therefore restrict 
attention to the case that $M$ is a cmx semimartingale. 
\bb\ni
If the quantum Duhamel formula is valid for all real $p$ then the 
functional quantum Ito formula may be proved as in \cite[Theorem 
6.2.]{Vin2} using the Fourier functional calculus 
$$
f(M_t)= \int_{\BBR} \hat{f}(p)e^{ipM_t}\,dp.
$$
Corollary \ref{C 79.01} only implies that $pM$ satisfies the 
quantum Duhamel formula for $|p|<\Ge$. The following corollary 
shows this is sufficient. Repeated applications of the quantum Ito product 
formula of Hudson and Parthasarathy are used to show that $pM$ 
satisfies the quantum Duhamel formula for all real $p$. 
This is true in particular if $M$ is one of the 
quantum semimartingales described in Proposition \ref{PP 76.01}.
\bb\ni
\begin{Cor}\label{C 80.221} 
Let $\hat{M}=\{\hat{M}_t:t\in [0,1]\}$ be a symmetric quantum 
semimartingale 
$$
\hat{M}_t=\int_0^t (\hat{E}_s \,d\GL_s+ \hat{F}_s \,dA_s
+\hat{F}\st\,dA\d_s+\hat{H}_s\,ds)
$$
with corresponding symmetric cmx semimartingale $M=\{M_t:t\in 
[0,1]\}$: 
$$
M_t=\int_0^t (E_s \,d\GL_s+ F_s \,dA_s+F\st\,dA\d_s+H_s\,ds).
$$
If there exists $\Ge>0$ such that $\ell_{00}\subset 
\CA_{\Ge}(\LGk)$ where $\LGk$ is the control matrix of $M$ then
\bb\ni
{\rm (i)}  $\hat{M}$ is essentially self-adjoint with core $\FH_{00}$.
\bb\ni
{\rm (ii)} If $\ol{M}_t$ is the closure of $\hat{M}_t$ then 
$e^{ipM}$ is a regular quantum semimartingale and $pM$ satisfies  
the quantum Duhamel formula 
\beqn\label{E 80.140} 
e^{ip\ol{M}_t}= I\!+\!\int_0^t (\hat{E}_{{\rm exp}(ipM)}\,d\GL
\!+\!\hat{F}_{{\rm exp}(ipM)}\,dA
\!+\!\hat{G}_{{\rm exp}(ipM)}\,dA\d
\!+\!\hat{H}_{{\rm exp}(ipM)}\,ds),
\eeqn
for all real $p$, where 
\begin{eqnarray}
\nonumber\hat{E}_{{\rm exp}(ipM)}(t)&
=& e^{ip(\ol{M}_t+\hat{E}_t)}-e^{ip\ol{M}_t}
\\
\nonumber\hat{F}_{{\rm exp}(ipM)}(t)&=& ip\int_0^1 e^{i(1-u)p\ol{M}_t} 
\hat{F}_t e^{iup(\ol{M}_t+\hat{E}_t)}\,du
\\
\label{E 79.999}\hat{G}_{{\rm exp}(ipM)}(t)&=& ip\int_0^1 e^{i(1-
u)p(\ol{M}_t+\hat{E}_t)} \hat{F}\st_t e^{iu p\ol{M}_t}\,du
\\
\nonumber\hat{H}_{{\rm exp}(ipM)}(t)&=& ip\int_0^1 e^{i(1-u)p\ol{M}_t} 
\hat{H}_t e^{iu p\ol{M}_t}\,du
\\
\nonumber &&+(ip)^2\int_0^1\int_0^1
ue^{i(1-u)p\ol{M}_t}\hat{F}_te^{iu(1-v)p(\ol{M}_t+\hat{E}_t)} 
\hat{F}\st_t e^{iuv p\ol{M}_t}\,du\,dv
\end{eqnarray} 
\end{Cor}
\pf There is no loss of generality in assuming $\Ge=1$; simply 
replace $M$ by $\Ge M$. 
\bb\ni
It follows from \cite[Theorems VIII.25, VIII.21]{RS1} that 
$u\mapsto e^{i(1-u)p \ol{M}_t}$ and $u\mapsto 
e^{iup(\ol{M}_t+\tilde{E}_t)}$ 
are strongly continuous. Therefore the function
$u\mapsto e^{i(1-u)p M_t}F_t e^{iup(M_t+E_t)}$ belongs to 
$L\cmx\iiy(\CH)\so$ for almost all $t\in [0,1]$.
Theorems \ref{T 73.002} and \ref{T 79.797} (vi), (vii) now 
imply that $F_{{\rm exp}(ipM)}(t)$ is a cmx representation of 
$\hat{F}_{{\rm exp}(ipM)}(t)$ for almost all $t$.
The other integrands in 
(\ref{E 79.977}) and (\ref{E 79.999}) are similarly related.
By Theorem \ref{T 74.444} Equation (\ref{E 80.140}) 
is true for $-1\leq p\leq 1$.
\bb\ni
Since $e^{i\ol{M}}$ is a unitary process it follows that 
$\hat{E}_{{\rm exp}(iM)}$ is in $L\iiy_{\rm cmx}(\FH)\so$ while 
$\hat{F}_{{\rm exp}(iM)}$, $\hat{G}_{{\rm exp}(iM)}$ belong to  
$L^2_{\rm cmx}(\FH)\so$ 
and $\hat{H}_{{\rm exp}(iM)}$ is in $L^1_{\rm cmx}(\FH)\so$ so that 
$e^{i\ol{M}}$ is a regular quantum semimartingale.
\bb\ni
We now prove (\ref{E 80.140}) is valid for $p=2$. 
\bb\ni
By the quantum Ito product formula of Hudson and 
Parthasarathy, \cite[Theorem 5.]{Att}, \cite[Theorem 4.5]{HP1},  
\beqn\label{E 80.320}
e^{2i\ol{M}_t}&=&I+\int_0^t e^{i\ol{M_s}}d\left(e^{i\ol{M_s}}\right) +  
d\left(e^{i\ol{M_s}}\right)e^{i\ol{M_s}}+
d\left(e^{i\ol{M_s}}\right)d\left(e^{i\ol{M_s}}\right).
\eeqn
The coefficient of $d\GL$ is $e^{i\ol{M}}\hat{E}_{{\rm exp}(i\ol{M})}
+\hat{E}_{{\rm exp}(i\ol{M})}e^{i\ol{M}}
+\hat{E}_{{\rm exp}(i\ol{M})}\hat{E}_{{\rm exp}(i\ol{M})}$ which equals
$$
e^{i\ol{M}}e^{i(\ol{M}+\hat{E})}-e^{i2\ol{M}}+e^{i(\ol{M}
+\hat{E})}e^{i\ol{M}}
-e^{i2\ol{M}}+e^{i(2\ol{M}+2\hat{E})}+e^{i2\ol{M}}
-e^{i\ol{M}}e^{i(\ol{M}+\hat{E})}
-e^{i(\ol{M}+\hat{E})} e^{i\ol{M}}
$$
$$
=\hat{E}_{{\rm exp}(i2\ol{M})}
$$
The coefficient of $dA$ is
\begin{eqnarray*}
e^{i\ol{M}}\hat{F}_{{\rm exp}(i\ol{M})} 
+  \hat{F}_{{\rm exp}(i\ol{M})}e^{i\ol{M}}
+\hat{F}_{{\rm exp}(i\ol{M})}\hat{E}_{{\rm exp}(i\ol{M})} 
\hspace{-186pt}&&
\\
&=&  i \int_0^1 e^{i(2-u)\ol{M}} \hat{F} e^{iu(\ol{M}+\hat{E})}\,du +
 i \int_0^1 e^{i(1-u)\ol{M}} \hat{F} e^{iu(\ol{M}+\hat{E})}e^{i\ol{M}}\,du
\\
&& \hspace{0.5in}+ i \int_0^1 e^{i(1-u)\ol{M}} \hat{F} 
e^{iu(\ol{M}+\hat{E})}\left(e^{iu(\ol{M}+\hat{E})}-e^{i\ol{M}}\right) \,du
\\
&=&i\int_0^1\left(e^{i(2-u)\ol{M}} \hat{F} e^{iu(\ol{M}+\hat{E})}
+e^{i(1-u)\ol{M}} \hat{F} e^{i(1+u)(\ol{M}+\hat{E})}\right)\,du
\\
&=& i\int_0^{\frac12} \left(e^{i(1-v)2\ol{M}} 2\hat{F} 
e^{iv(2\ol{M}+2\hat{E})}
+e^{i(1/2-v)2\ol{M}} 2\hat{F} e^{i(1/2+v)(2\ol{M}+2\hat{E})}\right)\,dv
\\
&=& i \int_0^{\frac12} e^{i(1-v)2\ol{M}} 2\hat{F} 
e^{iv(2\ol{M}+2\hat{E})}\,dv
+\int_{\frac12}^1 e^{i(1-w)2\ol{M}} 2\hat{F} e^{iw(2\ol{M}+2\hat{E})}\,dw
\\
&=& i \int_0^{1} e^{i(1-v)2\ol{M}} 2\hat{F} e^{iv(2\ol{M}+2\hat{E})}\,dv
\\
&=& \hat{F}_{{\rm exp}(i2\ol{M})}
\end{eqnarray*}
Taking adjoints, or calculating as above, the coefficient of $dA\d$ is
$$
i\int_0^1 e^{i(1-u)(2\ol{M}+2\hat{E})} (2\hat{F})\st 
e^{iu 2\ol{M}}\,du=\hat{G}_{{\rm exp}(i2\ol{M})}. 
$$
The contribution of $\hat{H}$ to the coefficient of $dt$ is
\begin{eqnarray*}
i\int_0^1 \left(e^{i(2-u)\ol{M}}\hat{H}e^{iu\ol{M}}
+e^{i(1-u)\ol{M}}\hat{H}e^{i(1+u)\ol{M}}\right)\,du
\hspace{-2in}&&
\\
&=&
i\int_0^{\frac12} \left(e^{i(1-v)2\ol{M}}2\hat{H}e^{iv2\ol{M}}
+e^{i(1/2-v)2\ol{M}}2\hat{H}e^{i(1/2+v)2\ol{M}}\right)\,dv
\\
&=&
i\int_0^{\frac12}e^{i(1-v)2\ol{M}}2\hat{H}e^{iv2\ol{M}}\,dv
+i\int_{\frac12}^1 e^{i(1-w)2\ol{M}}2\hat{H}e^{iw2\ol{M}}\,dw
\\
&=&
i\int_0^1 e^{i(1-u)2\ol{M}}2\hat{H}e^{iu2\ol{M}}\,du
\end{eqnarray*}
We now calculate the contribution of $\hat{F}$ and 
$\hat{F}\st$ to the coefficient of $dt$:
\begin{eqnarray}
\nonumber 
\left(i\int_0^1 e^{i(1-u)\ol{M}} \hat{F}e^{iu(\ol{M}+\hat{E})}\,du\right)
\left(i\int_0^1 e^{i(1-v)(\ol{M}+\hat{E})} 
\hat{F}\st e^{iv \ol{M}}\,dv\right)\hspace{-3.5in} &&
\\
&=&- \int_0^1\int_0^1 e^{i(1-u)\ol{M}} 
\hat{F}e^{i(u+1-v)(\ol{M}+\hat{E})} \hat{F}\st e^{iv \ol{M}}\,dudv
\\
\label{E 80.311}
&=&
- \int_0^1\int_1^2 e^{i(2-u)\ol{M}} 
\hat{F}e^{i(u-v)(\ol{M}+\hat{E})} \hat{F}\st e^{iv \ol{M}}\,dudv
\end{eqnarray}
The change of variable $uv\mapsto v$ gives
\begin{eqnarray}
\nonumber 
e^{i\ol{M}}\int_0^1\int_0^1 ue^{i(1-u)\ol{M}}
\hat{F}e^{iu(1-v)\ol{M}} \hat{F}\st e^{iuv\ol{M}}\,dudv
\hspace{-150pt}
&&
\\
\label{E 80.312}
&=& \int_0^1\int_0^u e^{i(2-u)\ol{M}} 
\hat{F}e^{i(u-v)(\ol{M}+\hat{E})} \hat{F}\st e^{iv \ol{M}}\,dvdu
\end{eqnarray}
Similarly
\begin{eqnarray}
\nonumber
\left(\int_0^1\int_0^1 ue^{i(1-u)\ol{M}}\hat{F}e^{iu(1-v)\ol{M}} \hat{F}\st 
e^{iuv\ol{M}}\,dudv\right)e^{i\ol{M}}
\hspace{-180pt}&&
\\
\nonumber
&=&
\int_0^1\int_0^u e^{i(1-u)\ol{M}} 
\hat{F}e^{i(u-v)(\ol{M}+\hat{E})} \hat{F}\st e^{i(1+v) \ol{M}}\,dvdu
\\
\label{E 80.313}
&=&
\int_1^2\int_1^u e^{i(2-u)\ol{M}} 
\hat{F}e^{i(u-v)(\ol{M}+\hat{E})} \hat{F}\st e^{iv \ol{M}}\,dvdu
\end{eqnarray}
Subtracting Equations (\ref{E 80.312}) and (\ref{E 80.313}) from Equation 
(\ref{E 80.311}), the 
contribution of $\hat{F}$ and $\hat{F}\st$ to the coefficient of $dt$ is
\beqn \nonumber
-\int_0^2\int_0^u e^{i(2-u)\ol{M}} 
\hat{F}e^{i(u-v)(\ol{M}+\hat{E})} \hat{F}\st e^{iv \ol{M}}\,dvdu 
\hspace{-2.5in}&&
\\
&=&
-\int_0^1\int_0^u e^{i(1-u)(2\ol{M})}
(2\hat{F})e^{i(u-v)(2\ol{M}+2\hat{E})}(2\hat{F})\st e^{iv (2\ol{M})}\,dvdu 
\eeqn
Adding this to the contribution of $\tilde{H}$ the coefficient of $dt$ is
$\hat{H}_{{\rm exp}(i2\ol{M})}$. It follows from Equation (\ref{E 80.320}) that
\beqn\label{E 80.300}
e^{i2\ol{M}_t}= I\!+\!\int_0^t (\hat{E}_{{\rm exp}(i2\ol{M})}\,d\GL
\!+\!\hat{F}_{{\rm exp}(i2\ol{M})}\,dA
\!+\!\hat{G}_{{\rm exp}(i2\ol{M})}\,dA\d\!+\!\hat{H}_{{\rm exp}(i2\ol{M})}\,ds).
\eeqn
and the quantum Duhamel formula is valid for $2M$. It follows by induction 
that the formula is valid for $2^k M$ for all $k\in \BBN$. 
\bb\ni
We have shown above that $pM$ satisfies (\ref{E 80.140}) if
 Therefore  $2^k pM$ satisfies (\ref{E 
80.140}) 
for all $k\in \BBN$ and $-1\leq p\leq 1$. Therefore (\ref{E 80.140}) is valid for all real $p$.\,\,\epf
\bb\ni
It follows from the concluding remarks of Section \ref{PDK} that 
there is a large class of $(2k+1)$-diagonal quantum semimartingales
whose integrands are constructed from $L^2$-kernels. These 
semimartingales satisfy the quantum 
Duhamel formula.  
\begin{Cor}\label{C 79.77} Let $M$ be a symmetric quantum  
semimartingale
$$
\hat{M}_t= \int_0^t \hat{E}\,d\GL+\hat{F}\,dA+\hat{F}\st 
\,dA\d+\hat{H}\,ds.
$$ 
with corresponding cmx semimartingale
$$
M_t= \int_0^t ( {E}\,d\GL+{F}\,dA+{F}\st \,dA\d+{H}\,ds).
$$ 
If the integrands $E$, $F$, $F\st$ and $H$ are 
processes of $(2k+1)$-diagonal matrices then  
$\hat{M}$ is essentially self-adjoint with closure $\ol{M}$
and $p\ol{M}$ satisfies the quantum Duhamel formula for all real $p$. 
\end{Cor}
\pf If $\psi\in \FH^j$ then $E_t{\,}^i_j \psi=(\hat{E}_t \psi)^i$ and 
$\nm{E_t{\,}^i_j \psi} =\nm{(\hat{E}_t \psi)^i}\leq 
\nm{\hat{E}_t\psi}$. Therefore $\|E_t{\,}^i_j\|\leq \nm{E_t}$ 
for all $t\in [0,1]$ and 
$\Lnu^i_j(E)\leq \nm{E}\iy$. Similarly $\Lnu^i_j(F)\leq\nm{F}_2$ 
and $\Lnu^i_j(H)\leq\nm{H}_1$ for all $i,j$.  
If $\xi=\nm{E}\iy+2\nm{F}_2+\nm{H}_1$ then $\LGk$ is $(2k+1)$-diagonal 
with $\LGk^i_j\leq \xi(i+j)$ for all $i,j$.
It follows from Proposition \ref{P 76.02} that $\ell_{00}\subset 
\CA_{\Ge}(\Xi)$ for $\Ge=\xi\inv (4k^2+2k)\inv$.
\bb\ni
The required result now follows from Proposition \ref{P 76.02} 
and Corollary \ref{C 80.221}.\,\,\epf
\bb\ni
The argument in Corollary \ref{C 80.221} only uses the fact that 
(\ref{E 80.140}) is valid 
whenever $-\Ge\leq p\leq \Ge$ and does not rely on the properties of 
$\LGk$. Therefore it may be used to 
prove the following theorem.
\begin{Thm} Let $\hat{M}=\{\hat{M}_t:t\in [0,1]\}$ be a symmetric 
quantum semimartingale. If $\hat{M}_t$ is essentially 
self-adjoint and the quantum Duhamel formula 
{\rm (\ref{E 80.140})} is valid for $-1\leq p\leq 1$ then 
it is valid for all real $p$.
\end{Thm} 
The {\em Brownian quantum semimartingale} 
$B=A+A\d$ is an essentially self-adjoint process in $\CS\pr$ 
which does not belong to $\CS_{\rm sa}$. The operators $B_t$ are 
unbounded with common core $\CD(N^{1/2})$, where $N$ is the 
number operator.
Let $\CW:L^2(\GW,\BBP^{W})\hraa \FH$ be the Wiener--Ito 
isomorphism. Let $\hat{B}=(\hat{B}_t:t\in [0,1])$ be Brownian motion. 
Abusing notation, let $\hat{B}_t$ also denote the unbounded 
self-adjoint operator in $L^2(\GW,\BBP^{W})$ of multiplication by $\hat{B}_t$. 
Then $\CW\circ \hat{B}_t\circ\CW\inv=\ol{B}_t$ for $t\in [0,1]$.
\bb\ni
The quantum semimartingale $\hat{B}=A+A\d$ with quadruple $(0,I,I,0)$ 
satisfies the condition of Corollary
\ref{C 79.77} so that $\ol{B}$ satisfies the 
quantum Duhamel formula 
\beqn\label{E 79.700}
e^{i\ol{B}_t}= I+i\int_0^t e^{i\ol{B}_s}\,dB_s-\frac12\int_0^t 
e^{i\ol{B}_s}\,ds.
\eeqn
This is equivalent, {\em via} the Wiener--Ito isomorphism, to the classical
Ito formula for the function $x\mapsto e^{ix}$ and classical 
Brownian motion $W=\{W_t:t\in [0,1]\}$:
$$
e^{iW_t}= I+i\int_0^t e^{iW_s}\,dW_s-\frac12\int_0^t 
e^{iW_s}\,ds,                         
$$

\section{\label{QIF} Quantum Ito Formula}
As in \cite{Vin2} the quantum Duhamel formula in Section \ref{QDF} 
leads, {\em via} the Fourier functional calculus, to a functional 
Ito formula. 
\bb\ni
If $\hat{T}$ is a 
self-adjoint linear transformation in $\FH$ and 
$f$ and its Fourier transform $\hat{f}$ lie in $L^1(\BBR)$ then 
define $f(\hat{T})$ 
in $\CB(\FH)$ by the Bochner integral
\beqn\label{E 81,120}
f(\hat{T})=\int_{-\infty}\iiy \hat{f}(p)e^{ip\hat{T}}\,dp.
\eeqn
This definition of $f(\hat{T})$ agrees with the usual definition {\em via } the 
functional calculus.
\bb\ni
It follows from the proof of Corollary \ref{C 79.77} and the remarks 
preceding it that there is a class of 
non-trivial quantum semimartingales, whose integrands are 
$(2k+1)$-diagonal processes defined by $L^2$-kernels, 
satisfying the conditions of the following theorem.
\bb\ni
Recall that $C^{2+}(\BBR)=\{f\in 
L^1(\BBR):p\mapsto p^2\hat{f}(p) \mbox{ belongs to }L^1(\BBR)\}$.
\begin{Thm}\label{T 81.361}
Let $\hat{M}=\{\hat{M}_t:t\in [0,1]\}$ be a symmetric quantum 
semimartingale 
$$
\hat{M}_t=\int_0^t (\hat{E}_s \,d\GL_s+ \hat{F}_s \,dA_s
+\hat{F}\st\,dA\d_s+\hat{H}_s\,ds)
$$
with corresponding symmetric cmx semimartingale $M=\{M_t:t\in 
[0,1]\}$: 
$$
M_t=\int_0^t (E_s \,d\GL_s+ F_s \,dA_s+F\st\,dA\d_s+H_s\,ds).
$$
If there exists $\Ge>0$ such that $\ell_{00}\subset 
\CA_{\Ge}(\LGk)$ where $\LGk$ is the control matrix of $M$ then
$\hat{M}$ is essentially self-adjoint with core $\FH_{00}$.
\bb\ni
For each $f\in C^{2+}(\BBR)$ the process $f(\ol{M})=\{f(\ol{M}_t):t\in [0,1]\}$ 
defined by the formula {\rm (\ref{E 81,120})} is a regular 
quantum semimartingale which satisfies the quantum Ito formula:
\beqn\label{E 81.667}
f(\ol{M}_t)=f(0)+\int_0^t(\hat{E}_{f(M)}\,d\GL
+\hat{F}_{f(M)}\,dA+\hat{G}_{f(M)}\,dA\d+\hat{H}_{f(M)}\,ds)
\eeqn
where
\begin{eqnarray*}
\hat{E}_{f(M)}=\int_{-\infty}\iiy \hat{f}(p) \hat{E}_{{\rm exp}(ipM)}\,dp
&=& f(M+\hat{E})-f(M)
\\
\hat{F}_{f(M)}=\int_{-\infty}\iiy \hat{f}(p) \hat{F}_{{\rm exp}(ipM)}\,dp
&=&\int_{-\infty}\iiy \int_0^1 ip\hat{f}(p)e^{ip(1-u)M} 
\hat{F} e^{ipu(M+\hat{E})}\,du\,dp
\\
\hat{G}_{f(M)}=\int_{-\infty}\iiy \hat{f}(p) \hat{G}_{{\rm exp}(ipM)}\,dp
&=&\int_{-\infty}\iiy \int_0^1 ip\hat{f}(p)e^{ip(1-u)(M+\hat{E})} 
\hat{F}\st e^{ipu M}\,du\,dp
\\
\hat{H}\st_{f(M)}=\int_{-\infty}\iiy \hat{f}(p) 
\hat{H}_{{\rm exp}(ipM)}\,dp
&=&\int_{-\infty}\iiy \int_0^1 ip\hat{f}(p)e^{ip(1-u)M} \hat{H} 
e^{ipu M}\,du\,dp
\\
&&\hspace{-125pt} -\int_{-\infty}\iiy\int_0^1\int_0^1
p^2\hat{f}(p)ue^{ip(1-u)M}\hat{F}e^{ipu(1-v)(M+\hat{E})} \hat{F}\st 
e^{ipuv M}\,du\,dv\,dp
\end{eqnarray*} 
\end{Thm}
\pf
It follows from Corollary \ref{C 80.221} that 
\begin{eqnarray}
\nonumber
f(M_t)&=& f(0)+\int_{-\infty}\iiy\hat{f}(p)\int_0^t(E_{{\rm exp}(ipM)}\,d\GL
+F_{{\rm exp}(ipM)}\,dA
\\
\label{E 81.211}
&&\hspace{72pt}
+G_{{\rm exp}(ipM)}\,dA\d+H_{{\rm exp}(ipM)}\,ds)\,dp.
\end{eqnarray}
so that, formally, the theorem may be proved by interchanging the order of 
integration in (\ref{E 81.211}). This requires a Fubini type theorem for 
regular quantum semimartingales. The proof is the same as when $M$ itself 
is a regular quantum semimartingale \cite[Theorem 6.1, Theorem 
6.2]{Vin2}.\,\,\epf
\bb\ni
{\em Remark.} The proof of \cite[Theorem 6.1]{Vin2} is correct for 
the functions $e,f,g,h$ defined by the formulae \cite[(48)]{Vin2}. 
The proof under the given conditions is not quite correct. 
The sixth equality in that proof is not necessarily 
valid under hypothesis (d). However this may be remedied by 
observing that in the seventh equality the integral with respect to 
$p$ may also be taken outside of the inner product. 
\bb\ni
In order to mimic the classical Ito formula we recast Equation \ref{E 81.211}
in terms of the differential $Df(M)(\cdot)$ and the Ito second differential 
$D^2_I f(M)(\cdot,\cdot)$ of $f$ at $M$.
\bb\ni
If $\hat{H}\in \CB(\FH)$ is a bounded self-adjoint transformation then standard 
perturbation theory (Hille-Phillips)\cite{HiP} shows that 
$i(\hat{T}+\hat{H})$ is the generator of a 
one parameter unitary group $p\mapsto e^{ip(\hat{T}+\hat{H})}$ so that 
$f(\hat{T}+\hat{H})$ is 
defined by (\ref{E 81,120}) with $\hat{T}$ replaced by 
$\hat{T}+\hat{H}$. We show that $f$ has a differential at 
$\hat{T}$ in the 
direction of $H$.
\begin{Propn}\label{P 81.406}
Let $\hat{T}$ be a self-adjoint operator and let $\hat{H}$ and $\hat{K}$ be bounded self-adjoint 
operator in $\FH$. If $f\in C^{2+}(\BBR)$ then the first differential 
$Df(\hat{T})(\hat{H})$ and the second differential $D^2 f(\hat{T})(\hat{H},\hat{K})$ both exist and are 
given by the formulae
\begin{eqnarray}
\label{E 81.321} Df(\hat{T})(\hat{H})&=&\int_{-\infty}\iiy \int_0^1 ip\hat{f}(p)e^{ip(1-
u)\hat{T}}\hat{H}e^{ipu\hat{T}}\,dudp;        
\\
\nonumber
Df(\hat{T})(\hat{H},\hat{K})&=&\int_{-\infty}\iiy \int_0^1\int_0^1-p^2\hat{f}(p)
\left(e^{ip(1-u)\hat{T}}\hat{H}e^{ipu(1-v)\hat{T}}\hat{K}e^{ipuv \hat{T}}\right.
\\
\label{E 81.322} 
&&\hspace{0.5in}\left. +e^{ip(1-u)\hat{T}}\hat{K}e^{ipu(1-v)\hat{T}}\hat{H}e^{ipuv \hat{T}}\right) \,dudvdp.
\end{eqnarray}
\end{Propn}
\pf Consider the Duhamel formula
$$
e^{i(\hat{T}+\Ge \hat{H})}=e^{i\hat{T}}+\Ge \int_0^1 e^{i(1-u)\hat{T}} \hat{H} e^{iu\hat{T}}\,du+
\Ge^2\int_0^1 \int_0^1 e^{i(1-u)\hat{T}} \hat{H}e^{iu(1-v)\hat{T}}\hat{H} e^{iuv(\hat{T}+\Ge \hat{H})}\,dudv.
$$
Since $\nm{e^{i(1-u)\hat{T}}}$, $\nm{e^{iu(1-v)\hat{T}}}$ and 
$\nm{e^{iuv\hat{T}}}$ are all $1$ it follows that 
$$
\nm{ \frac{e^{i(\hat{T}+\Ge \hat{H})}-e^{i\hat{T}}}{\Ge}- \int_0^1 e^{i(1-u)\hat{T}} \hat{H} e^{iu\hat{T}}\,du}
\leq \Ge\nm{\hat{H}}^2.
$$
Therefore
$$
\nm{\frac{f(\hat{T}+\Ge \hat{H})-f(\hat{T})}{\Ge}-
\int_{-\infty}\iiy \int_0^1 ip\hat{f}(p)e^{ip(1-u)\hat{T}}\hat{H}e^{ipu\hat{T}}\,dudp}
\leq \Ge\nm{\hat{H}}^2\int_{-\infty}\iiy |p\hat{f}(p)|\,dp
$$
which tends to zero with $\Ge$. This proves (\ref{E 81.321}). A similar 
argument applied to the right hand side of Equation (\ref{E 81.321}) may be used
to prove (\ref{E 81.322}).\,\,\epf
\bb\ni
These differentials may be extended to all $\hat{H}$ and $\hat{K}$ in 
$\CB(\FH)$ as 
follows: If $\hat{H}=-\hat{H}\st$ and $\hat{K}=-\hat{K}\st$ then let
\begin{eqnarray*}
iDf(\hat{T})(\hat{H})&=&Df(\hat{T})(i\hat{H});
\\
D^2f(\hat{T})(\hat{H},\hat{K})=iD^2f(\hat{T})(\hat{H},i\hat{K})&=&D^2f(\hat{T})(i\hat{H},i\hat{K})=iD^2 f(i\hat{H},\hat{K}).
\end{eqnarray*}
Each operator in $\CB(\FH)$ is, uniquely, the sum of a self-adjoint 
and skew 
adjoint transformation so these definitions extend uniquely by 
linearity and bilinearity. The {\it Ito second differential }, 
$D^2_I$, is defined 
{\em via} the matrix operator equation
\begin{eqnarray}
\nonumber
\frac12 D^2 f
\left(\mx{\hat{X}}{0}{0}{\hat{X}}\right)\left(\mx{0}{\hat{H}}{\hat{K}}{0}, \mx{0}{\hat{H}}{\hat{K}}{0}\right)
\hspace{-170pt}&&
\\
\label{E 77.224}
&=& 
\mx{D^2_I f(\hat{X})(\hat{H},\hat{K})}{0}{0}{D^2_I f(\hat{X})(\hat{K},\hat{H})}. 
\end{eqnarray}
$D^2_I f$ is an unsymmetrised version of $D^2 f$ and $D^2 f(\hat{H},\hat{K})=
D^2_I f(\hat{H},\hat{K})+D^2_I f(\hat{K},\hat{H})$ \cite[\S 
2.]{Vin2}.
\bb\ni
We use the following notation:
\begin{eqnarray*}
Df(\ol{M})(d\hat{M})&=& Df(\ol{M})(\hat{E})\,d\GL+ Df(\ol{M})(\hat{F})\,dA
\\&&
\hspace{30pt}
+Df(\ol{M})(\hat{G})\,dA\d+Df(\ol{M})(\hat{H})\,dt,
\\
D^2_I f(\ol{M})(d\hat{M},d\hat{M})&=&D^2_I f(\ol{M})(\hat{E},\hat{E})\,d\GL
+D^2_I f(\hat{M})(\hat{F},\hat{E})\,dA
\\
&&\hspace{30pt} +D^2_I f(\ol{M})(\hat{E},\hat{G})\,dA\d
+ D^2_I f(\ol{M})(\hat{F},\hat{G})\,dt
\end{eqnarray*}
\bb\ni
The quantum Ito formula is most heuristically satisfying 
when $\hat{E}\equiv 0$. In this case it is a
generalisation of the classical Ito formula for 
Brownian (continuous) semimartingales. If $\hat{E}\not\equiv 0$ 
the more complicated formula 
\cite[Theorem 6.2.]{Vin2}, which generalises the classical Ito formula 
for some discontinuous semimartingales, is also true. 
\begin{Cor}[Quantum Ito Formula]\label{T 81.362}
Suppose the symmetric quantum semimartingale
$$
\hat{M}_t=\int_0^t (\hat{F}_s \,dA_s+\hat{F}\st\,dA\d_s+\hat{H}_s\,ds)
$$
satisfies the conditions of Theorem \ref{T 81.361}.
\bb\ni
If $f\in C^{2+}(\BBR)$ then $f(\ol{M})=\{f(\ol{M}_t):t\in [0,1]\}$ 
is a regular quantum semimartingale and satisfies the quantum Ito 
formula
\beqn\label{Ito}
f(\ol{M}_t)&=&f(0)+\int_0^t(Df(\ol{M}_s)(d\hat{M}_s)+D^2_I 
f(\ol{M}_s)(d\hat{M}_s,d\hat{M}_s))
\eeqn
\end{Cor}
\pf This theorem is a translation of Theorem \ref{T 81.361} using the above 
terminology together with Proposition \ref{P 81.406}\,\,\epf
\bb\ni
The Brownian quantum semimartingale 
$\hat{B}_t=A_t+A\d_t$, satisfies the conditions of 
Corollary \ref{T 81.362}.  It follows from 
(\ref{Ito}), or directly form Equation (\ref{E 79.700}), that 
$$
f(\ol{B}_t)=f(0)+\int_0^t f\pr(\ol{B}_s)\,d\hat{B}_s
+\frac12 \int_0^t f\prpr(\ol{B}_s)\,ds
$$
whenever $f\in C^{2+}(\BBR)$. The analysis in \cite[Section 7]{Vin2} and 
\cite{Att} shows that Equation (\ref{Ito}) is equivalent to 
the classical Ito formula for Brownian motion
$$
f(W_t)=f(0)+\int_0^t f\pr(W_s)\,dW_s
+\frac12 \int_0^t f\prpr(W_s)\,ds
$$
Thus we have extended the quantum Ito formula for regular quantum 
semimartingales to a non-trivial class of essentially self adjoint quantum 
semimartingales containing the Brownian quantum 
semimartingale.
\begin{Exl} \mbox{} \end{Exl}
For each real $\Gth$ the quantum semimartingale 
$P_t(\Gth)=e^{i\Gth(t)}A_t+e^{-i\Gth(t)}A\d_t$ is representable as a 
Brownian motion. It has a {\it conjugate } quantum 
semimartingale $Q_t(\Gth)=i(e^{i\Gth(t)}A_t-e^{-i\Gth(t)}A\d_t)$, 
sometimes called the ``upside-down'' quantum Brownian motion. 
$Q(\Gth)$ is the conjugate of $P(\Gth)$ by the Fourier--Wiener 
transform: $Q_t(\Gth)=\CF P_t(\Gth) \CF\inv$ \cite[IV \S 2.3, VI  
\S 1.13]{MEY}.
\bb\ni
The case where $\Gth(t)$ is a differentiable function of $t$ is 
interesting.  Let $P_t=P_t(\Gth(t))$ and $Q_t=Q_t(\Gth(t))$. 
Then $P=\{P_t:t\in [0,1]\}$ and $Q=\{Q_t:t\in [0,1]\}$ are  
essentially self-adjoint quantum semimartingales with core $\CE$. 
Using the characterisation (\ref{E 77.31}) of quantum stochastic 
integrals it is easy to check that
$$
dP_t= e^{i\Gth(t)}dA_t+e^{-i\Gth(t)}dA\d_t+
\left(e^{i\Gth(t)}A_t-e^{-i\Gth(t)}A\d_t\right)\Gth\pr(t)\,dt.
$$
$P$ satisfies  the conditions of Corollary \ref{C 79.77}. 
The differential form of the quantum Duhamel formula is
\begin{eqnarray*}
de^{iP_t} &=& e^{iP_t}dP_t- \frac12 e^{iP_t}\,dt
+i \Gth\pr(t)\int_0^1 e^{i(1-u)P_t}(e^{i\Gth(t)}A_t-e^{-i\Gth(t)}A\d_t)
e^{iuP_t}du \,dt \\
&=& e^{iP_t}dP_t- \frac12 e^{iP_t}\,dt
+ \Gth\pr(t)\int_0^1 e^{i(1-u)P_t}Q_t 
e^{iuP_t}du \,dt.
\end{eqnarray*}
The quantum Ito formula for $f\in C^{2+}(\BBR)$ is 
$$
d(f(\hat{P}_t))= f\pr(\hat{P}_s)\,d\hat{P}_s
+\frac12  f\prpr(\hat{P}_s)\,ds + \Gth\pr(s)Df(P_s)(Q_s)ds
$$
This is the classical Ito differential continuously perturbed by the last term
which is $\Gth\pr(t)$ times the differential of $f$ at $P_t$ 
with magnitude $Q_t$ in the direction of the conjugate 
Brownian quantum semimartingale to $P_t$.

\section{\label{QSF} Quantum Stratonovich Formula.}
Quantum Stratonovich product formulae arise naturally in the 
construction of models in quantum optics \cite[Section II]{GC}. 
They also appear in the analysis of the quantum Liouville equation 
\cite[{\it op. cit.}]{SA}.
\bb\ni
Care should be taken in interpreting the statement following
\cite[Equation 3.14.]{CG}. This could be thought to imply
that the correction term in the functional is the symmetrised 
second differential rather than the correct unsymmetrised Ito 
second differential. This interpretation would lead to an incorrect 
statement of the functional quantum Stratonovich formula.
\bb\ni
The functional quantum Stratonovich 
formula for regular quantum semimartingales may be deduced from the 
functional quantum Ito formula. The details will be published elsewhere.
\begin{Thm}\label{T 86.01} Let 
$$
\hat{M}_t=\int_0^t \hat{F}_s dA_s+\hat{F}\st_s\,dA\d_s+\hat{H}_s\,ds
$$
be a gauge-free regular self-adjoint quantum semimartingale. If $f\in 
C^{3+}(\BBR)$ then
\beqn\label{E 83.38}
f(\hat{M}_t)=f(0)+\int_0^t\!\!\!\!\!\!\!{\large \sf S} 
Df(\hat{M}_s)(d\hat{M}_s)
\eeqn
\end{Thm}
We have used the following notation. 
$$
C^{3+}(\BBR)=\{f\in L^1(\BBR):p\mapsto p^3\hat{f}(p)
\mbox{ belongs to } L^1(\BBR)\}.
$$
The {\em quantum Stratonovich integral}
$$
\int_0^t\!\!\!\!\!\!\!{\large \sf S} Df(\hat{M}_s)(d\hat{M}_s)
={\rm so-}\!\!\!\lim_{|\CP|\row 0}\sum_{k=0}^{n-1}
Df\left(\frac12\left(\hat{M}_{t_{k+1}}+\hat{M}_{t_k}\right)\right)
(\hat{M}_{t_{k+1}}-\hat{M}_{t_k})
$$
where $\CP$ is the partition $0=t_0<t_1 \cdots<t_{n}=t$ 
with $ |{\cal P}|$ $=\max\{|t_{k+1}-t_{k}|:k=0,n-1\}$ and the 
convergence is in the strong operator topology.
\bb\ni
Care should be taken in interpretation of the Stratonovich 
integral since, in general, $(A_{t_{k+1}}-A_{t_k})$ and
$(A\d_{t_{k+1}}-A\d_{t_k})$ do not commute with 
$(\hat{M}_{t_{k+1}}+\hat{M}_{t_k})$ and 
$$
Df(\hat{M}_s)(d\hat{M}_s)\not\equiv
Df(\hat{M}_s)(\hat{F}_s)\,dA_s+
Df(\hat{M}_s)(\hat{F}\d_s)\,dA\d_s+Df(\hat{M}_s)(\hat{H}_s)\,ds.
$$
The proof of Theorem \ref{T 86.01} which uses the Taylor series 
expansion of $f:\CB(\FH)\row \CB(\FH)$ to three terms with remainder 
depends on the regularity of $\hat{M}$. The control of the remainder
when $\hat{M}$ is irregular becomes a serious obstruction.
\bb\ni
The product formula may be recovered from the functional 
formula, with $f(x)=x^2$, using $2\times 2$ matrix valued processes 
as in \cite[\S 3]{Vin2}.
\bb\ni
It may also be shown that if 
$M=\int_0^t {F}_s dA_s+{F}\st_s\,dA\d_s+{H}_s\,ds$ is a 
symmetric cmx semimartingale satisfying the conditions of Theorem
\ref{T 81.361} and $f$ is a polynomial then 
\beqn\label{E 86.767}
f({M}_t)=f(0)+\int_0^t\!\!\!\!\!\!\!{\large \sf S} Df({M}_s)(d{M}_s)
\eeqn
where
$$
\int_0^t\!\!\!\!\!\!\!{\large \sf S} Df({M}_s)(d{M}_s)
={\rm se-}\!\!\!\lim_{|\CP|\row 0}\sum_{k=0}^{n-1}
Df\left(\frac12\left({M}_{t_{k+1}}+{M}_{t_k}\right)\right)
({M}_{t_{k+1}}-{M}_{t_k})
$$
and the se-limit is the {\em strong entrywise} limit: the limit 
in the strong operator topology in each matrix entry.
\bb\ni
{\it Problem 1. Is the cmx Stratonovich formula (\ref{E 
86.767}) valid for $f\in C^4(\BBR)\cap L^1(\BBR)$ and does it 
give rise to a corresponding quantum Stratonovich formula ?} 

\section{\label{PCM} Perturbation of Classical Semimartingales}
For the remainder of this article we will not be concerned with 
chaos matrices and the `` $\hat{}$ '' used to 
distinguish operator from chaos matrix will be dropped.
\bb\ni
It follows from Section \ref{PDK} and Corollary \ref{C 
79.77} that there is a large class of irregular 
essentially self-adjoint quantum
semimartingales which satisfy the quantum Duhamel formula and 
consequently the quantum Ito formula. The next three sections 
are devoted to showing that this class is closed under 
perturbation by regular self-adjoint quantum semimartingales.
\bb\ni
Theorem \ref{T 81.717}, the main result of this section,
is a special case of Theorem \ref{TT 85.01}. However the proof 
of the former is more sophisticated than 
the proof of the latter, which uses rather crude expansion 
methods.
\bb\ni
Let ${W}=\{{W}_t:t\in [0,1]\}$ be Brownian motion on 
$(\GW,\BBP)$, the classical Wiener space for the time interval 
$[0,1]$ endowed with the Brownian filtration. If $f\in L^2[0,1]$ 
let ${\psi}(f)$ denote the exponential martingale:
$$
{\psi}(f)=e^{\int_0^1 f(s)\,d{W}_s-\frac12 \int_0^1 
f(s)^2\,ds}
$$
If $\CW:L^2(\GW,\BBP)\hraa \FH$ is the 
Wiener--Ito isomorphism then $e(f)=\CW{\psi}(f)$.  
\bb\ni
Let ${\sf F}$ 
be a bounded real-valued adapted stochastic 
process on $(\GW,\BBP)$. 
Then the martingale ${\sf M}_t =\int_0^t {\sf F}_s\,d{W}_s$  
satisfies 
the classical Ito formula:
\beqn\label{E 80.080}
f({\sf M}_t)&=&f(0)+\int_0^t f\pr({\sf M}_s)\,d{\sf M}_s
+\frac12 \int_0^t f\prpr({\sf M}_s\,)d\la {\sf M} \ra_s,
\eeqn
where $d{\sf M}_s={\sf F}_s \,d{W}_s$ and $d\la {\sf M}\ra_s= 
{\sf F}_s^2\,ds$.
\bb\ni
If ${\sf F}=I$ then $M=\CW\circ {\sf M}\circ \CW\inv$ is 
the Brownian quantum semimartingale 
For more general $\sf F$ the cmx representation of $F=\CW\circ {\sf F}\circ 
\CW\inv$ is complicated and does not satisfy the conditions of 
Corollary \ref{T 81.362} and the classical Ito formula 
(\ref{E 80.080}) does not directly follow from (\ref{Ito}). 
\bb\ni
Using a martingale moment inequality \cite[Exercise 3.25]{KS} 
$$
\BBE[{\sf M}_t^4]
\leq 6\BBE\left[\int_0^t |{\sf F}_s^4|\,ds\right]\leq
6\BBE\left[\int_0^1 |{\sf F}_s^4|\,ds\right]
$$
so that ${\sf M}_t\in L^4(\GW,\BBP)$ for $t\in [0,1]$. It 
follows from direct calculation that ${\psi}(f)\in 
L^4(\GW,\BBP)$.

However ${\sf M}$ is a martingale and $M=\{M_t:t\in [0,1]\}$ 
defined by
$$
M_t=\int_0^t F\,dA+F\,dA\d
$$
is a quantum semimartingale that satisfies the Ito formula 
$$
f(\ol{M}_t)=f(0)+\int_0^t f\pr(\ol{M}_s) 
F_s(dA_s+dA\d_s)+\frac12\int_0^t f\prpr(M_s)F_s^2 \,ds
$$
for all bounded $C^2$ functions $f$. 
By \cite[VIII.3 Proposition 2]{RS1} with $p=4$ it follows that 
for each $t\in [0,1]$ the operator of multiplication by 
${\sf M}_t$ is self-adjoint on its maximal domain and that 
$\CE_{W}$ is a core for ${\sf M}_t$.
\bb\ni
These facts may be used to 
show that perturbations of $M$ by regular quantum semimartingales 
satisfy the quantum Ito formula.
\begin{Thm}\label{T 81.717} Let $N=N(R,S+F,S\st+F,U)$ be the quantum 
semimartingale $M+J$ where $J$ 
is a regular self-adjoint quantum semimartingale: 
$$
J_t=\int_0^t (R\,d\GL+S\,dA+S\st\,dA\d+U\,ds).
$$
Then $pN$ is essentially self-adjoint and satisfies the 
quantum Duhamel formula, {\rm (\ref{E 80.140})}, for each real $p$.
\end{Thm}
\begin{Cor}\label{C 81.747} If $f\in C^{2+}(\BBR)$ then $f(\ol{N})$ satisfies the 
quantum Ito formula {\rm (\ref{E 81.667})}.
\end{Cor}
We will only outline the proof, which follows that 
for the special case $F=I$ in \cite{Vin3}. We shall need a 
result from that paper.
\begin{Propn}\label{P 81.111} 
Let 
$N^{[n]} \equiv N^{[n]} (E^{[n]} ,F^{[n]} ,G^{[n]} ,H^{[n]} )$, 
be a sequence of regular quantum semimartingales and let
$N=N(E,F,G,H)$ be a quantum semimartingale such that:
\bb\ni
{\rm (i)} $E^{[n]} _t$, $F^{[n]} _t$, 
$G^{[n]} _t$ and $H^{[n]} _t$ converge strongly to   
$E_t$, $F_t$, $G_t$ and $H_t$ respectively for almost all $t\in 
[0,1]$;
\bb\ni
{\rm (ii)} $\sup\{\|N^{[n]} _t\|: t\in [0,1],\,\,n=1,2,\ldots\}=K<\infty$;
\bb\ni
{\rm (iii)} $\sup_n\|E^{[n]} _t\|= \Ga(t)\in L\iiy[0,1]$;
\bb\ni
{\rm (iv)} $\sup_n(\|F^{[n]} _t\|+\|G^{[n]} _t\|)= \Gb(t)\in L^2[0,1]$;
\bb\ni
{\rm (v)} $\sup_n\|E^{[n]} _t\|= \Gg(t)\in L^1[0,1]$.
\bb\ni
Then
\bb\ni 
{\rm (a)} $N$ is regular with $\nm{N_t}\leq K$ for all $t\in [0,1]$;
\bb\ni
{\rm (b)} $N^{[n]} _t$ converges strongly to $N_t$ for all $t\in [0,1]$.
\end{Propn}
We also use the following corollary of Trotter's theorem. It is 
a direct consequence of \cite[Theorems VIII.21, VIII.25(a)]{RS1}. 
\begin{Thm}\label{C 03.01}
Let $T_n$ and $T$ be self-adjoint operators and 
suppose that $\CD$ is a common core for $T$ and all $T_n$. If $T_n\Gv\row 
T\Gv$ for each $\Gv$ in $\CD$ then $e^{ipT_n}$ converges 
strongly to $e^{ipT}$ for each real $p$.
\end{Thm}
{\em Proof of Theorem \ref{T 81.717}}. 
It is enough to consider the case $p=1$ and $\nm{\hat{F}_t}\leq 1$ 
for all $t$. Then $N_t=M_t+J_t$ 
is essentially self--adjoint on $\CE$ by the Kato--Rellich theorem 
\cite[Theorem X.12]{RS1}. Since $J_t$ is bounded 
$\ol{N}_t=\ol{M}_t+\ol{J}_t$, where operators are added on their 
common domain.
\bb\ni
For each $n=1,2,\ldots$ choose a $C\iiy$ function 
$h_n:\BBR\row [0,\infty)$ with $h_n=0$ on $(-\infty,-4n]\cup 
[4n,\infty)$, with $h_n(x)=x$ for $x\in [-n,n]$ and with 
$-1\leq h_n\pr,h_n\prpr\leq 1$. By the classical Ito formula 
\cite[IV 3.]{RoW} 
$$
h_n({\sf M})_t=\int_0^t h\pr_n({\sf M}_s)F_s\,d{W}_s +
\frac12 \int_0^t h_n\prpr({\sf M}_s)\hat{F}_s^2\,ds. 
$$
Moreover the multiplication operators $h_n({\sf M})_t$, 
$h_n\pr({\sf M})_t$ and $h_n\prpr({\sf M})_t$ are bounded with
$\nm{h_n({\sf M})_t}\leq n$ and $\nm{h\pr_n({\sf M})_t},\,\,
\nm{h\prpr_n({\sf M})_t}\leq 1$. 
Applying the Wiener--Ito isomorphism, the restriction of 
$h_n(\ol{M})$ to $\CE$ is a regular essentially self-adjoint quantum 
semimartingale and, for $\psi\in \CE$, 
$$
h_n(\ol{M})_t\psi=\int_0^t (h\pr_n(\ol{M}_s)F_s(dA_s+dA\d_s) +
\frac12 h_n\prpr(\ol{M}_s)F_s^2\,ds)\psi. 
$$
The processes  
$$
F\ub{n}=S+h\pr_n(\ol{M})F,\quad G\ub{n}=S\st 
+h\pr_n(\ol{M})F,\quad H\ub{n}=U+ \frac12 h_n\prpr(\ol{M}_s)F^2
$$ 
are strongly measurable and 
$$
\nm{F\ub{n}_t}, \nm{G\ub{n}_t}\leq \nm{F_t}+1;\qquad
\nm{H\ub{n}_t}\leq \nm{U_t}+1/2;
$$
$$(F\ub{n})\st=G\ub{n}; \qquad(H\ub{n})\st=H\ub{n}.
$$ 
Therefore 
$$
N\ub{n}=\int_0^t 
(E_s\,d\GL_s+F\ub{n}_s\,dA_s+G\ub{n}_s\,dA\d_s+H\ub{n}_s\,ds)
$$
is a regular self-adjoint quantum semimartingale 
with 
$$
\ol{N}\ub{n}=\ol{J}+h_n(\ol{M})
$$ 
Let $\psi\in \CE$. Since $N\ub{n}$ and $R$ are bounded $\CE$ 
is a common core for $M_t$, $N_t$, $N\ub{n}_t$ and $R_t$ 
for all $t\in [0,1]$ and $n\in \BBN$.   
By the spectral theorem \cite[Theorem VIII.5 (c)]{RS1} 
$\lim_{n\roy}h_n(\ol{M}_t)\psi=\ol{M}_t\psi$. 
Since $J$ is regular $N\ub{n}_t\Gv=h_n(\ol{M}_t)\Gv+\ol{J}_t\Gv$ 
converges to $\ol{N}_t\Gv$. Similarly $(N\ub{n}_t+R_t)\Gv$ converges 
to $(\ol{N}_t+R_t)\Gv$ 
\bb\ni
It follows from Theorem \ref{C 03.01} that, for each real $p$, 
\beqn
\lim_{n\roy}e^{ipN\ub{n}_t} = e^{ip\ol{N}_t},&\qquad &
\lim_{n\roy}e^{ip(N\ub{n}_t+R_t)} = e^{ip(\ol{N}_t+R_t)},
\eeqn
in the strong operator topology. 
\bb\ni
Now $\nm{h_n\pr}\iy, \nm{h_n\prpr}\iy\leq 1$ and 
$h\pr_n$ and $h\prpr_n$ converge pointwise to $1$ and $0$ respectively. 
By the spectral theorem 
\cite[Theorem VIII.5 (d)]{RS1} $h_n\pr(\ol{M}_t)$ is strongly convergent 
to $I$ and $h_n\prpr(\ol{M}_t)$ is strongly convergent to $0$. Therefore
\beqn
\hspace{-0.25in}\lim_{n\roy}F\ub{n}_t =S_t+F_t,\quad
\lim_{n\roy}G\ub{n}_t =S\st_t+F_t,\quad
\lim_{n\roy}H\ub{n}_t =U_t,
\eeqn
in the strong operator topology. 
\bb\ni
The regular quantum semimartingale $N\ub{n}$ satisfies the 
quantum Duhamel formula, \cite[Proposition 5.1]{Vin2},
for $p=1$:
$$
e^{iN\ub{n}_t}=I+\int_0^t ( E^{[n]} \,d\GL+
F^{[n]} \,dA+G^{[n]} \,dA\d+
H^{[n]} \,ds),
$$
where
\begin{eqnarray*}
 E^{[n]}=E_{{\rm exp}(iN\ub{n})}(t)&=& e^{i(N_t\ub{n}+R_t)}-
e^{iN\ub{n}_t},
\\
 F^{[n]}=F_{{\rm exp}(iN\ub{n})}(t)
&=& i\int_0^1 e^{i(1-u)N\ub{n}_t} F\ub{n}_t 
e^{iu(N\ub{n}_t+R_t)}\,du,
\\
 G^{[n]}=G_{{\rm exp}(iN\ub{n})}(t)
&=& i\int_0^1 e^{i(1-u)(N\ub{n}_t+R_t)} G\ub{n}_t e^{iu 
N\ub{n}_t}\,du,
\\
 H^{[n]}=H_{{\rm exp}(iN\ub{n})}(t)
&=& i\int_0^1 e^{i(1-u)N\ub{n}_t} H\ub{n}_t e^{iu N\ub{n}_t}\,du
\\
&& \hspace{-0.5in}-\int_0^1\int_0^1
ue^{i(1-u)N\ub{n}_t}F\ub{n}_te^{iu(1-v)(N\ub{n}_t+R_t)} 
G\ub{n}_t e^{iuv N\ub{n}_t}\,du\,dv.
\end{eqnarray*} 
If $t,u\in [0,1]$ and $n\in \BBN$ then
$$
\nm{e^{i(1-u)N\ub{n}_t} F\ub{n}_t e^{iu(N\ub{n}_t+R_t)}}\leq 
\nm{F\ub{n}_t}\leq \nm{F_t}+1.
$$
Let $\psi\in \FH$. Since operator multiplication is continuous on 
norm bounded sets for the strong operator topology 
\cite[Proposition 2.4.1.]{BR1}
$$
\lim_{n\roy} e^{i(1-u)N\ub{n}_t} F\ub{n}_t e^{iu(N\ub{n}_t+R_t)}\psi
=e^{i(1-u)\ol{N}_t} F\ub{n}_t e^{iu(\ol{N}_t+R_t)}\psi.
$$
By the dominated convergence theorem for Bochner integrals 
\cite[Chapter II, Theorem 3]{DiU}
$$
\lim_{n\roy} F_{{\rm exp}(iN\ub{n})}(t)\psi=
\lim_{n\roy} F_{{\rm exp}(i\ol{N})}(t)\psi.
$$
Similarly  $E_{{\rm exp}(iN\ub{n})}(t)$, 
$G_{{\rm exp}(iN\ub{n})}(t)$ 
and $H_{{\rm exp}(iN\ub{n})}(t)$ converge strongly to
$E_{{\rm exp}(i\ol{N})}(t)$, 
$G_{{\rm exp}(i\ol{N})}(t)$ and $H_{{\rm exp}(i\ol{N})}(t)$ respectively, 
where the limits are defined by  (\ref{E 79.999}). Therefore condition (i) 
of Proposition \ref{P 81.111} is satisfied by the quantum 
semimartingales $e^{i\ol{N}}$ and $e^{iN\ub{n}}$. 
\bb\ni
Condition (ii)  of Proposition \ref{P 81.111}
is satisfied since $e^{iN\ub{n}_t}$ is unitary.
\bb\ni
Since $\nm{f_n\pr(M_t)}, \,\nm{f_n\prpr(M_t)}\leq 1$ 
the  $F\ub{n},G\ub{n}\in L^2([0,1], \CB(\FH)\so)$ and 
$H\ub{n}\in L^1([0,1], \CB(\FH)\so)$.
Since $J(R,S,S\st,U)$ is a quantum semimartingale
\begin{eqnarray*}
\nm{E_{{\rm exp}(iN\ub{n})}(t)}&\leq&2=\Ga(t)\in L\iiy[0,1],
\\
\nm{{F}_{{\rm exp}(iN\ub{n})}(t)}+
\nm{G_{{\rm exp}(iN\ub{n})}(t)}&\leq& \nm{F\ub{n}_t}+\nm{ G\ub{n}_t}
=\Gb(t)\in L^2[0,1],
\\
\nm{H_{{\rm exp}(iN\ub{n})}(t)}
&\leq& \nm{H\ub{n}_t} + \nm{F\ub{n}_t}\nm{G\ub{n}_t}=\Gg(t)\in 
L^1[0,1].
\end{eqnarray*} 
and conditions (iii), (iv) and (v)  of Proposition \ref{P 81.111}
are satisfied. 
\bb\ni
It follows from Proposition \ref{P 81.111} that  
$e^{i\ol{N}}$ is a regular quantum semimartingale and satisfies the 
quantum Duhamel formula (\ref{E 80.140}). This outlines the proof 
of Theorem \ref{T 81.717}.\,\,\epf
\bb\ni
Corollary \ref{C 81.747} follows from Theorem \ref{T 81.717} 
exactly as Theorem \ref{T 81.361} follows from Corollary
\ref{C 80.221}.
\section{\label{DHE} The Duhamel Expansion}
Another way to show that a perturbed quantum 
semimartingales satisfies the quantum Duhamel formula is to 
expand both sides of the formula using 
the classical Duhamel expansion.
\bb\ni
Let $J,J_1,\ldots,J_n$ be bounded self-adjoint operators in a 
Hilbert space $H$ and let $M$ be a self-adjoint operator in $H$. 
Define the kernel
$$
K\ub{n}(M:J_1,\ldots,J_n)\defeq
\int_{\GD\ub{n}}e^{i(1-u_1)M}J_1 e^{i(u_1-u_2)M}J_2 \cdots 
e^{i(u_{n-1}-u_n)M}J_n e^{iu_n M}du
$$
Where $\GD\ub{n}=\{1\geq u_1\geq \cdots u_j\geq \cdots\geq 
u_n\geq 0\}$ and $du= du_1\ldots du_n$.
\bb\ni
Then $K\ub{n}(M:J_1,\ldots,J_n)$ is a bounded linear operator and 
\beqn\label{EE 84.01}
\nm{K\ub{n}(M:J_1,\ldots,J_n)}&\leq &
\frac{\nm{J_1}\cdots\nm{J_n}}{n!}
\eeqn
Now put $K\ub{0}(M,J)=e^{iM}$ and $K\ub{n}(M,J)=K\ub{n}(M:J,\ldots,J)$.
\bb\ni
Using the bounds (\ref{EE 84.01}) it may be shown that
\beqn\label{EE 84.02}
e^{i(M+J)}=\sum_{n=0}\iiy i^n K\ub{n}(M,J),
\eeqn
the series being uniformly convergent in $\CB(H)$. The expansion
(\ref{EE 84.02}) is known as the Duhamel expansion of 
$e^{i(M+J)}$. A proof may be found in \cite{HiP}.
\bb\ni
If $\Phi_1,\Phi_2,\ldots $ belong to $\CB(H)$ and $k\leq n$ then 
put
$$
K\ub{n}(M,J:\stackrel{(j_1)}{\Phi}_{\!\!1},\ldots,\stackrel{(j_k)}
{\Phi}_{\!\!k})
=K\ub{n}(M:J_1,\ldots,J_n)
$$
where
$$
J_j=\left\{\ba{rl} \Phi_i&\quad \mbox{ if } j=j_i,\,\,i=1,\ldots,k, \\
J&\quad\mbox{otherwise}.
\ea\right.
$$
This definition is independent of the order of the $\Phi$s: if 
$\tau$ is a permutation of $\{1,\ldots,k\}$ then 
Thus
$$ 
K\ub{n}(M,J:\stackrel{(j_{\tau(1)})}{\Phi}_{\!\!\tau(1)},\ldots,
\stackrel{(j_{\tau(k)})}
{\Phi}_{\!\!\tau(k)})
=K\ub{n}(M,J:\stackrel{(j_1)}{\Phi}_{\!\!1},\ldots,\stackrel{(j_k)}
{\Phi}_{\!\!k})
$$
We sometimes compress this notation by 
suppressing $M$ and $J$. Thus
\begin{eqnarray*}
K\ub{n}&=&K\ub{n}(M,J), \\
K\ub{n}(w)&=&K\ub{n}(wM,wJ), \\
K\ub{n}(\stackrel{(j_1)}{\Phi}_{\!\!1},\ldots,\stackrel{(j_k)}{\Phi}_{\!\!k})
&=&K\ub{n}(M,J:\stackrel{(j_1)}{\Phi}_{\!\!1},\ldots,\stackrel{(j_k)}
{\Phi}_{\!\!k})\\
K\ub{n}(w;\stackrel{(j_1)}{\Phi}_{\!\!1},\ldots,\stackrel{(j_k)}{\Phi}_{\!\!k})
&=&K\ub{n}(wM,wJ:w\stackrel{(j_1)}{\Phi}_{\!\!1},\ldots,w\stackrel{(j_k)}
{\Phi}_{\!\!k})
\end{eqnarray*}
If $\Ga_1,\Ga_2,\Ga_3, \ldots$ are positive integers and 
$\Gs_j=\Ga_1+\cdots+\Ga_j+j$ then
\begin{eqnarray}
\label{EE 84.03}
K\ub{\Gs_2-1}(\stackrel{(\Gs_1)}{\Phi}_{\!\!1})
&= & \int_0^1(1-w)^{\Ga_1}w_1^{\Ga_2}K\ub{\Ga_1}(1-
w)\Phi_1K\ub{\Ga_2}(w)\,dw
\end{eqnarray}
\begin{eqnarray}\label{EE 84.04}
\lefteqn{K\ub{\Gs_3-
1}(\stackrel{(\Gs_1)}{\Phi}_{\!\!1},\stackrel{(\Gs_2)}{\Phi}_{\!\!2})}\\
\nonumber &=&
\int_{\GD\ub{2}}(1-w_1)^{\Ga_1}(w_1-w_2)^{\Ga_2}w_2^{\Gg}
K\ub{\Ga_1}(1-w_1)\Phi_1 K\ub{\Ga_2}(w_1-w_2)\Phi_2 K\ub{\Ga_3}(w_2)\,dw
\eeqn
These identities are special cases of the following proposition.
\begin{Propn}
If $\Ga_1,\ldots,\Ga_{n},\ldots$ are positive integers and 
$\Gs_j=\Ga_1+\cdots+\Ga_j+j$ then
\begin{eqnarray}
\label{EE 84.05} K\ub{\Gs_{n+1}-1}
(\stackrel{(\Gs_1)}{\Phi}_{\!\!1},\ldots,\stackrel{(\Gs_j+j)}{\Phi}_{\!\!j}
,\ldots,\stackrel{(\Gs_n)}{\Phi}_{\!\!n})&&
\\
\nonumber &\hspace{-2in}&\hspace{-2.5in}
=\int_{\GD\ub{n}}(1-w_1)^{\Ga_1}\cdots(w_{j-1}-w_{j})^{\Ga_j}
\cdots(w_{n-1}-w_{n})^{\Ga_n}w_n^{\Ga_{n+1}}
K\ub{\Ga_1}(1-w_1)\Phi_1\cdots 
\\
\nonumber &\hspace{-1.5in}&\hspace{-1.5in}
\cdots \Phi_{j-1}K\ub{\Ga_j}(w_{j-1}-
w_{j})\Phi_j\cdots\Phi_nK\ub{\Ga_{n+1}}(w_n)\,dw
\end{eqnarray}
\end{Propn}
\pf  The right hand side of (\ref{EE 84.05}) is
\begin{eqnarray*}
&& \int_{\GD\ub{n}\times\GD\ub{\Ga_1}\times\cdots\GD\ub{\Ga_n}}
(1-w_1)^{\Ga_1}\cdots(w_{n-1}-w_{n})^{\Ga_n}w_n^{\Ga_{n+1}}
e^{i(1-u^1_1)(1-w_1)M}J\cdots
\\
&&\hspace{0.3in}\cdots J e^{i u^1_{\Ga_1}(1-w_1)M} \Phi_1
e^{i(1-u^2_1)(w_1-w_2)M}J\cdots
J e^{iu^2_{\Ga_2})(w_1-w_2)M} \Phi_2
e^{i(1-u^3_1)(w_2-w_3)M}J\cdots
\\
&&\hspace{0.6in}\cdots J e^{iu^j_{\Ga_j})(w_{j-1}-w_j)M} \Phi_j
e^{i(1-u^{j+1}_1)(w_j-w_{j+1})M}J\cdots
\\
&&\hspace{1.6in} \cdots J e^{iu^{j+1}_{\Ga_{j+1}}(w_{j}-w_{j+1})M} \Phi_{j+1}
e^{i(1-u^{j+2}_1)(w_{j+1}-w_{j+2})M}J\cdots
\\
&&\hspace{0.3in}\cdots J e^{iu^{n}_{\Ga_{n}}(w_{n-1}-w_{n})M} \Phi_{n}
e^{i(1-u^{n+1}_1)w_n M}J\cdots Je^{iu^{n+1}_{\Ga_{n+1}}w_{n+1}M}
\,dw du^1 \ldots du^{n+1}
\end{eqnarray*}
By applying the change of variables
$$
\ba{lcll} 
\xi\ub{k}_j&=& w_k+u\ub{k}_j(w_{k-1}-w_k) &\quad w_0=1;\,\,
k=1,\ldots,n;\,\, j=1,\ldots,\Ga_k,
\\
\xi\ub{k}_{\Ga_k+1}&=& w_k &\quad k=1,\ldots,n,
\\
\xi\ub{n+1}_j&=& u\ub{n+1}_j w_n &\quad j=1,\ldots,\Ga_{n+1},
\ea
$$
we obtain the left hand side of (\ref{EE 84.05}).\,\,\epf
\bb\ni
We will also need a second order Duhamel expansion. If $J$ and 
$R$ belong to $\CB(H)$ then
$$
e^{i(M+J+R)}
=\sum_{n=0}\iiy \sum_{k=0}^n \sum_{1\leq j_1<\cdots<j_k\leq n}
K\ub{n}(\stackrel{(j_1)}{R},\ldots,\stackrel{(j_1)}{R})
$$
If $R=0$ this reduces to (\ref{EE 84.02}) so that 
\beqn\label{EE 84.06}
e^{i(M+J+R)}-e^{i(M+J)} 
=\sum_{n=0}\iiy \sum_{k=1}^n \sum_{1\leq j_1<\cdots<j_k\leq n}
K\ub{n}(\stackrel{(j_1)}{R},\ldots,\stackrel{(j_1)}{R})
\eeqn

\section{\label{PQS} Perturbation of Quantum Semimartingales}
We shall need a series version of Proposition \ref{P 81.111} and 
a lemma on quantum stochastic differentiation under the integral 
sign.
\begin{Propn}\label{PP 85.01} 
Let 
$K\ub{n}\equiv K\ub{n}(E\ub{n},F\ub{n},G\ub{n},H\ub{n})$, 
be a sequence of regular quantum semimartingales and let
$N=N(E,F,G,H)$ be a quantum semimartingale such that:
\bb\ni
{\rm (i)} $\sum_{n=0}\iiy E\ub{n}_t$, $\sum_{n=0}\iiy F\ub{n}_t$, 
$\sum_{n=0}\iiy G\ub{n}_t$ and $\sum_{n=0}\iiy H\ub{n}_t$ converge strongly to   
$E_t$, $F_t$, $G_t$ and $H_t$ respectively for almost all $t\in 
[0,1]$;
\bb\ni
{\rm (ii)} $\sum_{n=0}\iiy\|K\ub{n}_t\|=\kappa(t)\in L\iiy[0,1]$;
\bb\ni
{\rm (iii)} $\sum_{n=0}\iiy\|E\ub{n}_t\|= \Ga(t)\in L\iiy[0,1]$;
\bb\ni
{\rm (iv)} $\sum_{n=0}\iiy(\|F\ub{n}_t\|+\|G\ub{n}_t\|)= \Gb(t)\in L^2[0,1]$;
\bb\ni
{\rm (v)} $\sum_{n=0}\iiy\|E\ub{n}_t\|= \Gg(t)\in L^1[0,1]$.
\bb\ni
Then
\bb\ni 
{\rm (a)} $N$ is regular with $\nm{N_t}\leq \kappa$ for all $t\in [0,1]$;
\bb\ni
{\rm (b)} $\sum_{n=0}\iiy K\ub{n}_t$ converges strongly to $N_t$ for all $t\in [0,1]$.
\end{Propn}
\begin{Lemma}\label{LL 85.01} For each $u\in \GD\ub{n}$ let
$$
k_t (u)=\int_0^t \mbox{\scriptsize 
E} (u)\,d\GL+\mbox{\scriptsize F} (u)\,dA+\mbox{\scriptsize 
G} (u)\,dA\d+\mbox{\scriptsize H} (u)\,ds
$$
be a regular quantum semimartingale with $\nm{K_t(u)}\leq 
c<\infty$ for all $t\in [0,1]$ and $u\in \GD $. Suppose also that 
$u\mapsto k (u)$, $\mbox{\scriptsize E} (u)$, 
$u\mapsto \mbox{\scriptsize F} (u)$, $\mbox{\scriptsize G} (u)$ and 
$u\mapsto \mbox{\scriptsize H} (u)$ are 
strong operator continuous maps from $\GD\ub{n}$ to 
$L\iiy([0,1],\CB(\FH)\so)$, 
$L^2([0,1],\CB(\FH)\so)$ and $L^1([0,1],\CB(\FH)\so)$ 
respectively. Then $K =\{K_t:t\in [0,1\}$ defined by
$$
K_t =\int_{\GD } k_t (u)\,du,
$$
is a regular quantum semimartingale. Moreover 
\beqn\label{FF 85.01}
K_t=\int_0^t (E \,d\GL+F \,dA+ G \,dA\d+H \,ds),
\eeqn
where
$$
E =\int_{\GD\ub{n}}\mbox{\scriptsize E}(u) \,du,\quad
F =\int_{\GD\ub{n}}\mbox{\scriptsize F}(u) \,du
$$
$$
G =\int_{\GD\ub{n}}\mbox{\scriptsize G}(u) \,du,\quad
H =\int_{\GD\ub{n}}\mbox{\scriptsize H}(u) \,du
$$
\end{Lemma}
\pf It follows from (\ref{E 77.31})  that
\begin{eqnarray*}
\ip{K_t(u) e(f)}{e(g)} &=& \int_{\GD\ub{n}}\ip{k_t(u) 
e(f)}{e(g)}\,du
\\
&&\hspace{-0.8in}= \int_{\GD\ub{n}}\int_0^t \ip{fg \mbox{\scriptsize E}(u)+
f\mbox{\scriptsize F}(u)
+\ol{g}\mbox{\scriptsize G}(u)
+\mbox{\scriptsize H}(u))(s)e(f)}{e(g)}\,ds
\\
&&\hspace{-0.8in}= 
\int_0^t \ip{fg { E}(u)+f{ 
F}(u)+\ol{g}{G}(u)+{H}(u))(s)e(f)}{e(g)}\,ds.
\end{eqnarray*}
The conditions of the theorem have been chosen so that the 
interchange of integrals is valid. The identity (\ref{FF 85.01}) 
now follows from the characterisation of quantum stochastic 
integrals (\ref{E 77.31}).\,\,\epf
\bb\ni
Although the following theorem is true as stated below we 
shall only prove it when $E=0$.
\begin{Thm}\label{TT 85.01}
Let $M=M(E,F,F\st,H)$ be and essentially self-adjoint quantum semimartingale 
which satisfies the quantum Duhamel formula and let 
$J=J(R,S,S\st,U)$ be a 
regular quantum semimartingale. Then $M+J$ is an essentially 
self-adjoint quantum semimartingale which satisfies the quantum 
Duhamel formula.
\end{Thm}
\pf
By the Kato--Rellich theorem $M+J$ is an essentially self-adjoint 
quantum semimartingale. 
\bb\ni
To prove the quantum Duhamel formula we formally expand the 
terms on both sides using the Duhamel expansion and 
equate the resulting expressions. The formal calculations are 
then made rigorous. 
\bb\ni
The Duhamel expansion (\ref{EE 84.02}) of 
$e^{i(M_t+J_t)}$ is
\beqn \label{EE 85.1}
e^{i(M_t+J_t)}&=&\sum_{n=0}\iiy i^n K\ub{n}(M_t,J_t)
\eeqn
where $K\ub{n}=K\ub{n}(M,J)$ and 
$$
K\ub{n}=\int_{\GD\ub{n}}e^{i(1-u_1)M}Je^{i(u_1-u_2)M}\cdots 
e^{i(u_{n-1}-u_n)M}Je^{iu_n M}\,du.
$$
The integrand $k\ub{n}$ is the product of elements of $\CS$ 
and so belongs to $\CS$. The quantum stochastic differential
$dk\ub{n}$ may be calculated from the quantum Duhamel formula 
and the quantum Ito product formula (this is done  
in the computation of $dK\ub{n}$ below.) Since multiplication is 
continuous on bounded sets for the strong operator topology it 
follows that $k\ub{n}$ satisfies the conditions of Lemma \ref{LL 
85.01}. Thus $K\ub{n}$ is in $\CS$ and its differential may be 
calculated using the formula (\ref{FF 85.01}).
\bb\ni
The regular quantum semimartingale $t\mapsto \int_0^t H_s\,ds$ 
may be incorporated in $J$ so that $H$ may be assumed to be zero. 
To further reduce the notational complexity we proceed 
in two steps.
\bb\ni
Step 1. Assume that the coefficients $E$, $H$ and $R$ are 
zero. It then follows from the quantum Ito product formula that 
the product of any three quantum differentials taken from 
$\{dM,dJ\}$ is zero. Since $dJ$ and $dM$ have no terms in $d\GL$ 
we must prove the simple form of the quantum Duhamel formula:
\begin{eqnarray}
\label{EE 85.01}
e^{i(M_t+J_t)}&=&I+
\int_0^t i \int_0^1 e^{i(1-u)(M+J)}d(M+J) e^{iu(M+J)}\,du \\
\nonumber&&\hspace{-0.7in}+
\int_0^t i^2 \int_{\GD\ub{2}} e^{i(1-u_1)(M+J)}d(M+J) e^{i(u_1-u_2)(M+J)}
d(M+J)e^{iu_2(M+J)}\,du   
\end{eqnarray}
Fixing $u_0=1$ and $u_{n+1}=0$ we have, by formula
(\ref{FF 85.01}),
\begin{eqnarray}
\nonumber dK\ub{n}&=&\sum_{j=1}^n\int_{\GD\ub{n}}e^{i(1-u_1)M}J\cdots
\!e^{i(u_{j-1}-u_j)M}dJ\cdots Je^{iu_n M}\,du
\\
\nonumber && \hspace{-0.9in}+
\sum_{1\leq j<k\leq n}\int_{\GD\ub{n}}e^{i(1-u_1)M}J\cdots
\!e^{i(u_{j-1}-u_j)M}dJ\cdots e^{i(u_{k-1}-u_k)M}dJ\cdots 
 Je^{iu_n M}\,du
\\
\nonumber &&\hspace{-0.5in}+ \sum_{j=1}^{n+1}\int_{\GD\ub{n}}
\!\!\!e^{i(1-u_1)M}J\cdots
d\left(e^{i(u_{j-1}-u_j)M}\right) J\cdots Je^{iu_n M}\,du
\\
\nonumber &&\hspace{-0.5in}+
\sum_{1\leq j<k\leq n}\int_{\GD\ub{n}}\!\!\!e^{i(1-u_1)M}J\cdots
e^{i(u_{j-1}-u_j)M}dJ\cdots d\left(e^{i(u_{k-1}-u_k)M}\right)J\cdots 
 Je^{iu_n M}\,du
\\
\nonumber &&\hspace{-0.5in} +
\sum_{1\leq j<k\leq n}\int_{\GD\ub{n}}\!\!\!e^{i(1-u_1)M}J\cdots
d\left(e^{i(u_{j-1}-u_j)M}\right)J\cdots e^{i(u_{k-1}-u_k)M}dJ\cdots 
Je^{iu_n M}\,du
\\
\nonumber &&\hspace{-0.5in}+
\sum_{1\leq j<k\leq n+1}\int_{\GD\ub{n}}\!\!\!e^{i(1-u_1)M}J\cdots
d\left(e^{i(u_{j-1}-u_j)M}\right)J\cdots d\left(e^{i(u_{k-1}-
u_k)M}\right)J\cdots Je^{iu_n M}\,du
\\
\nonumber&&\hspace{-0.5in}= \sum_{j=1}^n K\ub{n}(\stackrel{(j)}{dJ})
+\sum_{1\leq j<k\leq n}K\ub{n}(\stackrel{(j)}{dJ},\stackrel{(k)}{dJ})
+ \sum_{j=1}^{n+1} \left(iK\ub{n+1}(\stackrel{(j)}{dM})+i^2 
K\ub{n+2}(\stackrel{(j)}{dM},\stackrel{(j+1)}{dM})\right)
\\
\nonumber&&\hspace{-0.5in}+ \sum_{1\leq j<k\leq n+1} 
i K\ub{n+1}(\stackrel{(j)}{dJ},\stackrel{(k)}{dM})
+ \sum_{1\leq j<k\leq n+1} 
i K\ub{n+1}(\stackrel{(j)}{dM},\stackrel{(k)}{dJ})
\\
\nonumber &&\hspace{-0.5in} +
\sum_{1\leq j<k\leq n+1} 
i^2 K\ub{n+2}(\stackrel{(j)}{dM},\stackrel{(k+1)}{dM}).
\\
\label{EE 85.000} &&\hspace{-0.5in}= \sum_{j=1}^n K\ub{n}(\stackrel{(j)}{dJ})
+\sum_{1\leq j<k\leq n}K\ub{n}(\stackrel{(j)}{dJ},\stackrel{(k)}{dJ})
\\
\label{EE 85.001} &&\hspace{-0.5in}+ 
 \sum_{j=1}^{n+1} iK\ub{n+1}(\stackrel{(j)}{F\st})dA\d
+\sum_{1\leq j<k\leq n+1} 
i K\ub{n+1}(\stackrel{(j)}{dJ},\stackrel{(k)}{dM})
\\
\label{EE 85.002} &&\hspace{-0.5in}+ 
 \sum_{j=1}^{n+1} iK\ub{n+1}(\stackrel{(j)}{F})dA
+ \sum_{1\leq j<k\leq n+1} 
i K\ub{n+1}(\stackrel{(j)}{dM},\stackrel{(k)}{dJ})
\\
\label{EE 85.003}
&&\hspace{-0.5in}+ 
\sum_{1\leq j<k\leq n+2} 
i^2 K\ub{n+2}(\stackrel{(j)}{dM},\stackrel{(k)}{dM}).
\end{eqnarray}
The terms above may be calculated by using the multilinearity of the 
$K\ub{n}$ then taking out the basic differentials preserving 
their order and then using the quantum Ito table. For example
\begin{eqnarray*}
K\ub{n}(\stackrel{(j)}{dJ})&=&  
K\ub{n}(\stackrel{(j)}{S})\,dA 
+K\ub{n}(\stackrel{(j)}{S}{}\!\!\st)\,dA\d
+K\ub{n}(\stackrel{(j)}{R})\,dt,
\\
K\ub{n+1}(\stackrel{(j)}{dM},\stackrel{(k)}{dJ})&=&
K\ub{n+1}(\stackrel{(j)}{F},\stackrel{(k)}{S}{}\!\!\st)
\end{eqnarray*}
Now write $dK\ub{n}=F\ub{n} dA+G\ub{n}dA\d +H\ub{n}\,dt$. 
It follows from the bounds (\ref{EE 84.01}) that
\beqn\label{B 85.01}
\nm{F_t\ub{n}}\leq \frac{\nm{J_t}^{n-1}}{(n-1)!}\nm{S_t+F_t},
\quad
\nm{G_t\ub{n}}\leq \frac{\nm{J_t}^{n-1}}{(n-1)!}\nm{S_t\st+F_t\st}
\eeqn
\beqn \label{B 85.02}
\nm{H_t\ub{n}}\leq
\frac{\nm{J_t}^{n-1}}{(n-1)!}\nm{U_t}\!+\!\frac12 
\frac{\nm{J_t}^{n-2}}{(n-2)!} \nm{S_t}^2\!+\!
\frac{\nm{J_t}^{n-1}}{(n-1)!} \nm{F_t}\cdot\nm{S_t}
\!+\!\frac12 \frac{\nm{J_t}^{n}}{n!} \nm{F_t}^2.
\eeqn
It follows that $F\ub{n}, G\ub{n}\in L^2([0,1],\CB(\FH)\so)$ and
$H\ub{n}\in L^1([0,1],\CB(\FH)\so)$.
\bb\ni
We now expand the terms on the right hand side of (\ref{EE 
85.01}).
\begin{eqnarray}
&&\hspace{-0.5in}\nonumber i\int_0^1 e^{i(1-u)(M+J)}d(M+J)e^{iu(M+J)} \,du 
\\
&&\nonumber =i\int_0^1 
\left(\sum_{\Ga=0}\iiy i^{\Ga}(1-u)^{\Ga}
K\ub{\Ga}((1-u))\right) d(M+J)
\left(\sum_{\Gb=0}\iiy i^{\Gb}u^{\Gb}K\ub{\Gb}\right)\,du
\\
&& \nonumber =\int_0^1 
\sum_{n=0}\iiy i^{n+1} \sum_{\Ga+\Gb=n}(1-u)^{\Ga}u^{\Gb}
K\ub{\Ga}((1-u)) d(M+J)K\ub{\Gb}\,du 
\\
&&\nonumber= \sum_{n=0}\iiy\sum_{\Ga=0}^n i^{n+1} 
K\ub{n+1}(\stackrel{(\Ga+1)}{d(M+J)}) 
\\
&& \label{EE 85.02}= \sum_{n=1}\iiy\sum_{\Ga=1}^n
i^{n} K\ub{n}(\stackrel{(\Ga+1)}{dJ})
+
\sum_{n=0}\iiy\sum_{\Ga=0}^n 
i^{n+1} K\ub{n+1}(\stackrel{(\Ga+1)}{dM})
\end{eqnarray}
Similarly
\begin{eqnarray}
\lefteqn{i^2 \int_{\GD\ub{2}} e^{i(1-u_1)(M+J)}d(M+J) e^{i(u_1-u_2)(M+J)}
d(M+J)e^{iu_2(M+J)}\,du}
\\
&&\nonumber= \sum_{n=0}\iiy \sum_{\Ga+\Gb=0}\iiy 
i^{n+2}K\ub{n+2}(\stackrel{(\Ga+1)}{d(M+J)},\stackrel{(\Ga+\Gb+2)}{d(M+J)} 
\\
&&\label{EE 85.03} = \sum_{n=2}\iiy\sum_{1\leq j<k\leq 
n}\!\!\!\!\!\!i^n K\ub{n}(\stackrel{(j)}{dJ},\stackrel{(k)}{dJ})
+
\sum_{n=1}\iiy\sum_{1\leq j<k\leq n+1}\!\!\!\!\!\!i^{n+1}
K\ub{n+1}(\stackrel{(j)}{dJ},\stackrel{(k)}{dM}) 
\\
&&\nonumber = \sum_{n=1}\iiy\sum_{1\leq j<k\leq 
n+1}\!\!\!\!\!\!i^{n+1} K\ub{n+1}(\stackrel{(j)}{dM},\stackrel{(k)}{dJ})
+
\sum_{n=0}\iiy\sum_{1\leq j<k\leq n+2}\!\!\!\!\!\!i^{n+2}
K\ub{n+2}(\stackrel{(j)}{dM},\stackrel{(k)}{dM})
\end{eqnarray}
Adding (\ref{EE 85.02}) and (\ref{EE 85.03}) and comparing with 
the expression (\ref{EE 85.001}) for $dK\ub{n}$ gives
\begin{eqnarray*}
\sum_{n=0}\iiy i^n dK\ub{n}&=&
i\int_0^1 e^{i(1-u)(M+J)}d(M+J) e^{iu(M+J)}\,du \\ 
&& \hspace{-0.5in}+ i^2 \int_{\GD\ub{2}} e^{i(1-u_1)(M+J)}d(M+J) e^{i(u_1-u_2)(M+J)}
d(M+J)e^{iu_2(M+J)}\,du 
\end{eqnarray*}
It follows from the bounds (\ref{B 85.01}) and (\ref{B 85.02}) 
that if 
$$
\kappa(t)=\nm{J_t}e^{\nm{J_t}},\quad 
\Gb(t)=2e^{\nm{J_t}}\nm{S_t+F_t},
$$
$$
\Gg(t)= e^{\nm{J_t}}(\nm{U}_t+\frac12(\nm{S_t}+\nm{F_t})^2
$$
then the series $K\ub{n}$ and the quantum semimartingale $N$ defined 
by
\begin{eqnarray*}
N_t&=&I+\int_0^t i\int_0^1 e^{i(1-u)(M+J)}d(M+J) e^{iu(M+J)}\,du \\
&&\hspace{-0.5in}+ 
i^2 \int_0^t\int_{\GD\ub{2}} e^{i(1-u_1)(M+J)}d(M+J) e^{i(u_1-u_2)(M+J)}
d(M+J)e^{iu_2(M+J)}\,du 
\end{eqnarray*}
satisfy the conditions of Proposition \ref{PP 85.01}. The 
quantum Duhamel formula (\ref{EE 85.01}) now follows 
from Proposition \ref{PP 85.01} (b) and (\ref{EE 85.1}).
\bb\ni
Step 2. The integrands $E$, $H$, $S$ and $U$ are zero. In this 
case a product of differentials from $\{dM,dJ\}$ containing more 
than two $dM$s is zero as is any product containing the string 
$dJ\cdot dM\cdot dJ$.
\bb\ni
Expanding $dK\ub{n}$ as in Step 1 
and using the identity $d\GL\cdot d\GL=d\GL$ the sums 
(\ref{EE 85.000}) are replaced by
$$
\sum_{k=1}^n\sum_{1\leq j_1<\cdots<j_k\leq n} 
\!\!\!\!K\ub{n}(\stackrel{(j_1)}{dJ},\ldots,\stackrel{(j_k)}{dJ})
=
\sum_{k=1}^n\sum_{1\leq j_1<\cdots<j_k\leq n} 
\!\!\!\!K\ub{n}(\stackrel{(j_1)}{R},\ldots,\stackrel{(j_k)}{R})\,d\GL=E\ub{n}{}\pr 
\,d\GL.
$$
It follows from the bounds (\ref{EE 84.01}) that
$$
\nm{E\ub{n}{}\pr_t}\leq \sum_{r=1}^n \frac{\nm{J_t}^{n-r}\nm{R_t}^r}{(n-
r)!r!}. 
$$
It follows from (\ref{EE 84.06}) that
$$
\sum_{0}\iiy E\ub{n}{}\pr=e^{i(M+J+R)}-e^{i(M+J)} 
$$
the coefficient of $d\GL$ in the quantum Duhamel formula.
\bb\ni
Using the identity $d\GL\cdot dA\d=dA\d$ the sum
(\ref{EE 85.001}) is replaced by
\begin{eqnarray*}
\lefteqn{i\sum_{k=0}^n
\sum_{1\leq j_1<\cdots<j_{k+1}\leq n+1} 
K\ub{n+1}(\stackrel{(j_1)}{dJ},\ldots,
\stackrel{(j_k)}{dJ},\stackrel{(j_{k+1})}{dM})} \\
&=&
i\sum_{k=0}^n
\sum_{1\leq j_1<\cdots<j_{k+1}\leq n+1} 
K\ub{n+1}(\stackrel{(j_1)}{R},\ldots,
\stackrel{(j_k)}{R},\stackrel{(j_{k+1})}{F}\!\!{}\st)=iG\ub{n}{}\pr dA\d.
\end{eqnarray*}
Using the identity $dA\cdot d\GL=dA$ the sum 
(\ref{EE 85.002}) is replaced by
\begin{eqnarray*}
\lefteqn{i\sum_{k=0}^n
\sum_{1\leq j_1<\cdots<j_{k+1}\leq n+1} 
K\ub{n+1}(\stackrel{(j_1)}{dM},\ldots,
\stackrel{(j_k)}{dJ},\stackrel{(j_{k+1})}{dJ})} \\
&=&
i\sum_{k=0}^n
\sum_{1\leq j_1<\cdots<j_{k+1}\leq n+1} 
K\ub{n+1}(\stackrel{(j_1)}{F},\ldots,
\stackrel{(j_k)}{R},\stackrel{(j_{k+1})}{R})=iF\ub{n}{}\pr dA.
\end{eqnarray*}
Using the identity $dA\cdot d\GL\cdots d\GL dA\d=dt$ the sum 
(\ref{EE 85.003}) is replaced by
\begin{eqnarray*}
\lefteqn{i^2  \sum_{k=1}^n K\ub{n+2}(\stackrel{(j_1)}{dM},\stackrel{(j_2)}{dJ},\ldots,
\stackrel{(j_k)}{dJ},\stackrel{(j_{k+1})}{dM})} \\
&=&
i^2  \sum_{k=1}^n 
K\ub{n+2}(\stackrel{(j_1)}{F},\stackrel{(j_2)}{R},\ldots,
\stackrel{(j_k)}{R},\stackrel{(j_{k+1})}{F}{}\!\!\st) =i^2 
H\ub{n}{}\pr dt
\end{eqnarray*}
It follows from the bounds (\ref{EE 84.01}) that
\begin{eqnarray*}
\nm{F\ub{n}_t{}\pr},\nm{G\ub{n}_t{}\pr} 
\leq \sum_{r=1}^n \frac{\nm{J_t}^{n-r}\nm{R_t}^r}{(n-r)!r!}\nm{F_t},
\quad  \nm{H\ub{n}_t{}\pr}\leq 
\frac12\sum_{r=1}^n \frac{\nm{J_t}^{n-r}\nm{R_t}^r}{(n-r)!r!}\nm{F_t}^2. 
\end{eqnarray*} 
The coefficient of $dA\d$ in the quantum Duhamel formula is 
\begin{eqnarray*}
i\int_0^1 e^{i(M+J+R)} F\st e^{i(M+J)}\,du &&
\\
&& \hspace{-2in}
=i\int_0^1 \sum_{\Ga=0}\iiy \sum_{k=0}^{\Ga} \sum_{1\leq 
j_1<\cdots< j_k\leq \Ga}\!\!\!\!\!\!i^{\Ga}(1-u)^{\Ga}K\ub{\Ga}((1-u);
\stackrel{(j_k)}{R},\stackrel{(j_{k+1})}{R})F\st
i^{\Gb} u^{\Gb} K\ub{u} \,du
\\
&& \hspace{-2in}
= \sum_{k=0}^n
\sum_{1\leq j_1<\cdots<j_{k+1}\leq n+1} 
\!\!\!\!\!\!\!i^{n+1}K\ub{n+1}(\stackrel{(j_1)}{R},\ldots,
\stackrel{(j_k)}{R},\stackrel{(j_{k+1})}{F}\!\!{}\st)
\end{eqnarray*}
Similarly the coefficient of $dA$ is
\begin{eqnarray*}
i\int_0^1 e^{i(M+J)} F e^{i(M+J+R)}\,du
&=&
i\sum_{k=0}^n
\sum_{1\leq j_1<\cdots<j_{k+1}\leq n+1} 
\!\!\!\!\!\!\!i^{n+1}K\ub{n+1}(\stackrel{(j_1)}{F},\ldots,
\stackrel{(j_k)}{R},\stackrel{(j_{k+1})}{R}).
\end{eqnarray*}
The coefficient of $dt$ is 
\begin{eqnarray*}
\lefteqn{i^2\int_{\Gd\ub{2}}\sum_{\Ga=0}\iiy K\ub{\Ga} F
\sum_{\Gb=0}\iiy \sum_{k=0}^{\Gb}
\sum_{1\leq j_1<\cdots<j_k\leq \Gb}
\!\!\!\!\!\!\!K\ub{\Gb}(u_1\!-\!u_2);\stackrel{(j_1)}{R},\ldots,\stackrel{(j_k)}{R})
F\st
\sum_{\Gg=0}\iiy K\ub{\Gg}(u_2)\,du} 
\\
&=& \sum_{n=0}\iiy i^{n+2} \sum_{\Ga+\Gb+\Gg=n} 
\sum_{\Ga+1\leq j_1<\cdots<j_k \leq\Ga+\Gb+1}
\!\!\!\!\!\!K\ub{\Ga+\Gb+\Gg+2}
(\stackrel{(\Ga)}{F},\stackrel{(j_1)}{R},\ldots,\stackrel{(j_k)}{R},
\stackrel{(\Ga+\Gb+2)}{F}\!\!\!\st),
\\
&=&
\sum_{n=0}\iiy \sum_{k=1}^n 
i^{n+2}K\ub{n+2}(\stackrel{(j_1)}{F},\stackrel{(j_2)}{R},\ldots,
\stackrel{(j_k)}{R},\stackrel{(j_{k+1})}{F}{}\!\!\st).
\end{eqnarray*}
\bb\ni
The quantum Duhamel formula now follows using the same argument 
as in Step 2.
\bb\ni
This method may also be used when $E$ is non-zero but the calculations 
are very tedious and we leave them to the conscientious reader. \,\,\epf

\section{\label{PROB}Further Problems}
We add five more problems to those posed in 
\cite[\S 12]{Vin2}.
\bb\ni
This paper makes some progress with the problems in 
\cite[\S 12]{Vin2} but provides nothing like a complete solution. 
The most important question is the following.
\bb\ni
{\it Problem 1.} {\em If $M$ is an essentially self-adjoint 
quantum semimartingale does $M$ satisfy the quantum Duhamel 
formula.}
\bb\ni
It is clear that both sides of (\ref{EE 72.01}) exist. Are they 
equal ?
\bb\ni
Theorem \ref{T 81.361} is not a generalisation of 
\cite[Theorem 6.2]{Vin2} since not all regular quantum 
semimartingales satisfy its conditions. The conditions in 
Theorem \ref{T 81.361} are on $(E,F,G,H)$ while those in 
\cite[Theorem 6.2]{Vin2} are on $(E,F,G,H)$ {\em and} $M$. 
Moreover, a self-adjoint operator $T\in \CB(\FH)$ need not have 
$\ell_{00}\in \CA(\Lnu(T))$.
\bb\ni
It is possible to construct a
synthetic theorem containing Theorem \ref{T 81.361} and 
Theorem 6.2 of \cite{Vin2}, but this would be unsatisfactory.
\bb\ni
{\it Problem 3.} {\em Find natural conditions on a quantum 
semimartingale $M$ for which the conclusions of both
Theorem {\rm\ref{T 81.361}} and {\rm \cite[Theorem 6.2.]{Vin2}}
remain true.}
\bb\ni
In general the cmx representation of the quantum stochastic 
process $F$ considered in Section \ref{PCM} is complicated and does not 
satisfy the conditions of 
Corollary \ref{T 81.362}. In that case the classical Ito formula 
(\ref{E 80.080}) does not follow from (\ref{Ito}). If $M$ is 
irregular the classical Ito formula for $M$ does not follow from
\cite[Theorem 6.2]{Vin2}.
\bb\ni
{\it Problem 3.} {\em Find conditions on an
integrable quadruple $(0,F,F\st,0)$, which is not necessarily 
commutative, which ensure that the 
quantum Ito formula {\rm (\ref{Ito})} is valid and implies the 
classical Ito formula {\rm (\ref{E 80.080})}.}
\bb\ni
Formally, the analysis in Section \ref{PCM} carries over 
to the 
case when $M$ is an essentially self-adjoint quantum 
semimartingale on $\CE$ which satisfies the quantum Ito formula.
In this case $h\pr$ and $h\prpr$ must be replaced by $Dh$ and 
$D^2_I h$. There are serious technical difficulties in finding 
bounds on $Dh$ and $D^2_I h$ which we have been unable to 
overcome.  
See for example \cite{PEL}.
\bb\ni
{\it Problem 4.} {\em Is it possible to modify the method used 
in Section \ref{PCM} to give a more analytic proof of Theorem
\ref{TT 85.01}\,?}
\bb\ni
The final problem is about the Stratonovich formula.
\bb\ni
{\it Problem 5. Is the cmx Stratonovich formula (\ref{E 
86.767}) valid for $f\in C^{3+}(\BBR)$ and does it 
give rise to the corresponding quantum Stratonovich formula 
(\ref{E 83.38}) ?}

\newpage
\btb{99}

\bibitem{Att} S. ATTAL, `An algebra of non-commutative bounded semimartingales: 
square and angle quantum brackets', Journal of Functional Analysis, 124 (1994) 
292--332.
\bibitem{Att1} S. ATTAL, `Classical and quantum stochastic 
calculus', {\em Prepublication de l'Institute Fourier, Grenoble, 
358 (1998)} 
\bibitem{AM} S. ATTAL and P. A. MEYER, `Interpr$\acute{\mbox{e}}$tation 
probabiliste et extension 
des int$\acute{\mbox{e}}$grales stochastiques non commutatives', 
Springer Lecture Notes  in 
Mathematics, 1557, (1993) 312--327. 
\bibitem{BR1} O. BRATTELI and D. W. ROBINSON, {\em Operator algebras and quantum 
statistical mechanics I} (Springer--Verlag, 1979).
\bibitem{BR2} O. BRATTELI and D. W. ROBINSON, {\em Operator algebras and quantum 
statistical mechanics II} (Springer--Verlag, 1979).
\bibitem{DiU} J. DIESTEL and J. J. UHL, {\em Vector measures} (American 
Mathematical Society, 1977).
\bibitem{Ev1} M. P. EVANS,  `Existence of quantum diffusions', 
Probab. Theory Related Fields 81 (1989) 473--483.
\bibitem{GC} C. W. GARDINER and M. J. COLLETT, Input and output 
in damped quantum systems: `Quantum stochastic differential 
equations and the master equation', Phys. Rev. A 31 (6) (1985) 
3761--3773.  
\bibitem{GJ} J. GLIMM and A JAFFE, {\em Quantum Physics, A 
Functional Integral Point of View,} (Springer-Verlag, New York 
1981).
\bibitem{HiP} E. HILLE and R. S. PHILLIPS, {\em Functional 
Analysis and Semigroups,} (Colloq. Publ. Amer. Math. Soc. 1957). 
\bibitem{HP1} R. L. HUDSON and K. R. PARTHASARATHY,  `Quantum Ito's 
formula and stochastic evolutions', Comm. Math. Phys. 93, (1984) 301--323.
\bibitem{KS} I. KARATZAS and S. E. SHREVE, {\em Brownian motion 
and stochastic calculus} (Springer-Verlag 1988).
\bibitem{KuW} H. KUNITA and S. WATANABE, `On square integrable 
martingales', Nagoya Math. J. 30 (1967) 209--245.
\bibitem{Lin1} J. M. LINDSAY, `Quantum and non-causal stochastic 
calculus', Probab. Theory Related Fields 97 (1993) 65--80.
\bibitem{Lin2} J. M. LINDSAY, `Independence for quantum 
stochastic integrators', {\em Quantum Probability and Related 
Topics, Vol. VI,} (ed. L. Accardi) (World Scientific, Singapore, 
1991) 325--332.
\bibitem{MEY} P. A. MEYER, {\em Quantum probability for probabilists} 
(Lecture Notes in Mathematics 1538, Springer-Verlag 1993).
\bibitem{Nel} E. NELSON, `Analytic vectors', Ann. Math. 70 (1959) 
572--615.
\bibitem{OB1} N. OBATA, {\em White noise calculus and Fock space} 
(Lecture Notes in Mathematics 1577, Springer-Verlag 1993).
\bibitem{OB1} N. OBATA, `Generalised quantum stochastic processes 
on Fock space, Publ. RIMS, Kyoto Univ. 31 (1995) 667--702.
\bibitem{P1} K. R. PARTHASARATHY, {\em An introduction to quantum stochastic 
calculus} (Birkhauser, 1992).
\bibitem{P2} K. R. PARTHASARATHY, `Quantum stochastic calculus', 
in {\em Proceedings of Stochastic processes and their 
applications, Nagoya 1985,} (ed. K Ito and T Hida) Springer Lecture Notes in 
Mathematics 1203, (1986) 177--196.
\bibitem{PEL} V. V. Peller, Lecture Notes in Pure and Applied 
Mathematics 122, (1990).
\bibitem{RS1} M. REED and B. SIMON, {\em Methods of Modern 
Mathematical Physics, I and II} (Academic Press, 1972 and 1975).
\bibitem{RoW} L. C. G. ROGERS and D. WILLIAMS, {\em Diffusions, Markov 
processes, and martingales: Ito calculus} (John Wiley and Sons, 1987).
\bibitem{Rud} W. RUDIN, {\em Real and Complex Analysis} 
(McGraw-Hill 1970).
\bibitem{SA} T. SAITO and T. ARIMITSU, `Quantum stochastic 
Liouville equation of Ito type', Modern Phys. Lett. B, 7 (29-30) 
(1993) 1951--1959.
\bibitem{Vin} G. F. VINCENT--SMITH, `Unitary quantum stochastic evolutions', 
Proc. London Math. Soc. (3) 63 (1991) 401--425. 
\bibitem{Vin2} G. F. VINCENT--SMITH, `The Ito formula for quantum 
semimartingales', Proc. London Math. Soc. (3) 75 (1997) 671--720.
\bibitem{Vin3} G. F. VINCENT--SMITH, `The Ito formula for 
perturbed Brownian martingales', Bull. London Math. Soc. 32 (2000) 736--745.
\bibitem{Wei} J. WEIDMANN, {\em Linear Operators in Hilbert 
Spaces} (Springer-Verlag, 1980).

\etb          
\vspace{1in}
ORIEL COLLEGE
\bb\ni
OXFORD OX1 4EW
\bb\ni
UNITED KINGDOM
\bb\ni
graham.vincent-smith@oriel.ox.ac.uk

\end{document}